\documentclass[letterpaper]{article}

\usepackage{arxiv}

\usepackage{amsmath, amssymb, latexsym, multicol, amsthm, bm, placeins}
\usepackage{color}
\usepackage{cleveref}

\usepackage{style}

\makeatletter\@addtoreset{equation}{section}\makeatother
\newcommand{\Hr}{{\rm H}}
\newcommand{\Wr}{{\rm W}}
\newcommand{\Lr}{{\rm L}}

\title{Meshfree Variational Physics Informed Neural Networks (MF-VPINN): an adaptive training strategy}

\author{Stefano Berrone\thanks{Dipartimento di Scienze Matematiche, Politecnico di Torino, Corso Duca degli Abruzzi 24, 10129 Torino, Italy. stefano.berrone@polito.it}
	\And Moreno Pintore\thanks{MEGAVOLT team, Inria -- Laboratoire Jacques-Louis Lions, SCAI, Sorbonne Universit\'e, 4 place Jussieu, 75005 Paris, France. moreno.pintore@inria.fr}
}

\renewcommand{\shorttitle}{Meshfree Variational Physics Informed Neural Networks (MF-VPINN)}

\begin{document}

\maketitle

\begin{abstract}
In this paper, we introduce a Meshfree Variational-Physics-Informed Neural Network. It is a Variational-Physics-Informed Neural Network that does not require the generation of the triangulation of the entire domain and that can be trained with an adaptive set of test functions. In order to generate the test space, we exploit an a posteriori error indicator and add test functions only where the error is higher. Four training strategies are proposed and compared. Numerical results show that the accuracy is higher than the one of a Variational-Physics-Informed Neural Network trained with the same number of test functions but defined on a quasi-uniform mesh.
\end{abstract}

\vspace{0.3cm}
\noindent \textbf{\textit{MSC-class}} 65N12; 65N15; 65N50; 68T05; 92B20

\section{Introduction}\label{sec:introduction}
Physics-Informed Neural Networks (PINNs) are a rapidly emerging numerical technique used to solve partial Differential equations (PDEs) by means of a deep neural network. The first idea can be traced back to the works of Lagaris et al. \cite{lagaris1997artificial, lagaris1998artificial, lagaris2000neural}, but, thanks to the hardware advancements and the existence of deep learning packages like Tensorflow \cite{tensorflow2015-whitepaper}, Pytorch \cite{NEURIPS2019_9015} {{and JAX} 
 \cite{jax2018github}}, they have recently became popular since the works of Raissi et al. \cite{raissi2017physics, raissi2017bphysics}, published in \cite{raissi2019physics}. In its original formulation, the approximate solution is computed as the output of a neural network trained to minimize the PDE residual on a set of collocation points inside the domain and on its boundary.

The growing interest in PINNs is strictly related to their flexibility.  In fact, with minor changes to the implementation, it is possible to solve a huge variety of problems. For example, exploiting the nonlinear nature of the involved neural network, \mbox{nonlinear \cite{pu2021solving, yuan2022pinn}} and high-dimensional PDEs \cite{guo2023high} can be solved without the need for globalization methods or additional nonlinear solvers. Moreover,  by changing the neural network's input dimensions or suitably adapting the loss function, it is possible to solve parametric \cite{demo2021extended, gao2021phygeonet} or \mbox{inverse \cite{chen2020electromagnetic, tartakovsky2018learning}} problems. When external data are available, they can also be used to guide the optimization phase and improve the PINN accuracy \cite{zhao2021physics}.

In order to improve the original PINN proposed in \cite{raissi2019physics} and to adapt it to solve specific problems, several generalizations have been proposed. For example, the deep Ritz method {(DRM) \cite{yu2018deep, muller2022error, lu2021priori} looks for a minimizer of the PDE energy functional and, in the deep Galerkin method (DGM) \cite{sirignano2018dgm, al2022extensions, li2021deep}, an approximation of the $L^2$ norm of the PDE residual is minimized. It is also possible to exploit domain decomposition strategies \cite{smith1997domain, toselli2006domain} as in the conservative PINN (CPINN) \cite{ameya2020conservative}, in the parallel PINN \cite{shukla2021parallel}, in the extended PINN (XPINN) \cite{jagtap2020extended}, or in the Finite Basis PINN (FBPINN) \cite{moseley2023finite}. Moreover, it is even possible to change the neural network architecture or the training strategy as in \cite{gao2021phygeonet,viana2021estimating,yang2021bpinns,yang2020physics,yucesan2021hybrid,zhu2019physics, pang2019fpinns}; between the methods based on different architectures, we highlight some works based on the novel Kolmogorov--Arnold Network (KAN) \cite{liu2024kan} architecture \cite{koenig2024kan, qian2024investigating} and on a Large Language Model (LLM) \cite{kumar2023mycrunchgpt}.  More extensive overviews of the existing approaches can be found in \cite{beck2022overview, cuomo2022scientific,lawal2022physics, viana2021survey}. In the context of the current work, an important extension is the Variational-Physics-Informed Neural Network (VPINN) \cite{kharazmi2019variational, kharazmi2021hp}, where the weak formulation of the problem is used to construct the loss function. }

In this work, we focus on VPINNs. As discussed in \cite{berrone2022solving,berrone2022variational,kharazmi2019variational, kharazmi2021hp}, in order to train a VPINN, one needs to choose a suitable space of test functions, compute the variational residuals against all the test functions on the basis of such a space, and minimize a linear combination of these residuals. {Since a spatial mesh is required to define the test functions, the VPINN cannot be considered a meshfree method, even though it is an extension of the PINN, which is meshfree.} In this work, we present an adaptive Meshfree VPINN (MF-VPINN) that does not require a global triangulation of the domain but is trained with the same loss function and neural network architecture of a standard VPINN. {Note that the MF-VPINN and the original VPINN can solve the same differential problems because the neural network is trained with the same loss functions. We also highlight that they can solve problems where the solution has low regularity that cannot be solved with standard PINNs, for example, in the presence of singular forcing terms, thanks to the weak formulation of the PDE without introducing further approximations or regularizations.} {However, one of the VPINN's limitations is that a triangulation of the entire domain is required to define the test functions. Generating it may be very expensive or even impractical for very complex geometries (like, for example, the ones in \cite{berrone2016towards}) and in moderate- or high-dimensional problems, for which automatic mesh-generation algorithms do  not exist. For such domains, it is therefore highly advisable or computationally necessary to use a meshfree method such as the original PINN or the proposed MF-VPINN. Moreover, when dealing with complex geometries for which a mesh can be hardly generated, the refinement of the mesh for adaptive methods can be very difficult. In this paper, we describe an algorithm that solves the problem and provides a reliable solution.}

The paper is organized as follows. In Section \ref{sec:problem_formulation}, we introduce the problem we are interested in. In particular, we focus on the problem discretization in Section \ref{sec:problem_discretization} and on the MF-VPINN loss function in Section \ref{sec:loss_function}. Then, an a posteriori error estimator is presented in Section \ref{sec:aposteriori} and used in Section \ref{sec:choice_of_M} to iteratively generate the required test functions. Numerical results are presented in Section \ref{sec:numerical_results}. In Section \ref{sec:implementation}, we describe the model implementation and some strategies to improve the model efficiency, in \mbox{Section \ref{sec:performance}} we compare different approaches to generate the test functions and compare their performance and, in Section \ref{sec:estimator_role}, we analyze the role of the error estimator introduced in Section \ref{sec:aposteriori}. {Similar tests are performed on a different problem in Section \ref{sec:holes} to describe possible extensions on more complex domains.} Finally, we conclude the paper in Section \ref{sec:conclusion} and discuss future perspectives and ideas.

\section{Problem Formulation}\label{sec:problem_formulation}

 Let us consider the following second-order elliptic problem, defined on a polygonal or polyhedral domain $\Omega\subset\R^n$ with a Lipshitz boundary $\Gamma=\partial\Omega$:
\begin{equation}\label{eq:model-pb}
\begin{cases}
Lu:=-\nabla \cdot (\mu \nabla u) + \boldsymbol{\beta}\cdot \nabla u + \sigma u =f & \text{in \ } \Omega\,, \\
u=g & \text{on \ } \Gamma \,, \end{cases}
\end{equation}
where $\mu, \sigma \in \Lr^\infty(\Omega)$, {$ \boldsymbol{\beta}$}$~\in (\Wr^{1,\infty}(\Omega))^n$ satisfy $\mu \geq \mu_0$, {$\sigma - \frac12 \nabla \cdot \boldsymbol{\beta} \geq 0$} 
 in $\Omega$ for some constant $\mu_0>0$, whereas $f \in L^2(\Omega)$ and $g=\bar{u}_{\vert\Gamma}$ for some $\bar{u} \in \Hr^1(\Omega)$.

In order to derive the corresponding variational formulation, we define the bilinear form $a$ and the linear form $F$ as
\begin{equation}\label{eq:form a}
a:V\times V \to \mathbb{R}\,, \qquad a(w,v)=\int_\Omega \mu \nabla w \cdot \nabla v + \boldsymbol{\beta}\cdot \nabla w \, v + \sigma w \, v\,,
\end{equation}
\begin{equation}\label{eq:forms F}
F:V\to \mathbb{R}\,, \qquad F(v)=\int_\Omega f \, v  \,;
\end{equation}
where $V$ is the function space $V=H_0^1(\Omega)$. We denote by $\alpha\ge\mu_0$ the coercivity constant of $a$ and by $\Vert a\Vert$ and $\Vert F\Vert$ the continuity constants of $a$ and $F$. Then, the variational formulation of Problem \eqref{eq:model-pb} reads as follows:  {Find $u \in\bar{u}+ V $ such that}
\begin{equation}\label{eq:model-pb-var}
a(u,v)=F(v) \qquad \forall v \in V\,.
\end{equation}

\subsection{Problem Discretization}\label{sec:problem_discretization}
{In order to numerically solve Problem \eqref{eq:model-pb-var}, one needs to choose suitable finite-dimensional approximations of the trial space} $\bar{u}+V$ and of the test space $V$. A Galerkin formulation is considered when we consider a finite-dimensional space $V_h^{\rm{trial}}$ for the trial space $\bar{u} + V_h^{\rm{trial}}$ and a finite-dimensional test space $V_h^{\rm{test}}$, with $V_h^{\rm{trial}}=V_h^{\rm{test}}$; whereas a Petrov--Galerkin formulation is considered otherwise. In this work, we consider a Petrov--Galerkin formulation in which the trial space is approximated by a set of functions $V^\NN$ of the form $V^\NN=\bar{u}+V_h^{\rm{trial}}$, with $V_h^{\rm{trial}}$ represented by a neural network suitably modified to enforce the Dirichlet boundary conditions, and the test space is a space $V_h$ of piecewise \mbox{linear functions}.

The neural network considered in the following is a standard fully connected feed-forward neural network. Given the number $L$ of layers and a set of matrices $A_\ell\in\mathbb{R}^{N_\ell\times N_{\ell-1}}$ and vectors $b_\ell\in\mathbb{R}^{N_\ell}$, $\ell=1,\ldots,L$ containing the neural network's trainable weights, the function $w:\R^n\rightarrow\R$ associated with the considered neural network architecture is:  
\begin{equation} \label{eq:nn_formula}
  \begin{aligned}
  &x_0=\boldsymbol{x}, \\
 &x_\ell = \rho(A_\ell x_{\ell-1} + b_\ell), \hspace{1.cm} \ell = 1,\ldots,L-1, \\
 &w(\boldsymbol{x}) = A_{L} x_{L-1} + b_L.
  \end{aligned}
\end{equation} 
where $\rho:\R\rightarrow\R$ is a nonlinear function applied element-wise  to the vector \mbox{$A_\ell x_{\ell-1} + b_\ell$}. In this section, we use $\rho(x)=\tanh(x)$; other common choices include, but are not limited to, $\rho(x)=\text{ReLU}(x)=\max\{0,x\}$, $\rho(x)=\text{RePU}(x)=\max\{0,x^p\}$ for $1<p\in\N$, \mbox{$\rho(x)=1/(1+e^{-x})$} and $\rho(x)=\log(1+e^x)$. Note that, in order to represent a function $w:\R^n\rightarrow\R$, the layer widths $N_{\ell}$ of the first and last layers are chosen as $N_0=n$ and $N_L=1$. We denote by $W^\NN$ the set of functions that can be represented as in \eqref{eq:nn_formula} for any combination of the neural network weights and by ${\bm w}^\NN$ the vector containing all the trainable weights of the neural network.

The function $w$ defined in \eqref{eq:nn_formula} is independent of the differential problem that has to be solved and is, in most   papers on PINNs or related models, trained to minimize both the residual of the equation and a term penalizing the discrepancy between $w_{|\Gamma}$ and $g$. Instead, we add a non-trainable layer $B$ to the neural network architecture in order to automatically enforce the required boundary conditions without the need to learn them during the training. As described in \cite{sukumar2022exact}, the operator $B$ acts on the neural network \mbox{output as}
\begin{equation}\label{eq:B_definition}
Bw = \phi w + \bar{g},
\end{equation}
where $\phi:\Omega\rightarrow\R$ is a function vanishing on $\Gamma$ and strictly positive inside $\Omega$, and $\bar{g}:\Omega\rightarrow\R$ is a suitable extension of $g:\Gamma\rightarrow\R$. The advantages of such an approach are also described in \cite{berrone2022enforcing}. Then, the discrete trial space approximating $\bar{u}+V$ can be defined as
$$
V^\NN = \{ v^\NN \in\bar{u}+ V : v^\NN=Bw \text{ for some }w\in W^\NN \}\,.
$$

On the other hand, the discrete test space $V_h$ is not associated with the neural network and only contains known test functions. In standard VPINNs, one generates a triangulation $\cal{T}$ of the domain $\Omega$ and then defines $V_h$ as the space of functions that coincide with a polynomial of order $p\in\N$ inside each element of $\cal{T}$. Instead, we want to construct a discrete space $V_h$ of functions independent from a global triangulation $\cal{T}$. Moreover, since in \cite{berrone2022variational} it has been proven that the VPINN convergence rate with respect to mesh refinement decreases when the order of the test functions is increased, we are interested in a space $V_h$ that only contains piecewise linear functions. For the sake of simplicity, we only consider the case where $n=2$; the discussion can be directly generalized to the more general case $n\in\N$.

Let $\hat P\subset\R^n$ be a reference patch. In the following discussion, $\hat P$ can be any arbitrary star-shaped polygon with $N_{\hat P}$ vertices and the dimension of its kernel strictly greater than zero. Nevertheless, in the numerical experiments, we only consider the reference patch $\hat P=[0,1]^2$ to avoid any unnecessary computational overhead. Let ${\cal{M}}=\{M_i\}_{i=1}^{n_\text{patches}}$ be a set of affine mappings such that $M_i:\hat P\rightarrow P_i\subset \Omega$, where we denote as $P_i$ the patch obtained transforming the reference patch $\hat P$ through the map $M_i$. We assume that ${\cal P}=\{P_i\}_{i=1}^{n_\text{patches}}$ is a cover of $\Omega$, i.e., $\cup_{i=1}^{n_\text{patches}} P_i = \Omega$, and we admit overlapping patches.

Let us consider the triangulation $\hat{\cal{T}}=\{\hat T_j:1\le j\le N_{\hat P}\}$ of $\hat P$ obtained by connecting each vertex with a single point $\bm{c}_{\hat P}$ in its kernel. It is then possible to define a piecewise linear function $\hat \varphi$ vanishing on the border of $\hat P$   such that $\hat\varphi(\bm{c}_{\hat P})=1$ and $\hat\varphi_{|\hat T_j}\in \P_1(\hat T_j)$, for any $j=1,\dots,N_{\hat P}$. Then, we define the discrete test space $V_h$ as $V_h = \text{span}\{\varphi_i:i=1,\dots,n_\text{patches}\}$, where $\varphi_i\in V$ is the piecewise linear function:
\begin{equation}\label{eq:phi_def}
 \varphi_i(\bm x) = \begin{cases} 
      \hat\varphi(M_i^{-1}(\bm x)), & \bm x\in P_i, \\
      0, & \bm x\notin P_i.
   \end{cases}
\end{equation}

\begin{figure}[t!]
\centering 
  \includegraphics[width=0.5\linewidth]{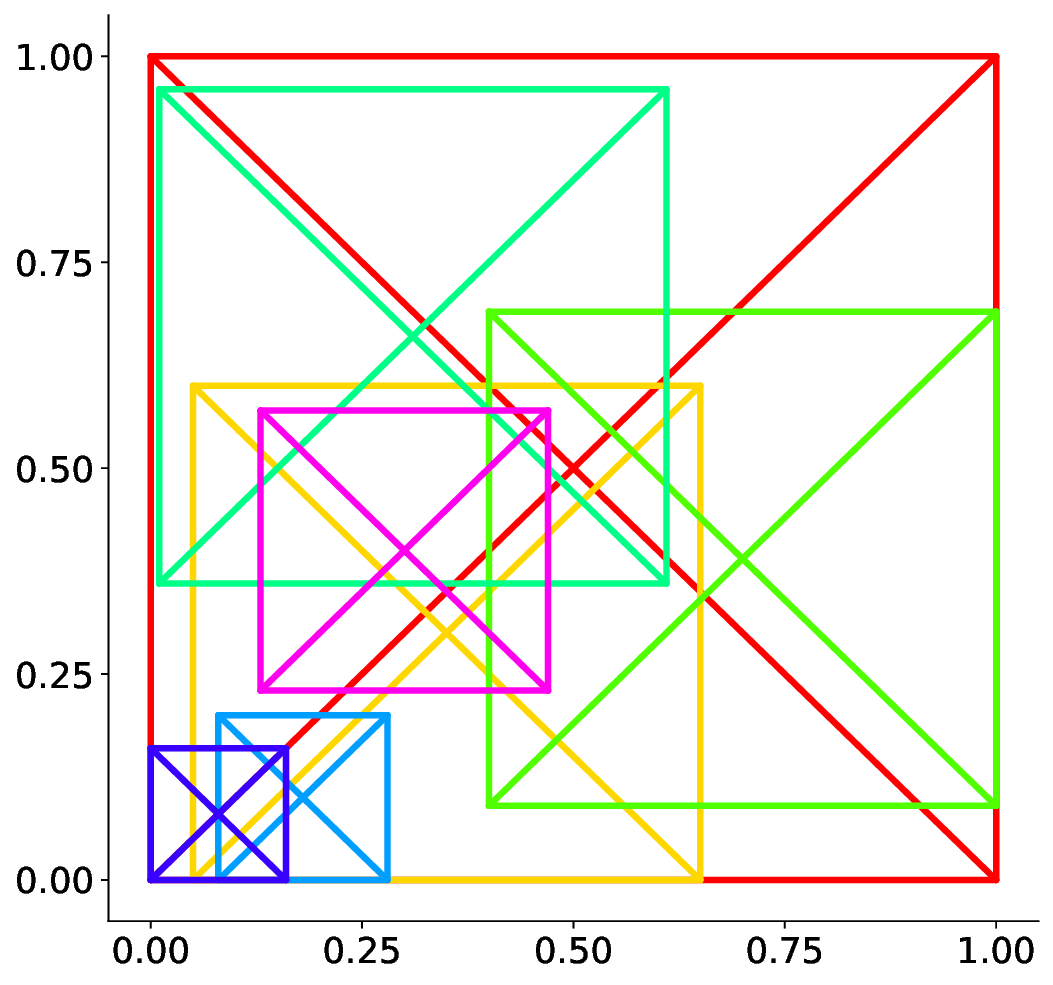} 
  \caption{Graphical representation of a set $\{P_i\}_{i=1}^{n_\text{patches}}$ obtained from a squared reference patch $\hat P$ with {$\bm{c_{\hat P}}$} in its center covering the domain $\Omega=(0,1)^2$.}
  \label{fig:example_of_patches}
\end{figure}

We remark that the only required triangulation is $\hat{\cal{T}}$, which contains only $N_{\hat P}$ triangles (in the numerical tests in this paper, $N_{\hat P}=4$). Instead, there exists no mesh on $\Omega$ and the test functions $\varphi_i$ and their supports $P_i$ are all independent. Therefore, the proposed method is said to be meshfree. A simple example of a set of patches ${\cal P}$ with $n_\text{patches}=7$ on the domain $\Omega=[0,1]^2$ is shown in Figure \ref{fig:example_of_patches}. For the sake of simplicity, in this work, we consider a squared reference patch $\hat P$ with {$\bm{c_{\hat P}}$} 
 coinciding with its center,  and we  let each mapping $M_i$ represent a combination of scalings and translations.

Using the introduced finite-dimensional set of functions $V^\NN$ and $V_h$, it is possible to discretize Problem \eqref{eq:model-pb-var} as follows: {Find $u^\NN \in V^\NN $ such that}
\begin{equation}\label{eq:discrete_model-pb-var}
a(u^\NN,v)=F(v) \qquad \forall v \in V_h\,.
\end{equation}

\subsection{Loss Function}\label{sec:loss_function}
In this section, we derive the loss function used to train the neural network. It has to be computable, and its minimizer has to be an approximate solution of {Problem \eqref{eq:model-pb}. We highlight that, when a standard PINN is used, the loss function can be seen as a discrete cost penalizing the residual of  \eqref{eq:model-pb} directly. In this context, instead, the loss function penalizes the variational residuals of \eqref{eq:model-pb-var} as in standard VPINNs. This is the key difference that differentiates the VPINNs (and its extension proposed in this manuscript) from the other generalizations of the original PINN introduced in Section \ref{sec:introduction}.}

Let us consider a quadrature rule of order $q\ge2$ on each triangle $T_j\in\hat{\cal{T}}$, $j=1,\dots,N_{\hat P}$, uniquely identified by a set of nodes and weights $\{(\widetilde{\bm{\xi}}^j_\ell,\widetilde\omega^j_\ell) : \ell \in I^{T_j}\}$. The nodes and weights of a composite quadrature formula of order $q$ on $\hat P$ can be obtained as 
\begin{equation*}\label{eq:ref_patch_quad_rule}
\{(\hat{\bm{\xi}}_\ell,\hat\omega_\ell) : \ell \in I^{\hat P}\} = 
\bigcup_{j=1}^{N_{\hat P}}\{(\widetilde{\bm{\xi}}^j_\ell,\widetilde\omega^j_\ell) : \ell \in I^{T_j}\}.
\end{equation*}

{Then, the corresponding quadrature rule of order $q$ of an arbitrary patch $P_i$ is \mbox{defined as}}
\begin{equation}\label{eq:patch_quad_rule}
\left\{(\bm{\xi}^i_\ell,\omega^i_\ell) : \ell \in I^{\hat P} | \bm{\xi}^i_\ell=M_i(\hat{\bm{\xi}}_\ell), \omega^i_\ell = \hat\omega_\ell\dfrac{\text{area}(P_i)}{\text{area}(\hat P)}\right\}.
\end{equation}

Using the quadrature rule in \eqref{eq:patch_quad_rule}, it is possible to define an approximate restriction on each patch of the forms $a$ and $F$ as follows:
\begin{equation}\label{eq:def-ah}
a_h^i(w,v)= \sum_{\ell \in I^{\hat P}} [\mu \nabla w \cdot \nabla v + \boldsymbol{\beta}\cdot \nabla w \, v + \sigma w v](\bm{\xi}^i_\ell) \,\omega^i_\ell \approx a_{P_i}(w,v)\,, 
\end{equation}
\begin{equation}\label{eq:def-Fh}
F_h^i(v) =  \sum_{\ell \in I^{\hat P}} [ f v](\bm{\xi}^i_\ell) \,\omega^i_\ell \approx F_{P_i}(v)\,,
\end{equation}
where $a_{P_i}(w,v)$ and $F_{P_i}(v)$ are defined as in \eqref{eq:form a} and \eqref{eq:forms F} but restricting the supports of the integrals to $P_i$. We remark that, since it is not possible to compute integrals involving a neural network exactly, we can only use the forms $a_h^i$ and $F_h^i$ in the loss function. Exploiting the linearity of $a(w,v)$ and $F(v)$ with respect to $v$ to consider only the basis $\{\varphi_i\}_{i=1}^{n_\text{patches}}$ of $V_h$  as set of test functions, we approximate Problem \eqref{eq:discrete_model-pb-var} as follows: {Find $u^\NN \in V^\NN $ such that}
\begin{equation}\label{eq:discrete_h_model-pb-var}
a_h^i(u^\NN,\varphi_i)=F_h^i(\varphi_i) \qquad \forall i = 1,\dots,n_\text{patches}\,.
\end{equation}

{Then, in order to cast Problem \eqref{eq:discrete_h_model-pb-var} into an optimization problem, we define \mbox{the residuals}}
\begin{equation}\label{eq:residuals}
r_{h,i}(w)=F_h^i(\varphi_i)-a_h^i(w,\varphi_i)\,, \qquad i  = 1,\dots,n_\text{patches}
\end{equation}
and the loss function
\begin{equation}\label{eq:loss-function}
R_h^2(w;{\cal P}) = \dfrac{1}{n_\text{patches}}\sum_{i =1}^{n_\text{patches}}\gamma_i r_{h,i}^2(w) \,,
\end{equation}
where $\gamma_i$ are suitable positive scaling coefficients. In this work, we use $\gamma_i=\text{area}(P_i)^{-1}$ to give the same importance to each patch. Note that this is equivalent to normalizing the quadrature rules involved in \eqref{eq:def-ah} and \eqref{eq:def-Fh}; this way, each residual $r_{h,i}$ can be regarded as a linear combination of the MF-VPINN value and derivatives independent of the size of the support of the patch $P_i$. We also highlight that the loss function depends on the choice of $\cal M$ since all the used test functions are generated starting from the corresponding mappings $M_i\in\cal M$. We are now interested in a practical procedure to obtain a set $\widetilde{\cal P}$ such that the approximate solution computed minimizing $R_h^2(\cdot;\widetilde{\cal P})$ is as accurate as possible with $\widetilde{\cal P}$ being as small as possible. 

\subsection{The a Posteriori Error Estimator}\label{sec:aposteriori}
The goal of this section is to derive an error estimator associated with an arbitrary patch $P_i$, with $i\in\{1,\dots,n_\text{patches}\}$. To do so, we rely on the a posteriori error estimator proposed in \cite{berrone2022solving}. It has been proven to be efficient and reliable; therefore, such an estimator allows us to know where the error is larger without knowing the exact solution of the PDE. Let us consider the patch $P_i$, formed by the triangles $T_{i,1},\dots,T_{i,N_{\hat P}}$ and a triangulation ${\cal T}_i$ of $\Omega$ such that $T_{i,j}\in{\cal T}_i$, for every $j=1,\dots,N_{\hat P}$. We remark that the triangulation ${\cal T}_i$ does not have to be explicitly generated; it is only used to properly define all the quantities introduced in \cite{berrone2022solving} required to derive the proposed error estimator. 

Let $V_h^i=\text{span}\{\psi_j^i:1\le j\le \text{dim}\left(V_h^i\right)\}$ be the space of piecewise linear functions defined on ${\cal T}_i$. Where $\{\psi_j^i:1\le j\le \text{dim}\left(V_h^i\right)\}$ is a Lagrange basis of $V_h^i$. It is then possible to define two constants $c_h^i$ and $C_h^i$, with $0<c_h^i<C_h^i$, such that
\begin{equation}\label{eq:norm-equiv-Vh}
c_h^i \vert v \vert_{1, \Omega} \leq \Vert \boldsymbol{v} \Vert_2 \leq C_h^i \vert v \vert_{1, \Omega} \qquad \forall v \in V_h^i \,,
\end{equation}
where $v=\sum_{j=1}^{\text{dim}\left(V_h^i\right)}v_j\psi_j^i$ is an arbitrary element of $V_h^i$ associated with the expansion coefficients $\boldsymbol{v}=\left\{v_1,\dots,v_{\text{dim}\left(V_h^i\right)}\right\}$ and $\Vert \boldsymbol{v} \Vert_2 = \left( \sum_{j=1}^{\text{dim}\left(V_h^i\right)} v_i^2 \right)^{1/2}$.

Then, given an integer $k\geq0$, for any element $E\in{\cal T}_i$, we define the projection operator $\Pi_{E,k} : L^2(E) \to \mathbb{P}_k(E)$ such that
\begin{equation}\label{eq:piek_property}
\int_E \Pi_{E,k} \phi = \int_E \phi \qquad \forall \phi \in L^2(E) \,.
\end{equation}

We also denote by {$\{(\bm{\xi}_\ell^E,\omega_\ell^E) : \ell \in I^{E}\}$} a quadrature formula of order $q$ on $E$ and define the quadrature-based discrete seminorm:
\begin{equation}\label{eq:discr-norm}
\Vert v \Vert_{0, E, \omega}= \left(\sum_{\ell \in I^E} v^2(\bm{\xi}^E_\ell) \,\omega^E_\ell \right)^{1/2} \,.
\end{equation}

We require the weights and nodes of this quadrature rule to coincide with the ones introduced in \eqref{eq:patch_quad_rule} when $E$ is a triangle included in $P_i$ (i.e., when $E\in\{T_{i,1},\dots,T_{i,N_{\hat P}}\}$). We can now introduce all the terms involved in the a posteriori error estimator.

Let $\eta_{{\mathrm {rhs}},1}(E)$ and $\eta_{{\mathrm {rhs}},2}(E)$ be the quantities:
\begin{equation}\label{eq:eta-f}
\begin{split}
\eta_{{\mathrm {rhs}},1}(E) &= h_E \Vert f - \Pi_{E,q-1} f \Vert_{0,E} \,,\\[3pt]
\eta_{{\mathrm {rhs}},2}(E) &= h_E \Vert f - \Pi_{E,q-1} f \Vert_{0,E,\omega} + \Vert f - \Pi_{E,q} f \Vert_{0,E,\omega} \,.
\end{split}
\end{equation}

They measure the oscillations of the forcing term with respect to its polynomial projections in various norms. Similar oscillations are also measured for the diffusion, convection and reaction terms by the terms $\eta_{{\mathrm{coef}},i}(E)$ for $i=1,\dots,6$:
\begin{equation}\label{eq:eta-coef}
\begin{split}
\eta_{{\mathrm {coef}},1}(E)  &= \Vert  \mu \nabla u^\NN - \Pi_{E,q} (\mu \nabla u^\NN) \Vert_{0,E} \,, \\[3pt]
\eta_{{\mathrm {coef}},2}(E) &= h_E \Vert  \boldsymbol{\beta}\cdot \nabla u^\NN - \Pi_{E,q-1}( \boldsymbol{\beta}\cdot \nabla u^\NN)  \Vert_{0,E} \,, \\[3pt]
\eta_{{\mathrm {coef}},3}(E) &= h_E \Vert  \sigma u^\NN - \Pi_{E,q-1}( \sigma u^\NN)  \Vert_{0,E} \,, \\[3pt]
\eta_{{\mathrm  {coef}},4}(E)  &= \Vert  \mu \nabla u^\NN - \Pi_{E,q} (\mu \nabla u^\NN) \Vert_{0,E,\omega} \,, \\[3pt]
\eta_{{\mathrm  {coef}},5}(E) &= h_E \Vert  \boldsymbol{\beta}\cdot \nabla u^\NN - \Pi_{E,q-1}( \boldsymbol{\beta}\cdot \nabla u^\NN)  \Vert_{0,E,\omega} \,, \\[3pt]
& \qquad \qquad \qquad + \Vert  \boldsymbol{\beta}\cdot \nabla u^\NN - \Pi_{E,q}( \boldsymbol{\beta}\cdot \nabla u^\NN)  \Vert_{0,E,\omega} \\[3pt]
\eta_{{\mathrm  {coef}},6}(E)  &= h_E \Vert  \sigma u^\NN - \Pi_{E,q-1}( \sigma u^\NN)  \Vert_{0,E,\omega} \\[3pt]
& \qquad \qquad \qquad + \Vert  \sigma u^\NN - \Pi_{E,q}( \sigma u^\NN)  \Vert_{0,E,\omega}\,,
\end{split}
\end{equation}
where $u^\NN$ is the output of the neural network after the enforcement of the Dirichlet boundary conditions through the operator $B$ and $h_E$ is the diameter of $E$. Then, let us define the term $\eta_{{\mathrm res}}(E)$, which measures how well the equation is satisfied, as
\begin{equation}\label{eq:eta-res}
\eta_{{\mathrm {res}}}(E) = h_E \Vert \, {\mathrm{ bulk}}_E(u^\NN)  \, \Vert_{0,E}  +  h_E^{1/2} \sum_{e \subset \partial E} \Vert \,{\mathrm{ jump}}_e(u^\NN)   \, \Vert_{0,e} \,,
\end{equation}
where 
\begin{equation*}\label{eq:def-bulk}
{\mathrm {bulk}}_E(u^\NN) = \Pi_{E,q-1}f +\nabla \cdot \Pi_{E,q} (\mu \nabla u^\NN) - \Pi_{E,q-1}( \boldsymbol{\beta}\cdot \nabla u^\NN + \sigma u^\NN)
\end{equation*}
\begin{equation*}\label{eq:def-jump}
{\mathrm {jump}}_e(u^\NN) = \Pi_{E_1,q}(\mu \nabla u^\NN)\cdot \boldsymbol{n} - \Pi_{E_2,q}(\mu \nabla u^\NN)\cdot \boldsymbol{n}
\,.
\end{equation*}

Note that ${\mathrm {jump}}_e(u^\NN)$ measures the interelemental jumps of $\Pi_{E,q}(\mu \nabla u^\NN)$ across the edge $e$ with normal unit vector {$\boldsymbol{n}$} shared by the elements $E_1$ and $E_2$. 

Finally, we introduce the approximate elemental  forms:
\begin{equation}\label{eq:def-ahi}
a_h^{i,E}(w,v)= \sum_{\ell \in I^E} [\mu \nabla w \cdot \nabla v + \boldsymbol{\beta}\cdot \nabla w \, v + \sigma w v](\bm{\xi}^E_\ell) \,\omega^E_\ell \,, 
\end{equation}
\begin{equation}\label{eq:def-Fhi}
F_h^{i,E}(v) =  \sum_{\ell \in I^E} [ f v](\bm{\xi}^E_\ell) \,\omega^E_\ell \,,
\end{equation}
where {$\bm{\xi}^E_\ell$} and $\omega^E_\ell$, $\ell\in I^E$, are the nodes and weights used in Equation \eqref{eq:discr-norm}. With such forms, it is possible to define the residuals
\begin{equation*}\label{eq:residuals_mesh_i}
r_{h,i,j}(w)=\sum_{E\in {\cal T}_i}F_h^{i,E}(\psi_j^i)-a_h^{i,E}(w,\psi_j^i)\,, \qquad j  = 1,\dots,\text{dim}\left(V_h^i\right)
\end{equation*}
and the quantity $\eta_{\mathrm {loss}}(E)$ as
\begin{equation}\label{eq:eta-local-loss}
\eta_{\mathrm {loss}}(E) = C_h\sqrt{ \sum_{j\in I_h^E}  r_{h,i,j}^2(u^\NN)} \,.
\end{equation}

Here, denoting the support of the function $\psi_j^i\in V_h^i$ by $\mathrm {supp}\, \psi_j^i$, the elemental index set 
\[
I_h^E =\{ j \in I_h : E \subset {\mathrm {supp}}\, \psi_j^i\}
\]
 is the set containing the indices of the functions whose support contains $E$. It is then possible to estimate the error between the unknown exact solution $u$ and its MF-VPINN approximation $u^\NN$ by means of the computable quantities \mbox{in \cref{eq:eta-f,eq:eta-coef,eq:eta-res} and \eqref{eq:eta-local-loss} as}
\begin{equation}\label{eq:local_estimator}
\vert u - u^\NN \vert_{1,E} \lesssim  \left( \eta_{\mathrm {res}}^2(E)  + \eta_{\mathrm {loss}}^2(E) + \sum_{i=1}^6\eta_{{\mathrm {coef}},i}^2(E)  +  \sum_{i=1}^2\eta_{{\mathrm {rhs}},i}^2(E) \right)^{1/2} \,.
\end{equation}

Once more, we refer to \cite{berrone2022solving} for the proof of such a statement. 

We recall that our goal is to obtain a computable error estimator associated with a single patch $P_i$. When evaluated on an element $E\in P_i$, the quantity on the right-hand side of Equation \eqref{eq:local_estimator} implicitly depends on several elements in $V_h^i$ that do not belong to $P_i$ because of the presence of $\eta_{\mathrm {res}}^2(E)$ and $\eta_{\mathrm {loss}}^2(E)$. Therefore, such an estimator is not computable without generating the triangulation ${\cal T}_i$ and the corresponding space $V_h^i$. Instead, we look for an error estimator that does not control the error on the entire patch but only in a neighborhood ${\cal N}_i$ of its center {$\bm{c_{P_i}} = M_i\left(\bm{c_{\hat P}}\right)$}. This can be carried out by considering only the terms whose computation involves geometric elements containing {$\bm{c_{P_i}}$} and the only function $\psi_j^i$ that does not vanish on {$\bm{c_{P_i}}$}. Note that such a function is the function $\varphi_i$ defined in \eqref{eq:phi_def}. Therefore, the error estimator $\eta_i$ that controls the error in ${\cal N}_i$ can be computed as
\begin{equation}\label{eq:patch_estimator}
\eta_i =   \left[ \eta_{\mathrm {res},i}^2  + C_h^2r_{h,i}^2(u^\NN) + \sum_{j=1}^{N_{\hat P}}\left(\sum_{k=1}^6\eta_{{\mathrm {coef}},k}^2(T_{i,j})  +  \sum_{k=1}^2\eta_{{\mathrm {rhs}},k}^2(T_{i,j})\right) \right]^{1/2},
\end{equation}
where $\eta_{\mathrm {res},i}$ is defined as
\begin{equation}\label{eq:eta-res_i}
\eta_{{\mathrm {res}},i} = \sum_{j=1}^{N_{\hat P}}\left( h_{T_{i,j}} \Vert \, {\mathrm{ bulk}}_{T_{i,j}}(u^\NN)  \, \Vert_{0,{T_{i,j}}}  +  h_{P_i}^{1/2}  \Vert \,{\mathrm{ jump}}_{e_{i,j}}(u^\NN)   \, \Vert_{0,e_{i,j}}\right) \,.
\end{equation}

In \eqref{eq:eta-res_i}, we denote by $h_{P_i}$ the diameter of the patch $P_i$ and by $e_{i,j}, j=1,\dots,N_{\hat P}$ the edges connecting its vertices with {$\bm{c_{P_i}}$}. 

Since $\eta_i$ can be seen as an approximation of the right-hand side of \eqref{eq:local_estimator}, we use it as an indicator of the error $\vert u - u^\NN \vert_{1,{\cal N}_i}$. It is important to remark that $\eta_i$ can be computed without generating ${\cal T}_i$ and $V_h^i$. In fact, its computation involves only the function $\varphi_i$, the triangles partitioning $P_i$ and the edges connecting its vertices with its center.

\subsection{The Choice of $\cal M$ and $\cal P$}\label{sec:choice_of_M}
In this section, the procedure adopted to generate the set of test functions used to train the MF-VPINN is described. We propose an iterative approach, in which the MF-VPINN is initially trained with very few test functions, and then other test functions are added in the regions of the domain in which the $H_1$ norm of the error is larger. We anticipate that, as shown in Section \ref{sec:estimator_role}, generating test functions in regions where $r_{h,i}^2$ is large may not lead to accurate solutions because $r_{h,i}^2$ is not proportional to the $H_1$ error. Therefore, such a choice may increase the density of test functions where they are not required while maintaining only a few test functions in regions in which the error is large. Instead, we use the error indicator $\eta_i$ defined in \eqref{eq:patch_estimator}.

Let us initially consider a cover ${\cal P}_0 = \{P_i\}_{i=1}^{n_\text{patches}}$ of $\Omega$ comprising a few patches (i.e., $n_\text{patches}$ is a small integer) and the corresponding set of mappings ${\cal M}_0=\{M_i\}_{i=1}^{n_\text{patches}}$ and test functions $\{\varphi_i\}_{i=1}^{n_\text{patches}}$. These sets induce a loss function $R_h^2(w;{\cal P}_0)$ as defined in \eqref{eq:loss-function}, which is used to train an MF-VPINN. After this initial training, one computes $\eta_i^\gamma=\gamma_i\eta_i$ for each patch $P_i\in{\cal P}_0$ and stores the result in the array {$\bm \eta$~}$=\left[\eta_1^\gamma,\dots,\eta_{n_\text{patches}}^\gamma\right]$. Note that $\eta_i^\gamma$ is a suitable rescaling of $\eta_i$ to get rid of dependence from the size of $P_i$. Let us choose a threshold $1\le \tau_0\le n_\text{patches}$, sort {$\bm \eta$} in descending order obtaining {${\bm \eta}$}$_\text{sort}=\left[\eta_{s_1}^\gamma,\dots,\eta_{s_{n_\text{patches}}}^\gamma\right]$ (where we denote by $[s_1,\dots,s_{n_\text{patches}}]$ the index set corresponding to a suitable permutation of $[1,\dots,n_\text{patches}]$) and consider the vector {$\overline{\bm \eta}$}$_0=\left[\eta_{s_1}^\gamma,\dots,\eta_{s_{\tau_0}}^\gamma\right]$. It is possible to note that $\overline{\bm \eta}_0$ contains only the $\tau_0$ worst values of the indicator; it thus allows us to understand where the error is higher and where additional test functions are required to increase the \mbox{model accuracy}.

It is then possible to move forward with the second iteration of the iterative training. For each patch $P_i$ such that {$\eta_i^\gamma\in\overline{\bm\eta}_0$}, we generate $k_\text{new}$ new patches $P_i^k$, $k=1,\dots,k_\text{new}$ with centers inside $P_i$ and areas such that $\text{area}(P_i)<\cup_{k=1}^{k_\text{new}}\text{area}(P_i^k) < c\cdot \text{area}(P_i)$, where $c>1$ is a tunable parameter. In the numerical experiments, we use $c=1.25$. There exist different strategies to choose the number, the dimension, and the position of the centers of the new patches. Such strategies are described in Section \ref{sec:numerical_results} with particular attention to the effects of these choices on the MF-VPINN accuracy.

Let us denote by ${\cal P}_1$ the set ${\cal P}_1={\cal P}_0\cup\{P_{s_1}^k\}_{k=1}^{k_\text{new}}\cup\dots\cup\{P_{s_{\tau_0}}^k\}_{k=1}^{k_\text{new}}$ and by ${\cal M}_1$ the corresponding set of mappings. Then, it is possible to define the loss function $R_h^2(w;{\cal P}_1)$, continue the training of the previously trained MF-VPINN, compute the error indicator $\eta_i^\gamma$ for each patch $P_i\in{\cal P}_1$, and obtain the vector $\overline{\bm \eta}_1$ used to decide where to insert the new patches to generate ${\cal P}_2$. In general, iterating this procedure, it is possible to compute a set of patches ${\cal P}_m$ and of mappings ${\cal M}_m$ from the previously obtained sets ${\cal P}_{m-1}$ and ${\cal M}_{m-1}$. Technical optimization details are discussed in Section \ref{sec:implementation}.

\section{Numerical Results}\label{sec:numerical_results}
In this section, we provide several numerical results to show the performance of the training strategy described in Section \ref{sec:choice_of_M}. In Section \ref{sec:implementation}, we describe the structure of the MF-VPINN implementation and highlight some details that have to be taken into account in order to increase the efficiency of the training phase. Different strategies to choose the position of the new patches are discussed in Section \ref{sec:performance}. The importance of the use of the error indicator is remarked in Section \ref{sec:estimator_role} with additional numerical examples. {An example on a more complex domain is shown in Section \ref{sec:holes} to discuss some ideas to adapt the proposed strategies in more complex domains.}

\subsection{Implementation Details}\label{sec:implementation}
The computer code used to perform the experiment is implemented in Python using the Python package Tensorflow \cite{tensorflow2015-whitepaper} to generate the neural network architecture and train the MF-VPINN. Using the notation introduced in Section \ref{sec:problem_discretization}, the used neural network consists of $L=5$ layers with $N_\ell=50$ neurons in each hidden layer \mbox{(i.e., for $\ell=1,\dots,L-1$}); the activation function is the hyperbolic tangent in each hidden layer. For the first iteration of the iterative training, the neural network weights in the $\ell$-th layer are initialized with a glorot normal distribution, i.e., a truncated normal distribution with mean 0 and standard deviation equal to $\sqrt{2/(N_{\ell-1}+N_\ell)}$. Then, for the subsequent iterations, their are initialized with the weights obtained at the end of the previous one. 

During the first iteration of the training (during the minimization of $R_h^2(\cdot;{\cal P}_0)$), the optimization is carried out by exploiting the ADAM optimizer \cite{kingma2014adam} with an exponentially decaying learning rate from $10^{-2}$ to $10^{-4}$ and with the second-order L-BFGS optimizer \cite{wright1999numerical}. Then, from the second training iteration, we only use the L-BFGS optimizer. We remark that L-BFGS allows a very fast convergence but only if the initial starting point is close enough to the problem's solution. Therefore, in the first training iteration, we use ADAM to obtain a first approximation of the solution that is then improved via L-BFGS. Then, since the $m$-th training iteration starts from the solution computed during the $(m-1)$-th one, we assume that the starting point is close enough to the solution of the new optimization problem (associated with a difference loss function with more patches) and we only use L-BFGS to increase the training efficiency.

During the $m$-th iteration of the training, the training set consists of all the quadrature nodes ${\bm \xi}_\ell^i$, for any $\ell\in I^{\hat P}$ and for any patch $P_i\in{\cal P}_m$ as defined in \eqref{eq:patch_quad_rule}. The order of the chosen quadrature rule is $q=3$ inside each triangle. The Dirichlet boundary conditions are imposed by means of the operator $B$ defined in \eqref{eq:B_definition}. In this operator, for our first numerical test, the function $\phi$ is a polynomial bubble vanishing on $\Gamma$ and $\overline g$ is the output of a neural network trained to interpolate the boundary data. {For the numerical test in Section \ref{sec:holes}, instead, $\phi$ is computed as in \cite{berrone2022enforcing} and $g=0$.} To decrease the training time, the functions $\phi$, $\nabla\phi$, $\overline g$ and $\nabla\overline g$ are evaluated only once at the beginning of the $m$-th training iteration and they are then combined to evaluate $Bu^\NN$ and its gradient (where $u^\NN$ is the output of the last layer of the neural network). The derivatives of $u^\NN$ and $\overline g$ are computed via automatic differentiation \cite{baydin2018automatic} due to the complexity of their analytical expressions.

{The output of the model is the value of the function $Bu^\NN$ and its gradient evaluated at the input points. Such values are then suitably combined using sparse and dense tensors to compute the quantity $R_h^2(Bu^\NN;{\cal P}_m)$. The sparse tensors contain the evaluation of $\varphi_i$ and $\nabla\varphi_i$ at each input point, whereas the dense ones store the quadrature weights, the vector $\bm\gamma=\{\gamma_i\}_{i=1}^{n_\text{patches}}$ and the evaluation of $\mu$, {${\boldsymbol{\beta}}$}, $\sigma$ and $f$ at the input points. We highlight that all these tensors have to be computed once at the beginning of the $m$-th training iteration (updating the ones of the $(m-1)$-th iteration) to significantly decrease the training computational cost.}

As discussed in Section \ref{sec:problem_discretization}, we assume that all the patches and test functions can be generated from a reference patch $\hat P$. For each patch $P_i\in {\cal P}_m$, one has to generate all the data structures required to assemble the loss function and the error indicator $\eta_i$. To do so, it is possible to explicitly construct all the tensors required to assemble the term ${\hat a}_h(w,\hat \varphi)$  and all the terms involved in the computation of the reference error indicator $\hat\eta$ only once, at the beginning of the first iteration of the training. Then, all these tensors can be suitably rescaled to obtain the ones corresponding to the patches and test functions involved in the loss function and error indicators computations. 

To stabilize the MF-VPINN, we introduce the $L^2$ regularization term 
\begin{equation*}\label{eq:l2_reg}
\mathcal L_{\text{reg}}({u}^{\mathcal{N\!N}}) = \lambda_{\text{reg}}\Vert \mathbf{u}^{\mathcal{N\!N}}\Vert_2^2,
\end{equation*}
where {$\mathbf{u}^{\mathcal{N\!N}}$} is the set of weights of the neural network introduced in Section \ref{sec:problem_discretization}. In our numerical experiments, we use $\lambda_\text{reg}=10^{-5}$. During the $m$-th iteration of the training, such a quantity is added to $R_h^2(B{u}^{\mathcal{N\!N}};{\cal P}_m)$ to obtain the training loss function 
\begin{equation}\label{eq:loss+reg}
{\mathcal L}_m({u}^{\mathcal{N\!N}}) = R_h^2(B{u}^{\mathcal{N\!N}};{\cal P}_m) + \mathcal L_{\text{reg}}(\mathbf{w}^{\mathcal{N\!N}}),
\end{equation}
which has to be minimized accurately enough. Indeed, if ${\mathcal L}_m$ is minimized poorly, the new patches ${\cal P}_{m+1}\backslash{\cal P}_m$ may be added in regions where they are not necessary because the accuracy of $Bu^\NN$ may still improve during the training  and may not be inserted in areas where they are required. Note that, in order to compute the numerical solution, the MF-VPINN has to be trained multiple times with a different set of patches ${\cal P}_m$ to minimize the losses $\{{\mathcal L}_m\}$. Since such an iterative training may be expensive, we propose an early stopping strategy \cite{prechelt1998early} based on the discussed error indicator to reduce its computational cost. In its basic version, early stopping consists of evaluating a chosen metric on a validation set in order to know when the neural network accuracy on data that are not present in the training set start worsening. Interrupting the training there prevents overfitting and improves generalization. In our context, instead, we can directly track the behavior of the MF-VPINN $H^1$ error on each patch through the corresponding error indicator to understand when it stops decreasing. Therefore, given the set of patches ${\cal P}_m$, the chosen metric is the linear combination $ES_m=\sum_{i=1}^{\text{dim}({\cal P}_m)}\eta_i^\gamma$. Numerical results showing the performance of this strategy are presented in Sections \ref{sec:performance} {and \ref{sec:holes}}.

\subsection{Adaptive Training Strategies}\label{sec:performance}
Let us consider the Poisson problem:
\begin{equation}\label{eq:model-pb-poisson}
\begin{cases}
-\Delta u = f & \text{in \ } \Omega\,, \\
\ \ \, u=g & \text{on \ } \Gamma \,, \end{cases}
\end{equation}
defined on the unit square $\Omega=(0,1)^2$. The forcing term $f$ and the boundary condition $g$ are chosen such that the exact solution is, in polar coordinates,
\begin{equation}\label{eq:sol4}
u(r,\theta) = r^\frac 23 \sin\left(\frac 23\left(\theta+\frac \pi2\right)\right).
\end{equation}

\begin{figure}[t!]
\centering 
  \includegraphics[width=0.75\linewidth]{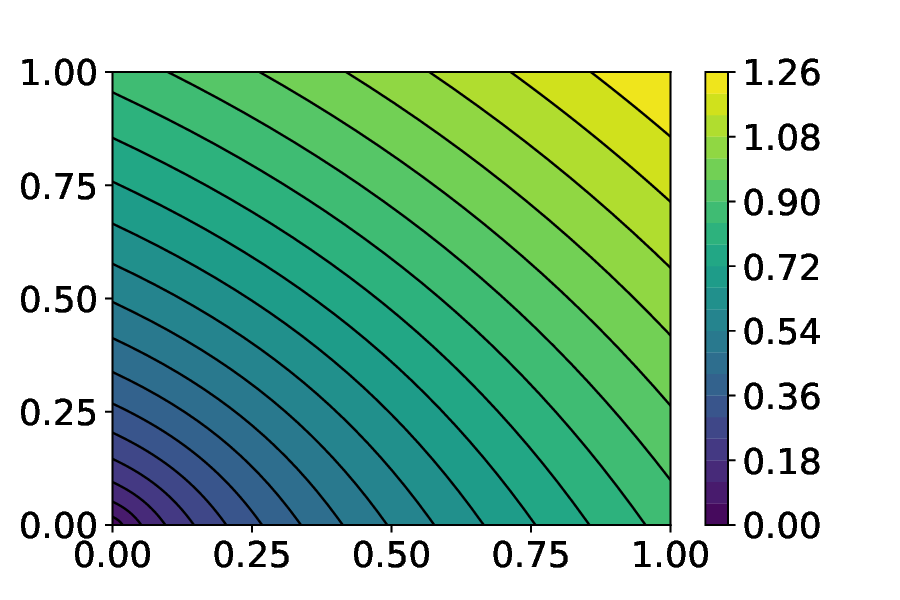} 
  \caption{Graphical representation of the solution $u$ in \eqref{eq:sol4}.}
  \label{fig:solution4}
\end{figure}

We use this function, represented in Figure~\ref{fig:solution4}, because the solution $u$ is such that $u\in H^{5/3-\varepsilon}(\Omega)$ but $u\in C^\infty(\Omega \backslash {\cal N}_0)$, where we denote by ${\cal N}_0$ a neighborhood of the origin. Therefore, we know that an efficient distribution of patches has to be characterized by a high density only near the origin.

{Below, we propose, in order of complexity, three alternatives to construct the new patches after having marked the ones with the higher error indicator. The first strategy is the most simple and intuitive, and the new patches are randomly generated with centers inside the marked patches, whereas the second strategy and third one place the new centers on a small local cartesian grid to ensure a more regular distribution. The difference between the second and the third strategies is that the marked patches are removed to increase the efficiency and we add a constraint to the marking procedure to ensure more regular distributions of the new patches.}

\begin{itemize}
\item[] \hspace{-0.75cm}\textit{{Strategy \#1: Random patch centers with uniform distribution} 
}
\end{itemize}

To solve Problem \eqref{eq:model-pb-poisson}, as a first strategy, we consider the reference patch $\hat P=(0,1)^2$ and generate a sequence of sets of patches. During the first training iteration, we use ${\cal P}_0=\{\hat P\}$ since this is already a cover of $\Omega$. During the second iteration, we enrich the set of patches as ${\cal P}_1={\cal P}_0 \cup \{P_1,P_2,P_3,P_4\}$ where $P_1,P_2,P_3$ and $P_4$ are squared patches with edge $h_i=0.6$, $i=1,\ldots,4$ and centers 
\[
\bm{c_{P_1}}=(0.3,0.3), \hspace{1cm} \bm{c_{P_2}}=(0.7,0.3),
\]\[
\bm{c_{P_3}}=(0.3,0.7), \hspace{1cm} \bm{c_{P_4}}=(0.7,0.7).
\]

This allows us to start from a homogeneous distribution of patches before utilizing the error indicator to choose the location of the new patches. Then, to decide how many patches have to be added to ${\cal P}_{m-1}$ to generate ${\cal P}_m$, we choose $\widetilde \tau_m$ such that
\begin{equation}\label{eq:tau_tilde_def}
\widetilde \tau_m = \text{dim}\left(\left\{\widetilde\tau\in\{1,\dots,\text{dim}({\cal P}_{m-1}\} : \dfrac{\sum_{i=1}^{\widetilde\tau}\eta_{s_i}^\gamma}{\sum_{i=1}^{\text{dim}({\cal P}_{m-1})}\eta_i^\gamma}<0.75\right\}\right) +1
\end{equation}
and fix
\begin{equation}\label{eq:tau_def}
\tau_m = \min(\lceil0.3\cdot\text{dim}({\cal P}_{m-1})\rceil, \widetilde\tau_m).
\end{equation}

Note that \eqref{eq:tau_tilde_def} allows us to consider the smallest set of patches such that the corresponding error indicators contribute at least $75\%$ of the global error indicator $ES_{m-1}$, whereas \eqref{eq:tau_def} is considered to limit the maximum number of patches that can be added for efficiency reasons.

Then, to generate the generic set of patches ${\cal P}_m$, we fix a multiplication factor $C_M$ to decide how many new patches have to be inserted inside each patch $P_i$ such that {$\eta_i^\gamma\in\overline{\bm \eta}_{m-1}$.} Inside each chosen patch $P_i$, $C_M$ centers {$\bm{\widetilde{c}_{P_i^k}}$}$~=(\widetilde{x}_i^k,\widetilde{y}_i^k)$, $k=1,\dots,C_M$, are randomly generated with a uniform distribution and the new patches' edges' lengths are chosen as $h_i^k=\lambda\frac{A_\text{ratio}}{\sqrt{C_M}}h_i$. Here, $\lambda$ is a random real value from the uniform distribution $U\left(\left[\frac{9}{10}, \frac{10}{9}\right]\right)$, and the scaling coefficient $\frac{A_\text{ratio}}{\sqrt{C_M}}$ is chosen such that the sum of the areas of the new patches is $A_\text{ratio}$ times the area of the original patch $P_i$. In the numerical experiments, we use $A_\text{ratio}=1.25$. This way, it is possible to allow the new patches to overlap and keep the area of the region $P_i\backslash\left(\cup_{k=1}^{C_M}P_i^k\right)$ reasonably small.

We remark that, with this strategy, it may happen that some patches are outside $\Omega$. In order to avoid this risk, we move the centers $\bm{\widetilde{c}_{P_i^k}}$ to obtain the actual patches centers {$\bm{c_{P_i^k}}$} \mbox{as follows}: 
\begin{equation}\label{eq:centers_update}
\bm{c_{P_i^k}}=(x_i^k,y_i^k) \leftarrow \left(\max\left\{\min\left\{\widetilde{x}_i^k,1-\frac{h_i^k}{2}\right\}, \frac{h_i^k}{2}\right\},    
\max\left\{\min\left\{\widetilde{y}_i^k,1-\frac{h_i^k}{2}\right\}, \frac{h_i^k}{2}\right\}\right).
\end{equation}

We remark that, when the patch $P_i$ is very close to a vertex of the domain, it is possible that multiple original centers {$\bm{\widetilde{c}_{P_i^k}}$} are such that the distance of both $\widetilde{x}_i^k$ and $\widetilde{y}_i^k$ from the $x$ and $y$ coordinates of the domain vertex is smaller than $h_i^k/2$. In this case, it is important to consider the random coefficient $\lambda$ in the definition of $h_i^k$ to avoid updating all these centers with the same point; otherwise, multiple new patches would coincide (because they would share the same center and size).

\begin{figure}[t!]
\centering 
  \includegraphics[width=0.6\linewidth]{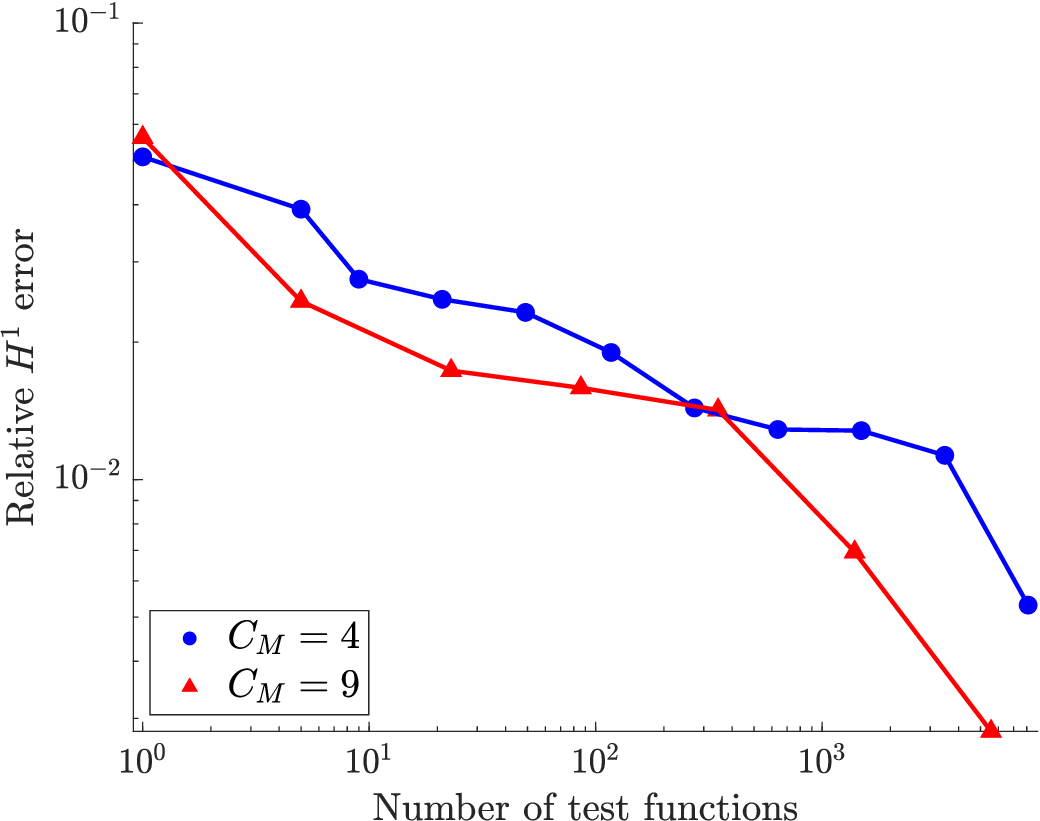}
  \caption{\textit{{Strategy \#1}
}: Relative $H^1$ errors obtained at the end of each training iteration for $C_M=4$ (blue circles) and $C_M=9$ (red triangles).}
  \label{fig:error_random_49_energy75}
\end{figure}

\begin{figure}[t!]
\centering 
\begin{subfigure}[H]{0.325\linewidth}
     \includegraphics[width=0.99\columnwidth,keepaspectratio,clip]{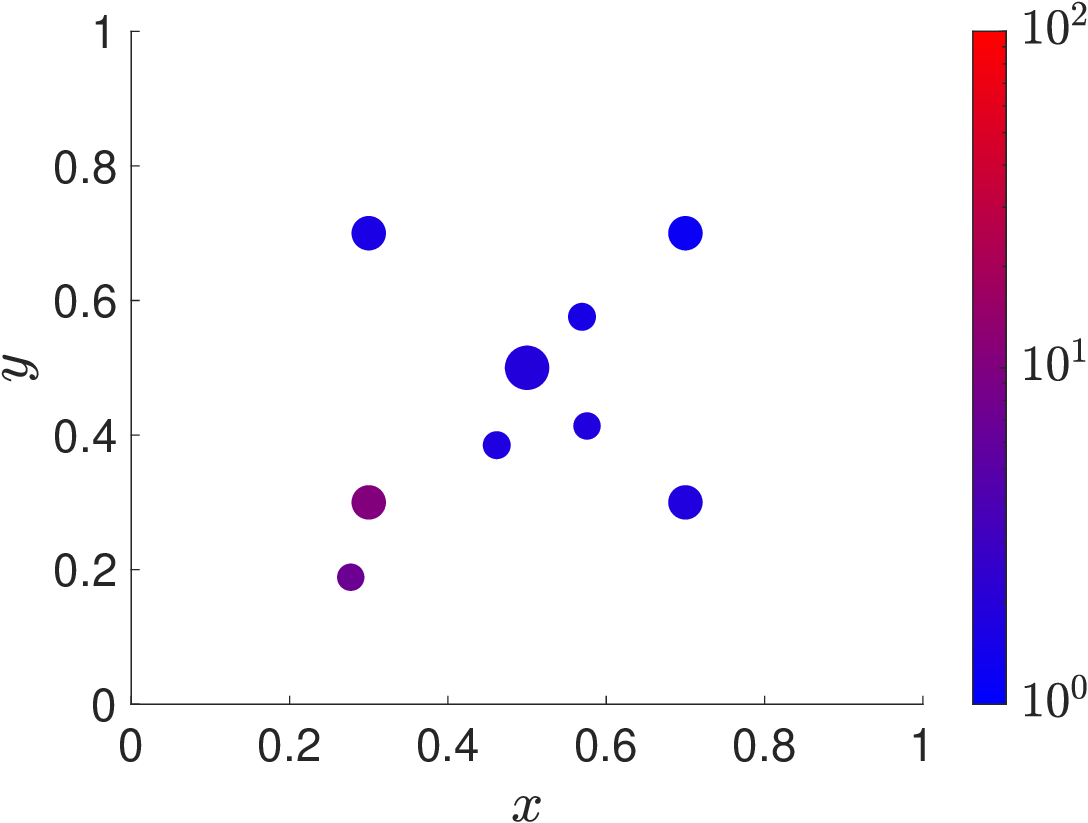} 
  \subcaption{\centering}
\end{subfigure}
\begin{subfigure}[H]{0.325\linewidth}
     \includegraphics[width=0.99\columnwidth,keepaspectratio,clip]{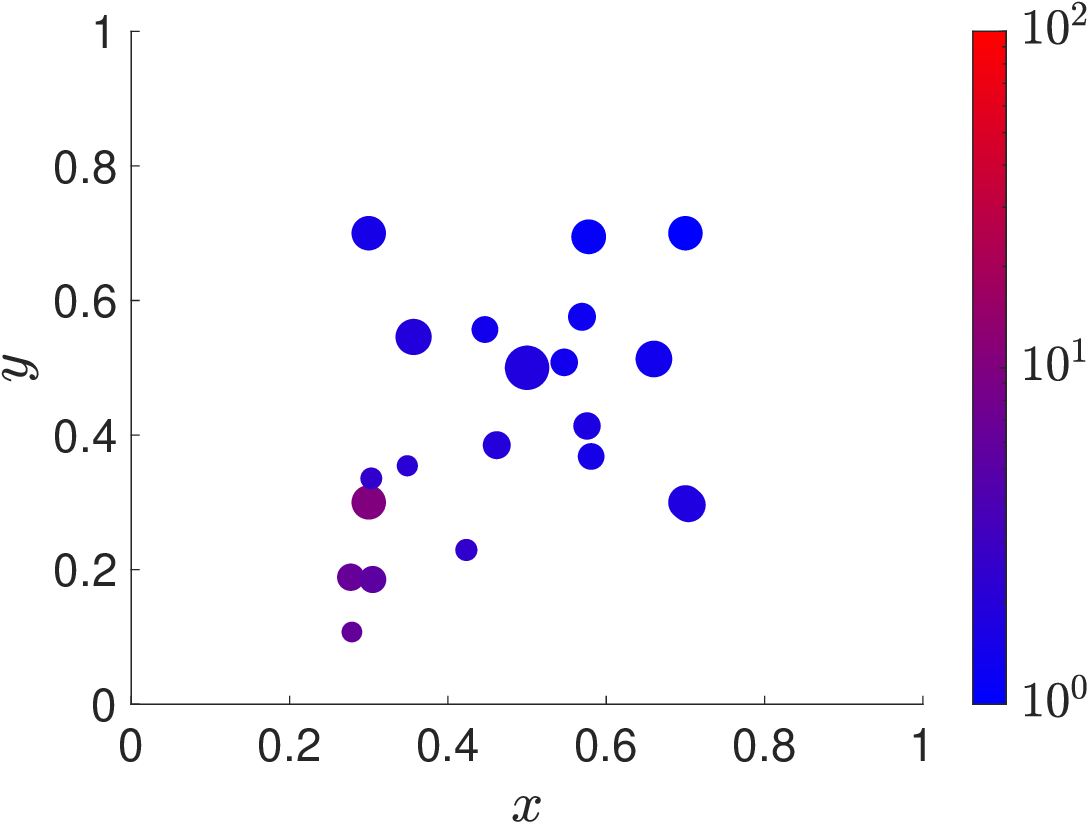} 
   \subcaption{\centering}
\end{subfigure}
\begin{subfigure}[H]{0.325\linewidth}
     \includegraphics[width=0.99\columnwidth,keepaspectratio,clip]{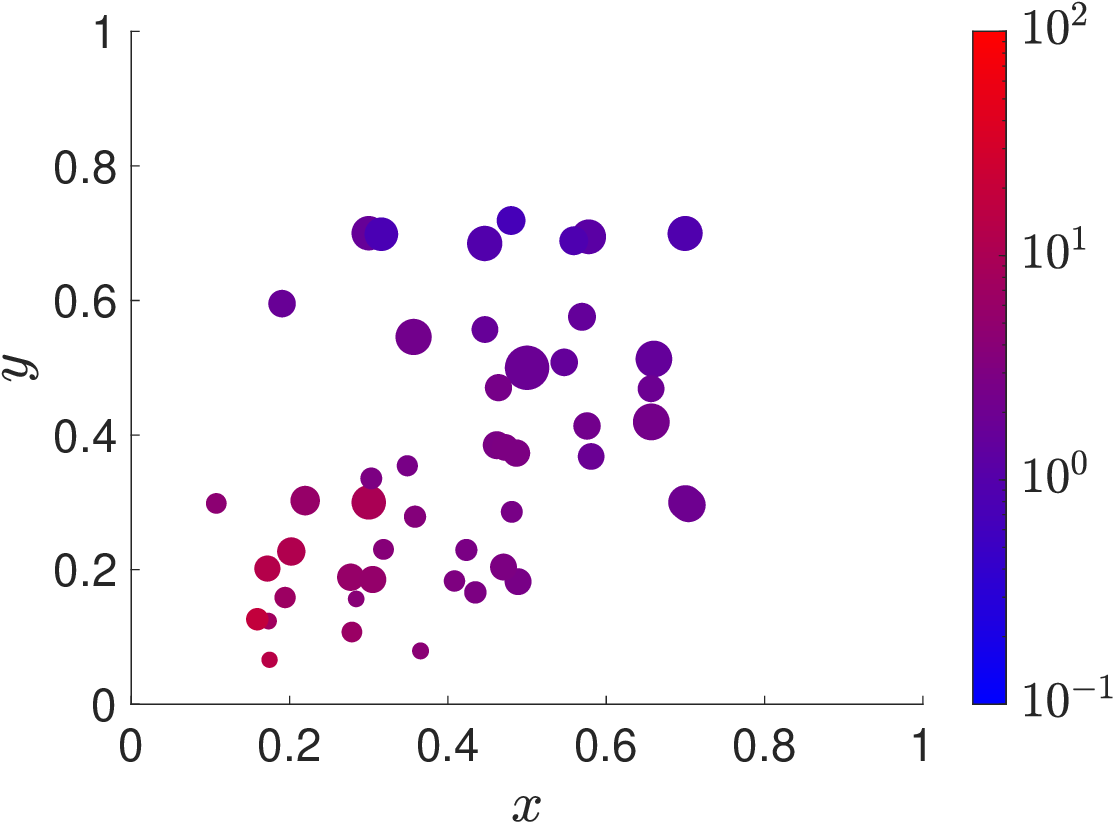} 
   \subcaption{\centering}
   \label{fig:random_4_energy75_estimators_step_3}
\end{subfigure}

\medskip

\begin{subfigure}[H]{0.325\linewidth}
     \includegraphics[width=0.99\columnwidth,keepaspectratio,clip]{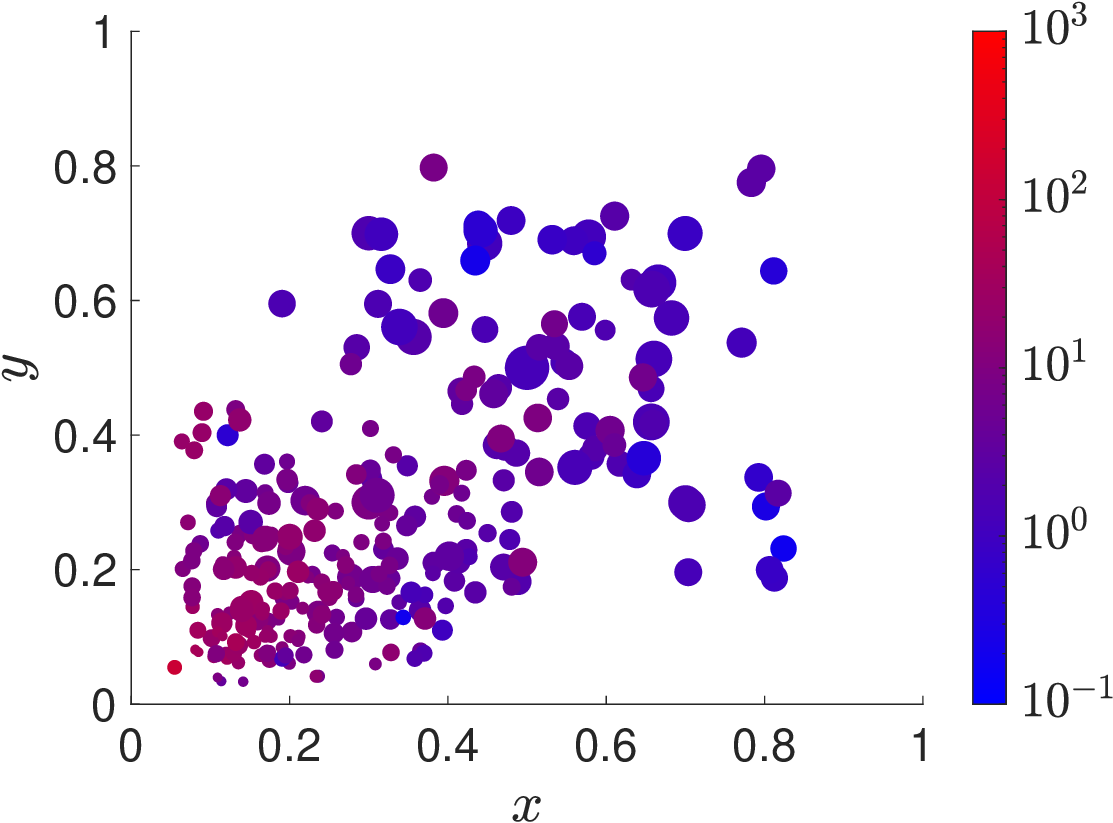} 
   \subcaption{\centering}
\end{subfigure}
\begin{subfigure}[H]{0.325\linewidth}
     \includegraphics[width=0.99\columnwidth,keepaspectratio,clip]{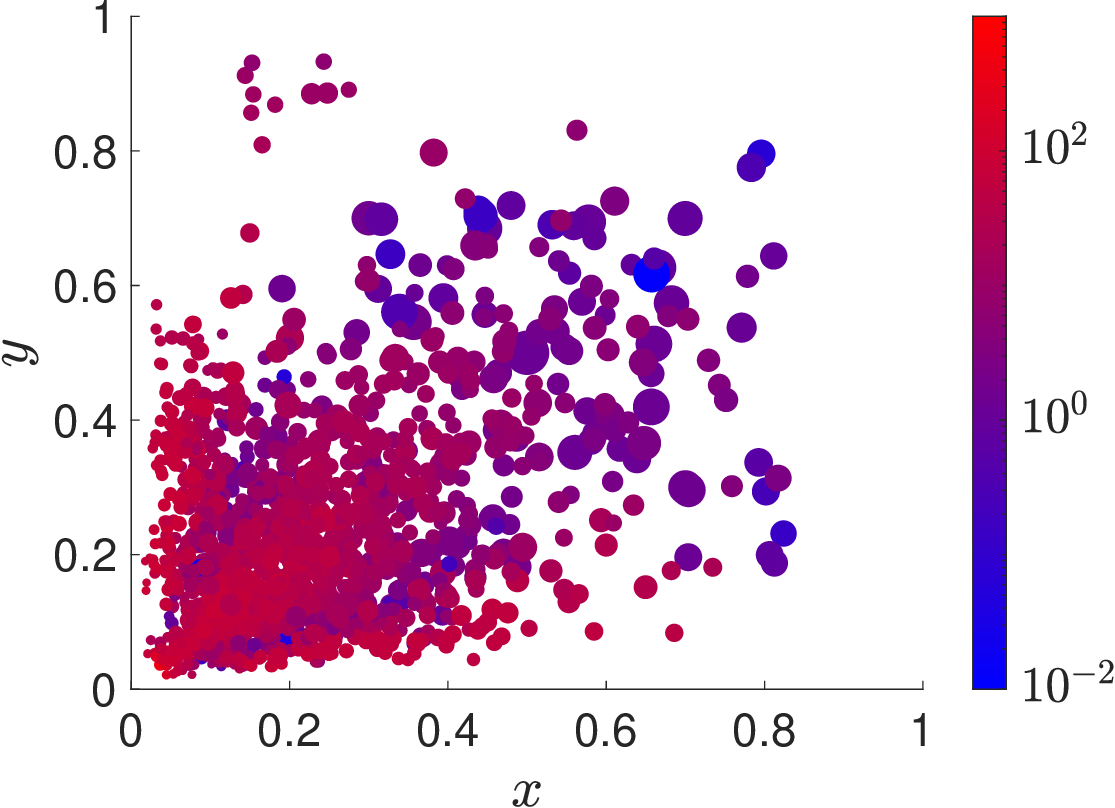} 
   \subcaption{\centering}
\end{subfigure}
\begin{subfigure}[H]{0.325\linewidth}
     \includegraphics[width=0.99\columnwidth,keepaspectratio,clip]{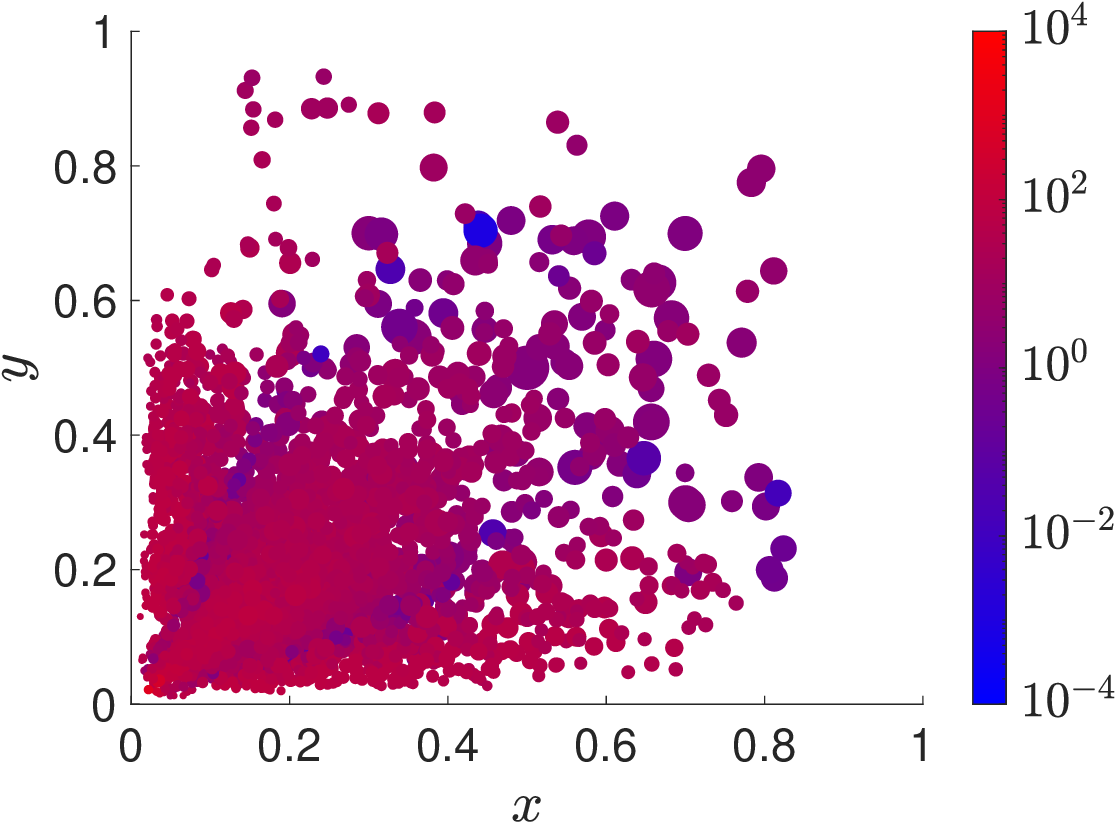} 
   \subcaption{\centering}
\end{subfigure}
\medskip
  \caption{\textit{{Strategy \#1}
}: Patches used to train the MF-VPINN with $C_M=4$. Each dot represents a patch $P_i$, its position is the center $\bm{c_{P_i}}$ of the patch, its size is proportional to the patch size $h_i^2$, and its color is associated with the quantity $\eta_i^\gamma$. {(\textbf{a}) Representation of ${\cal P}_2$; (\textbf{b}) Representation of ${\cal P}_3$; (\textbf{c}) Representation of ${\cal P}_4$; (\textbf{d}) Representation of ${\cal P}_6$; (\textbf{e}) Representation of ${\cal P}_8$; (\textbf{f}) Representation of ${\cal P}_9$.} 
}
  \label{fig:random_4_energy75_estimators}
\end{figure}

\begin{figure}[t!]
\centering 
\begin{subfigure}[H]{0.325\linewidth}
  \includegraphics[width=0.99\columnwidth,keepaspectratio,clip]{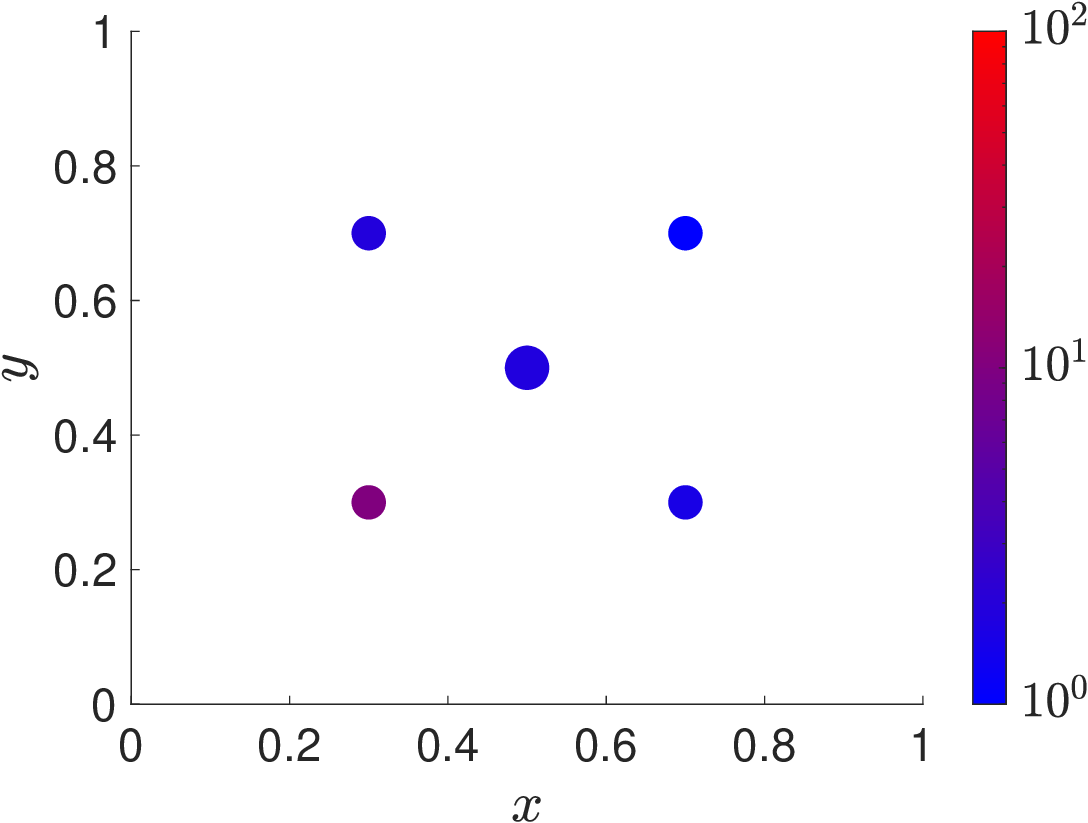} 
   \subcaption{\centering}
\end{subfigure}
\begin{subfigure}[H]{0.325\linewidth}
  \includegraphics[width=0.99\columnwidth,keepaspectratio,clip]{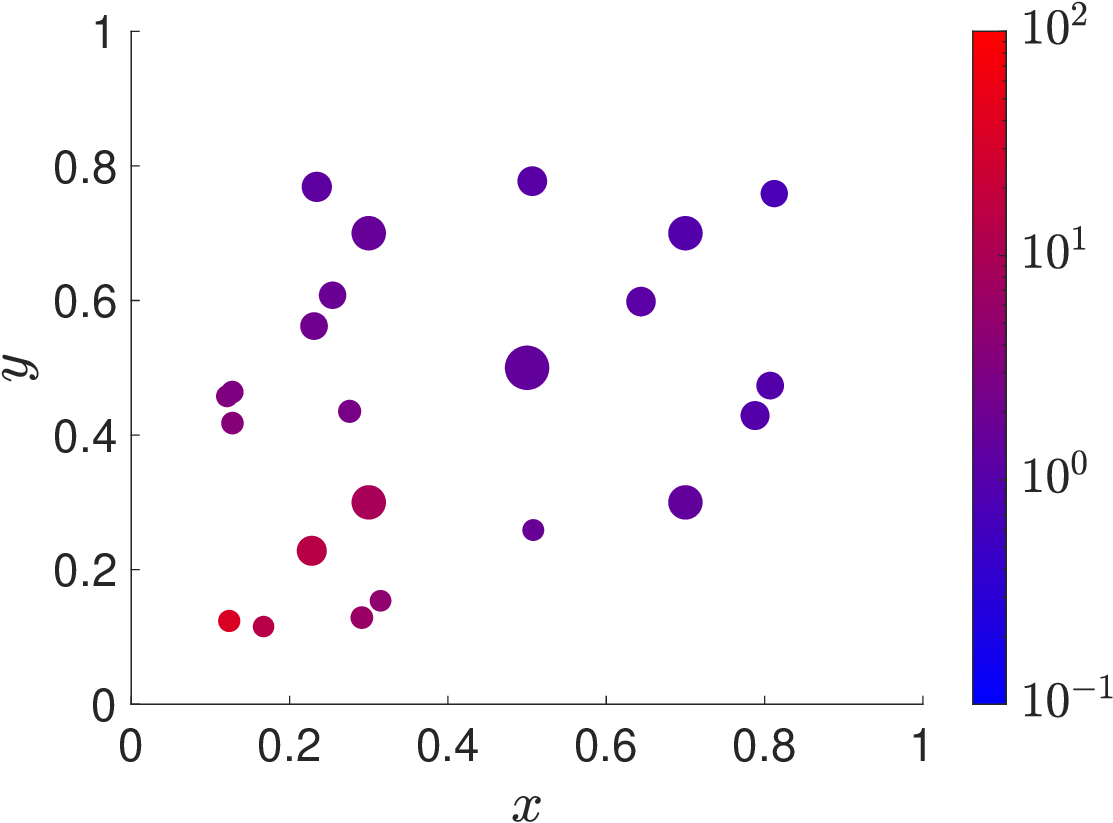} 
   \subcaption{\centering}
\end{subfigure}
\begin{subfigure}[H]{0.325\linewidth}
  \includegraphics[width=0.99\columnwidth,keepaspectratio,clip]{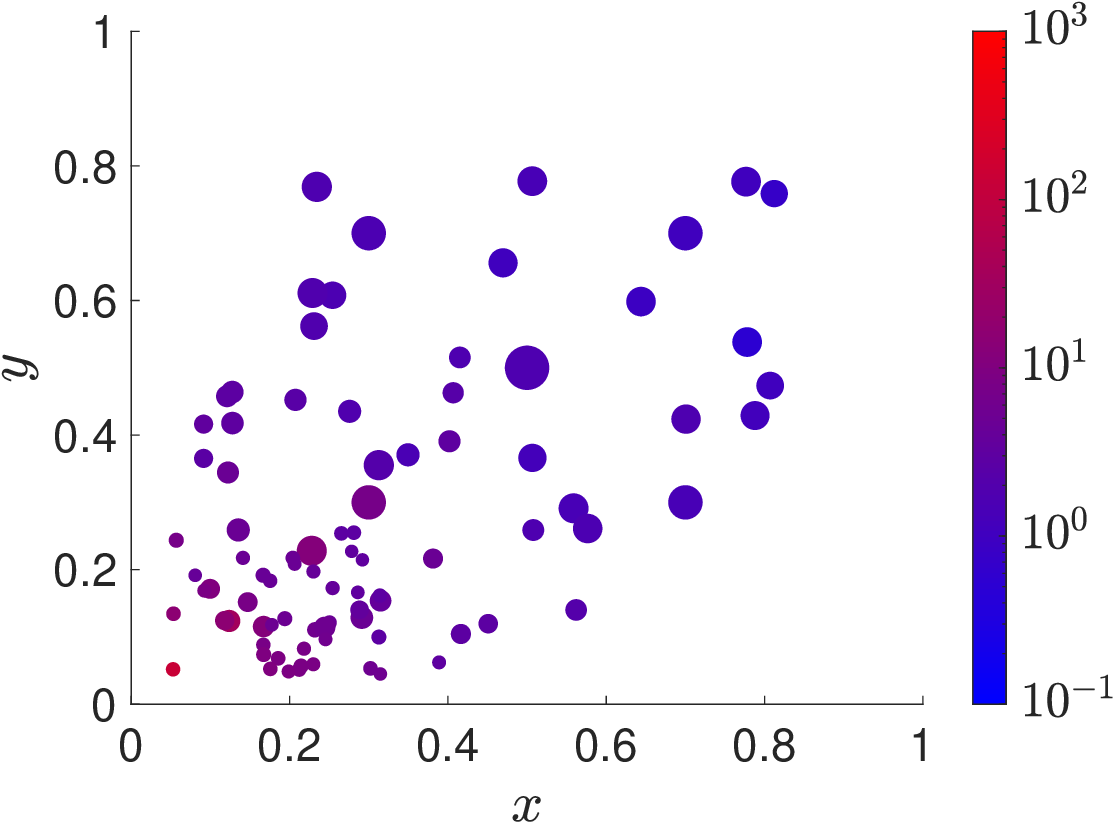}
   \subcaption{\centering}
\end{subfigure}

\medskip

\begin{subfigure}[H]{0.325\linewidth}
  \includegraphics[width=0.99\columnwidth,keepaspectratio,clip]{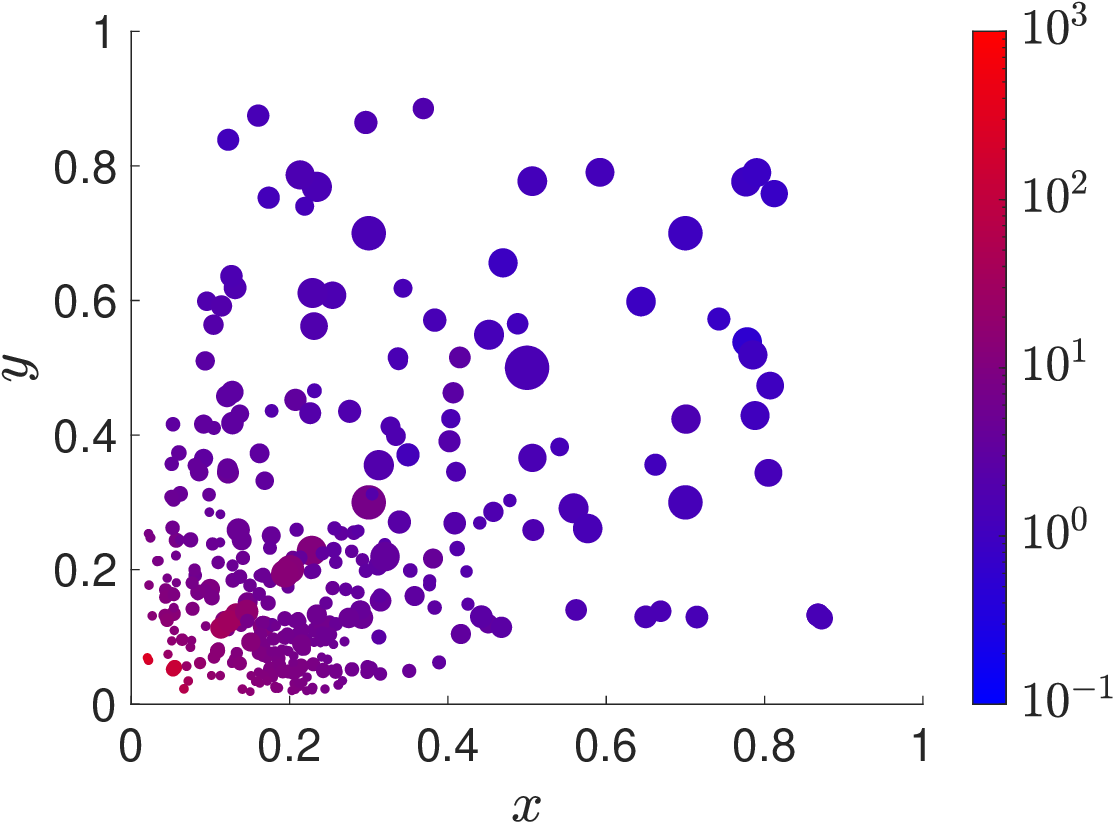} 
   \subcaption{\centering}
\end{subfigure}
\begin{subfigure}[H]{0.325\linewidth}
  \includegraphics[width=0.99\columnwidth,keepaspectratio,clip]{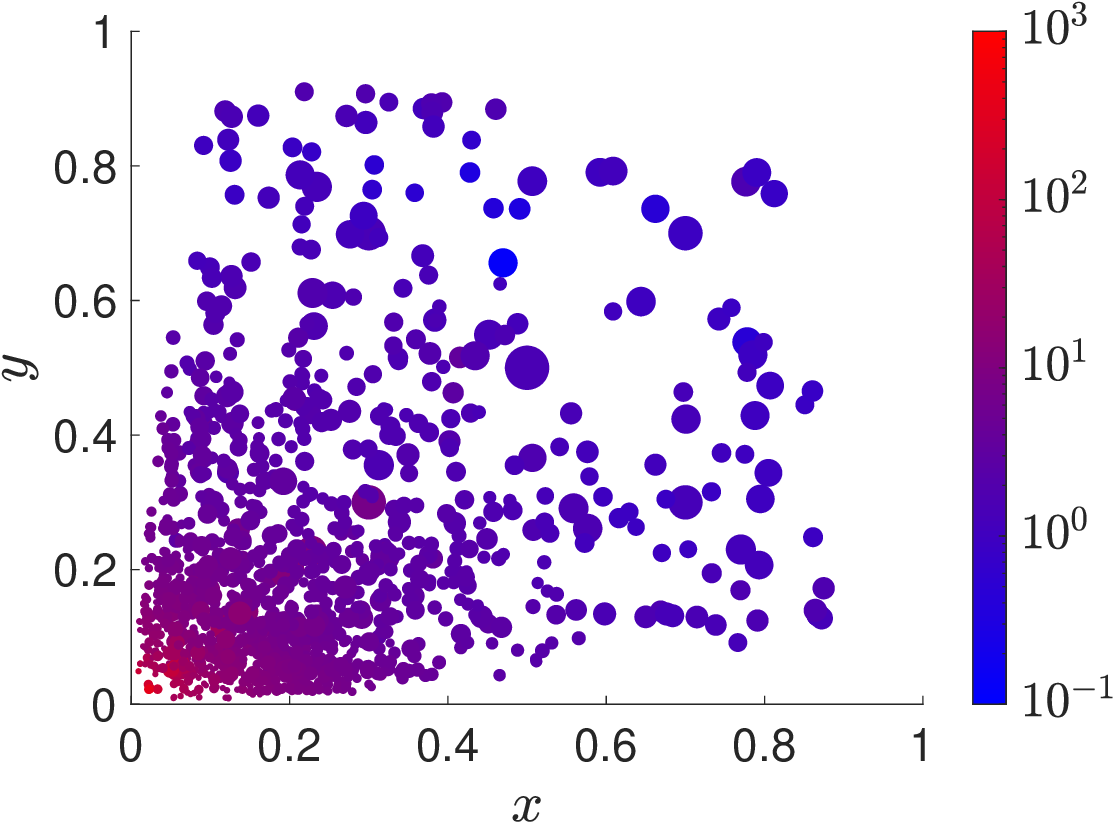} 
   \subcaption{\centering}
\end{subfigure}
\begin{subfigure}[H]{0.325\linewidth}
  \includegraphics[width=0.99\columnwidth,keepaspectratio,clip]{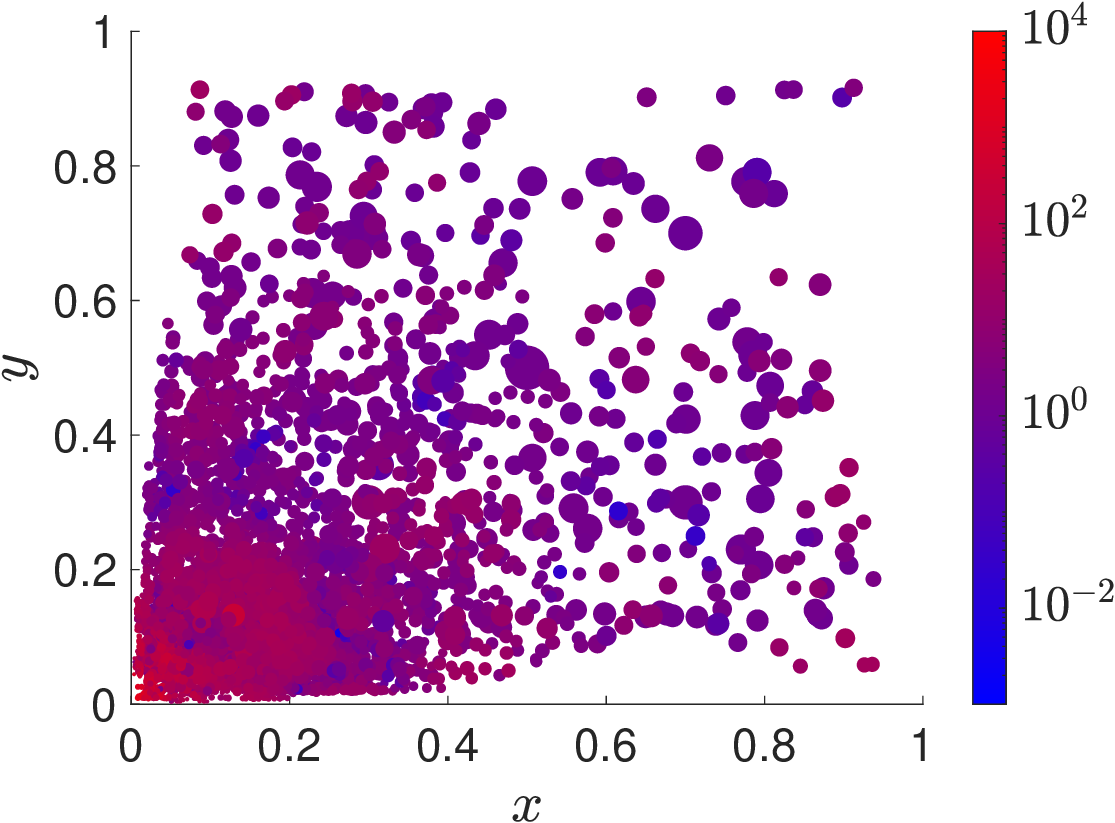}
   \subcaption{\centering}
   \label{fig:random_9_energy75_estimators_6}
\end{subfigure}
\medskip
  \caption{\textit{{Strategy \#1}
}: Patches used to train the MF-VPINN with $C_M=9$. Each dot represents a patch $P_i$, its position is the center $\bm{c_{P_i}}$ of the patch, its size is proportional to the patch size $h_i^2$, and its color is associated with the quantity $\eta_i^\gamma$. {(\textbf{a}) Representation of ${\cal P}_1$; (\textbf{b}) Representation of ${\cal P}_2$; (\textbf{c}) Representation of ${\cal P}_3$; (\textbf{d}) Representation of ${\cal P}_4$; (\textbf{e}) Representation of ${\cal P}_5$; (\textbf{f}) Representation of ${\cal P}_6$.} 
}
  \label{fig:random_9_energy75_estimators}
\end{figure}

For the numerical test, we consider $C_M=4$ and $C_M=9$. Using significantly more accurate quadrature rules, we compare the approximate solution with the exact one defined in \eqref{eq:sol4} and compute the relative $H^1$ error ${\Vert u - u^\NN\Vert_1}/{\Vert u \Vert_1}$ at the end of each training iteration. The obtained errors are shown as blue circles ($C_M=4$) and red triangles ($C_M=9$) in Figure \ref{fig:error_random_49_energy75}. It can be noted that, with both values of $C_M$, when more patches are used, the error is smaller, even though the convergence rate is limited by the low regularity of the solution. It is also interesting to observe the positions and sizes of the used patches; such information is summarized in Figures \ref{fig:random_4_energy75_estimators} and \ref{fig:random_9_energy75_estimators}. In such figures, each dot is in the center of a patch $P_i$, and its size and color represent the size $h_i^2$ and the scaled indicator $\eta_i^\gamma$ associated with $P_i$. It can be noted that, even if the new centers are chosen randomly in the few selected patches, the final distribution is the expected one. In fact, most of the patches cluster around the origin, whereas the rest of the domain is covered by fewer patches. Nevertheless, we highlight that, when $C_M=9$, there are more small and medium patches far from the origin, yielding a more uniform covering of the areas far from the singular point and a slightly better accuracy. 

\begin{itemize}
\item[] \hspace{-0.8cm}\textit{{Strategy \#2: Fixed patch centers}}
\end{itemize}

From the results discussed in \textit{Strategy \#1}, it can be observed that choosing the position of the new centers randomly may lead to non-uniform patch  distribution in regions far from the singular point. In order to obtain better distributions, let us fix a priori the position of the new centers. Let us consider the reference patch $\hat P=(0,1)^2$ and the points
\begin{equation}\label{eq:c4_centers}
\begin{aligned}
\bm{\hat c_{1}}=(0.25,0.25), \hspace{1cm} \bm{\hat c_{2}}=(0.75,0.25),\\
\bm{\hat c_{3}}=(0.25,0.75), \hspace{1cm} \bm{\hat c_{4}}=(0.75,0.75),
\end{aligned}
\end{equation}
when $C_M=4$ and 
\begin{equation}\label{eq:c9_centers}
\begin{aligned}
\bm{\hat c_{1}}=(0.2,0.2), \hspace{1cm} \bm{\hat c_{2}}=(0.2,0.5), \hspace{1cm} \bm{\hat c_{3}}=(0.2,0.8),\\
\bm{\hat c_{4}}=(0.5,0.2), \hspace{1cm} \bm{\hat c_{5}}=(0.5,0.5), \hspace{1cm} \bm{\hat c_{6}}=(0.5,0.8),\\
\bm{\hat c_{7}}=(0.8,0.2), \hspace{1cm} \bm{\hat c_{8}}=(0.8,0.5), \hspace{1cm} \bm{\hat c_{9}}=(0.8,0.8),
\end{aligned}
\end{equation}
when $C_M=9$. At the end of the $(m-1)$-th training iteration, if {$\eta_i^\gamma\in\overline{\bm \eta}_{m-1}$}, the $C_M$ centers inside $P_i$ are chosen as {$\bm{c_{P_i^k}}$}{$~= M_i(\bm{\hat c_{k}})$}, $k=1,\dots,C_M$. Once more, to avoid patches partially outside $\Omega$, we update such centers as in \eqref{eq:centers_update}. We highlight that defining the new centers as in \eqref{eq:c4_centers} and in \eqref{eq:c9_centers} and the length $h_i^k$ of the edges of the new patches as in \textit{Strategy \#1}, then the new patches with centers inside $P_i$ form a cover of $P_i$, i.e., $P_i\subsetneq\cup_{k=1}^{C_M}P_i^k$. Such a property does not hold if the new centers are randomly chosen.

\begin{figure}[t!]
\centering 
  \includegraphics[width=0.6\linewidth]{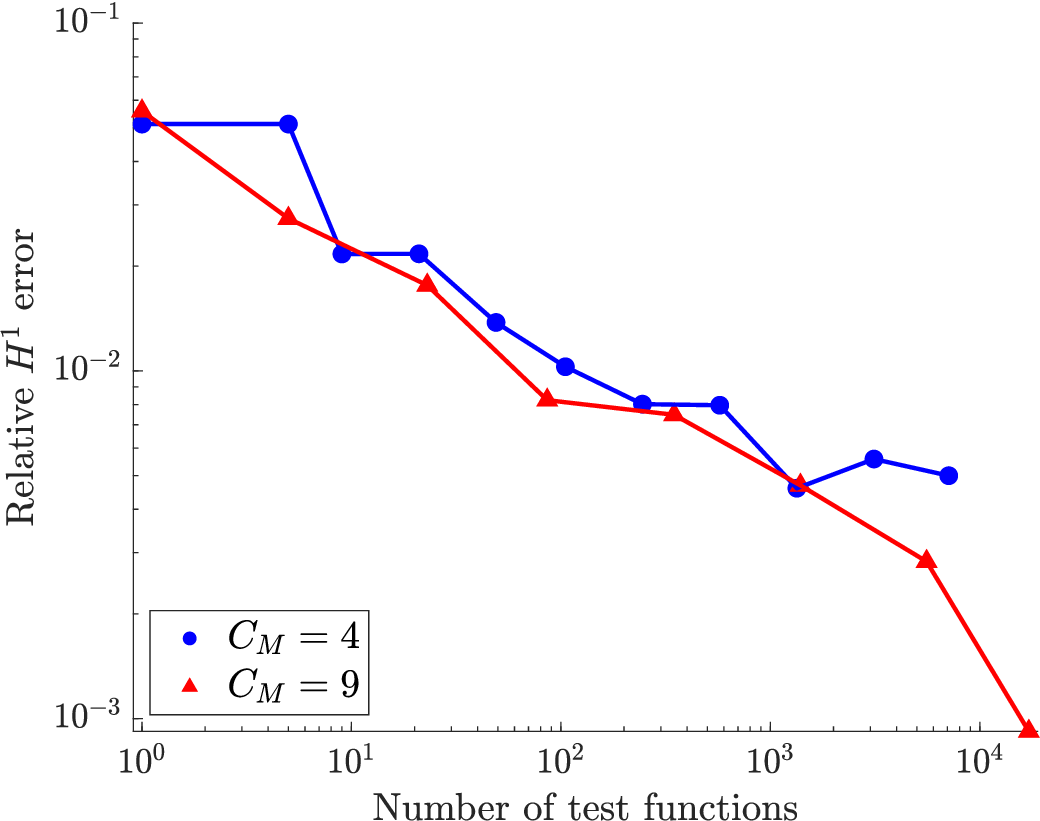} 
  \caption{\textit{{Strategy \#2}
}: Relative $H^1$ errors obtained at the end of each training iteration for $C_M=4$ (blue circles) and $C_M=9$ (red triangles).}
  \label{fig:error_struct_49_energy75}
\end{figure}

\begin{figure}[t!]
 \centering
\begin{subfigure}[H]{0.325\linewidth}
     \includegraphics[width=0.99\columnwidth,keepaspectratio,clip]{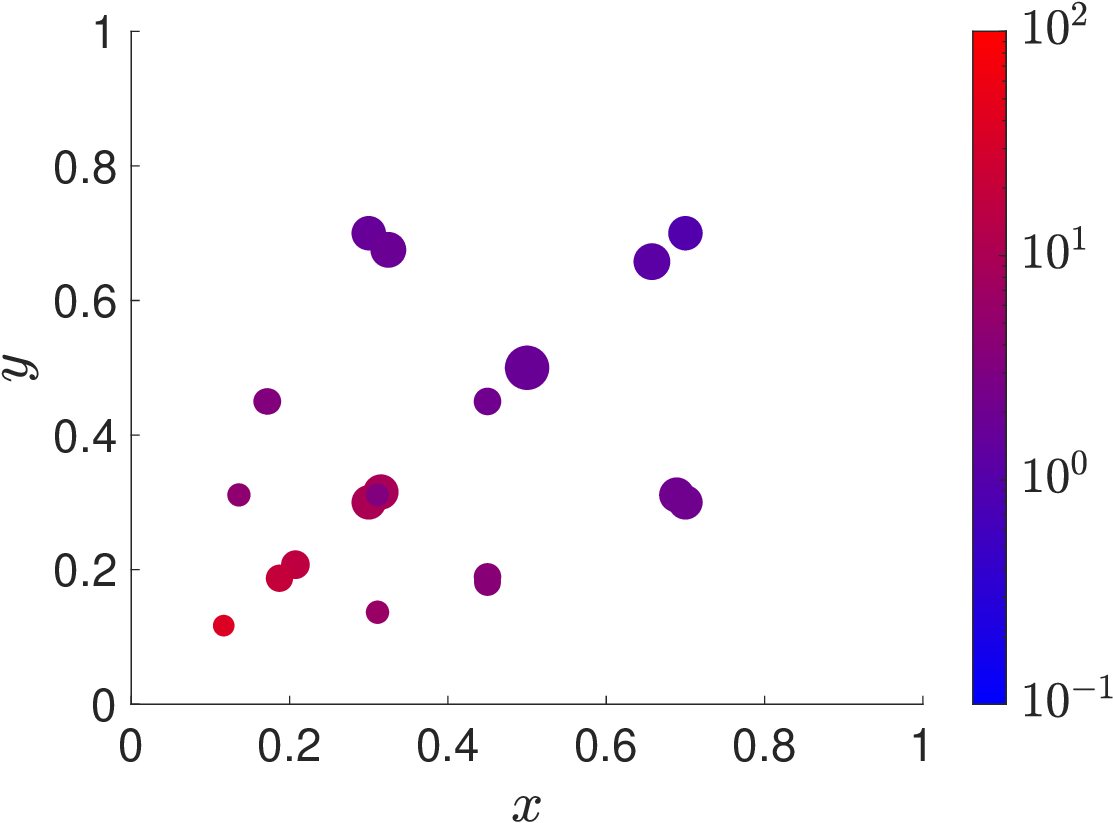} 
   \subcaption{\centering}
\end{subfigure}
\begin{subfigure}[H]{0.325\linewidth}
    \includegraphics[width=0.99\columnwidth,keepaspectratio,clip]{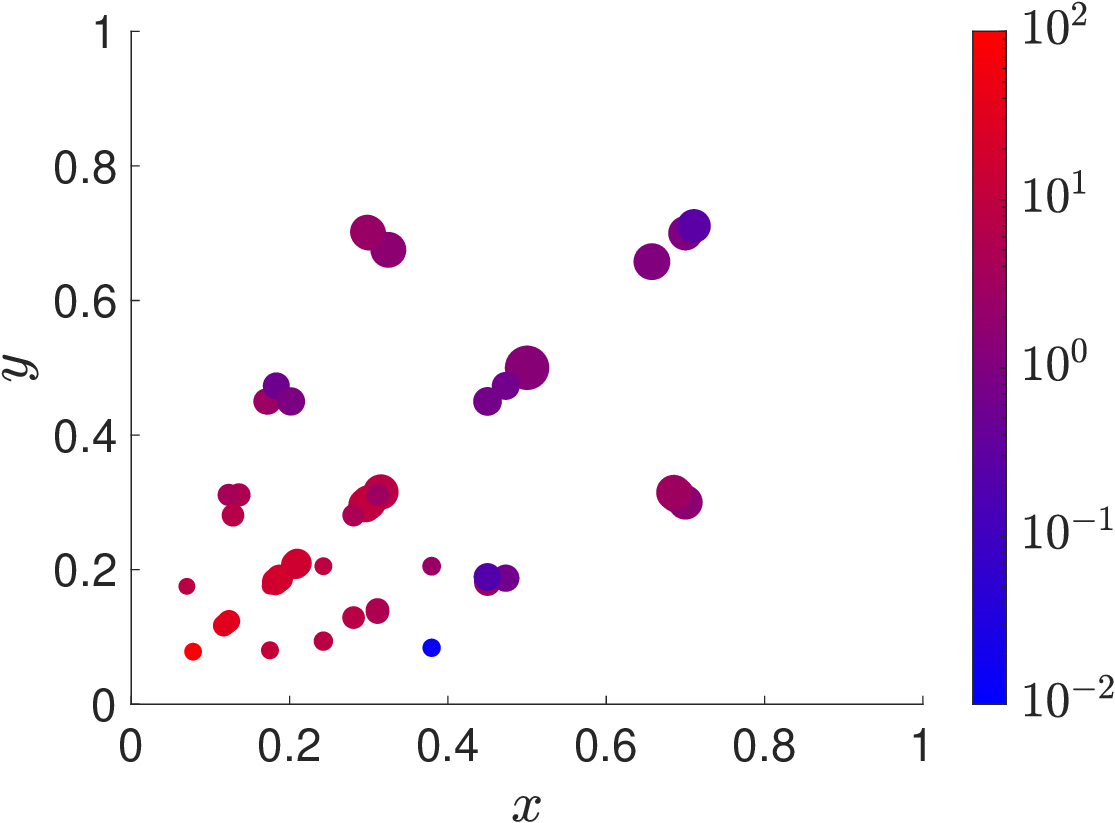} 
   \subcaption{\centering}
\end{subfigure}
\begin{subfigure}[H]{0.325\linewidth}
    \includegraphics[width=0.99\columnwidth,keepaspectratio,clip]{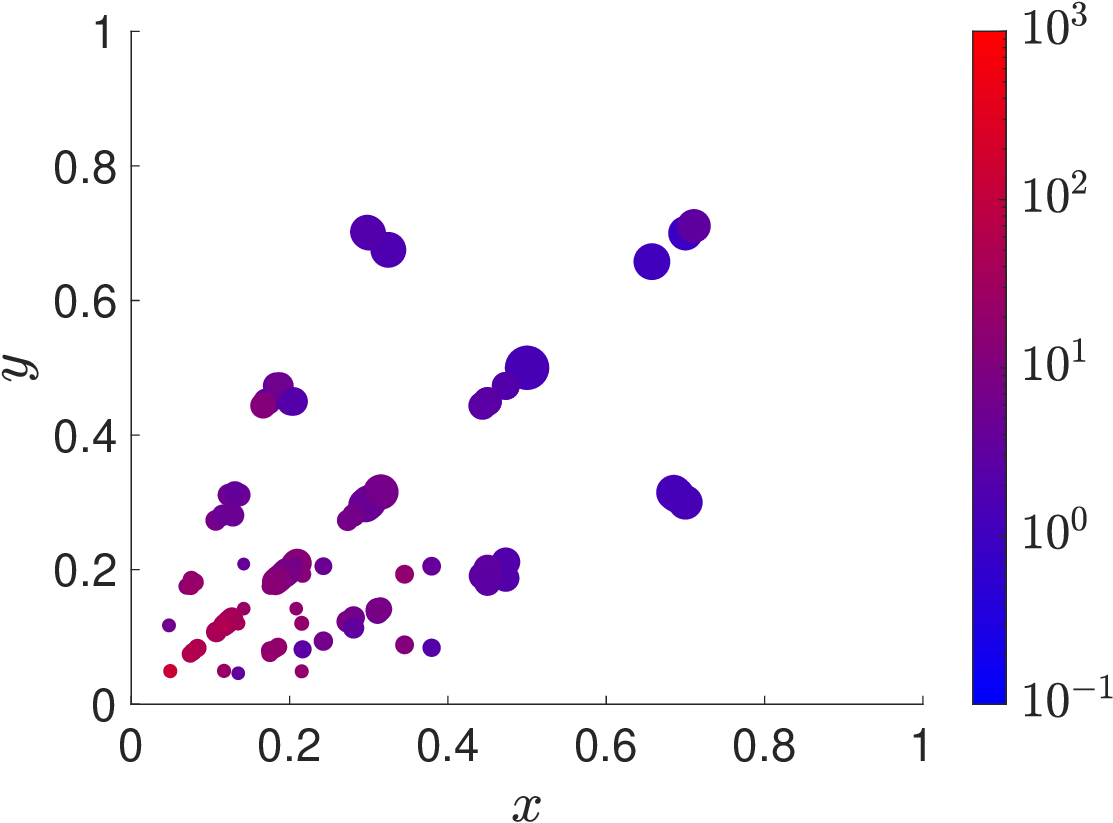}
   \subcaption{\centering}
\end{subfigure}

\medskip

\begin{subfigure}[H]{0.325\linewidth}
    \includegraphics[width=0.99\columnwidth,keepaspectratio,clip]{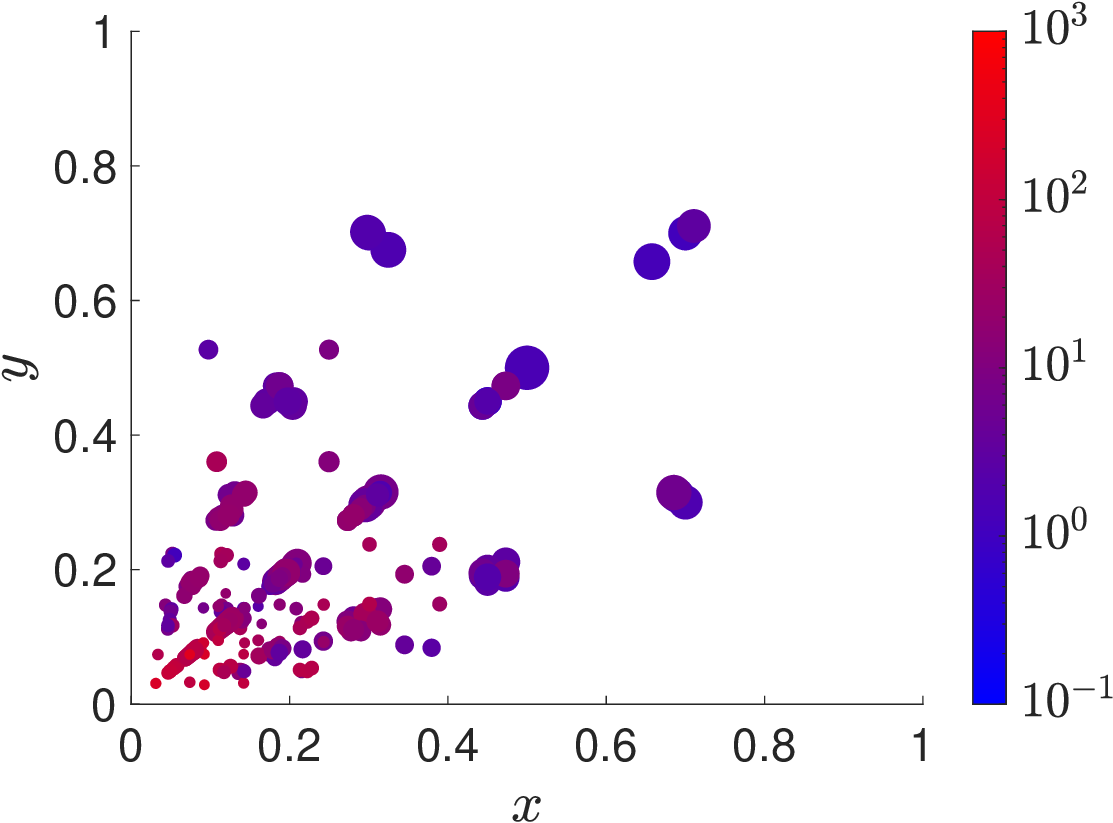} 
   \subcaption{\centering}
\end{subfigure}
\begin{subfigure}[H]{0.325\linewidth}
    \includegraphics[width=0.99\columnwidth,keepaspectratio,clip]{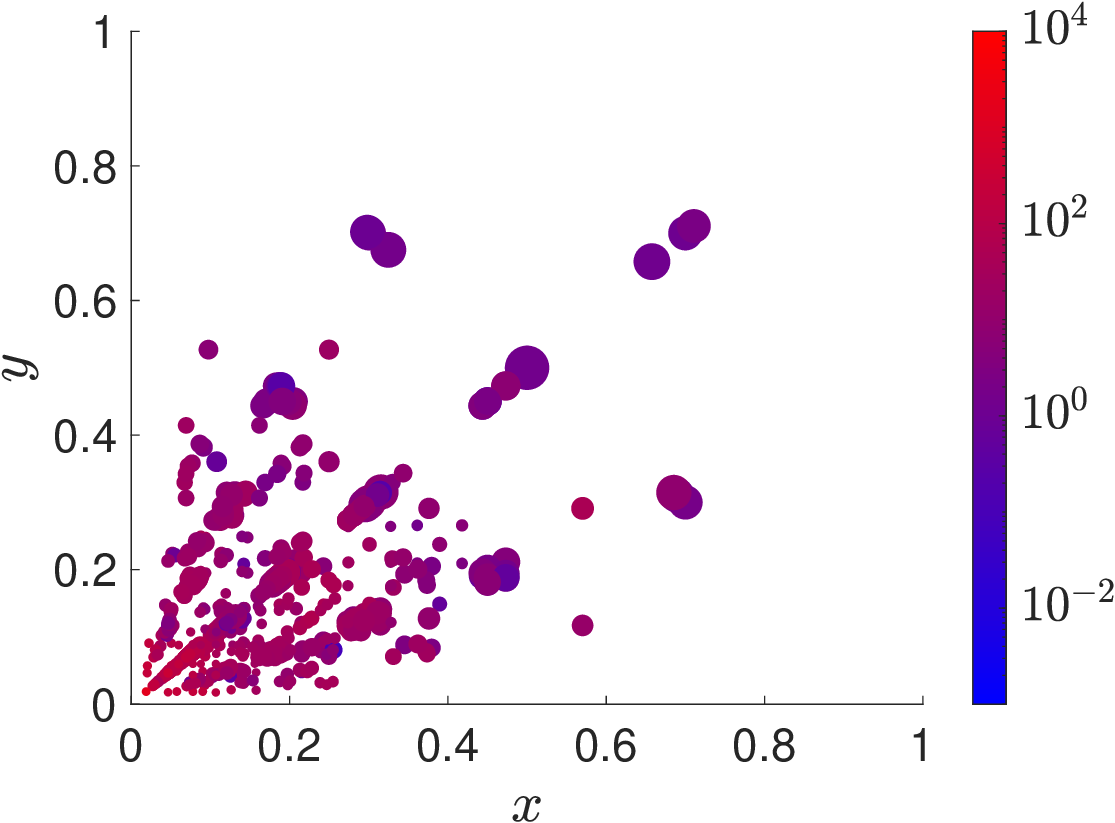} 
   \subcaption{\centering}
\end{subfigure}
\begin{subfigure}[H]{0.325\linewidth}
    \includegraphics[width=0.99\columnwidth,keepaspectratio,clip]{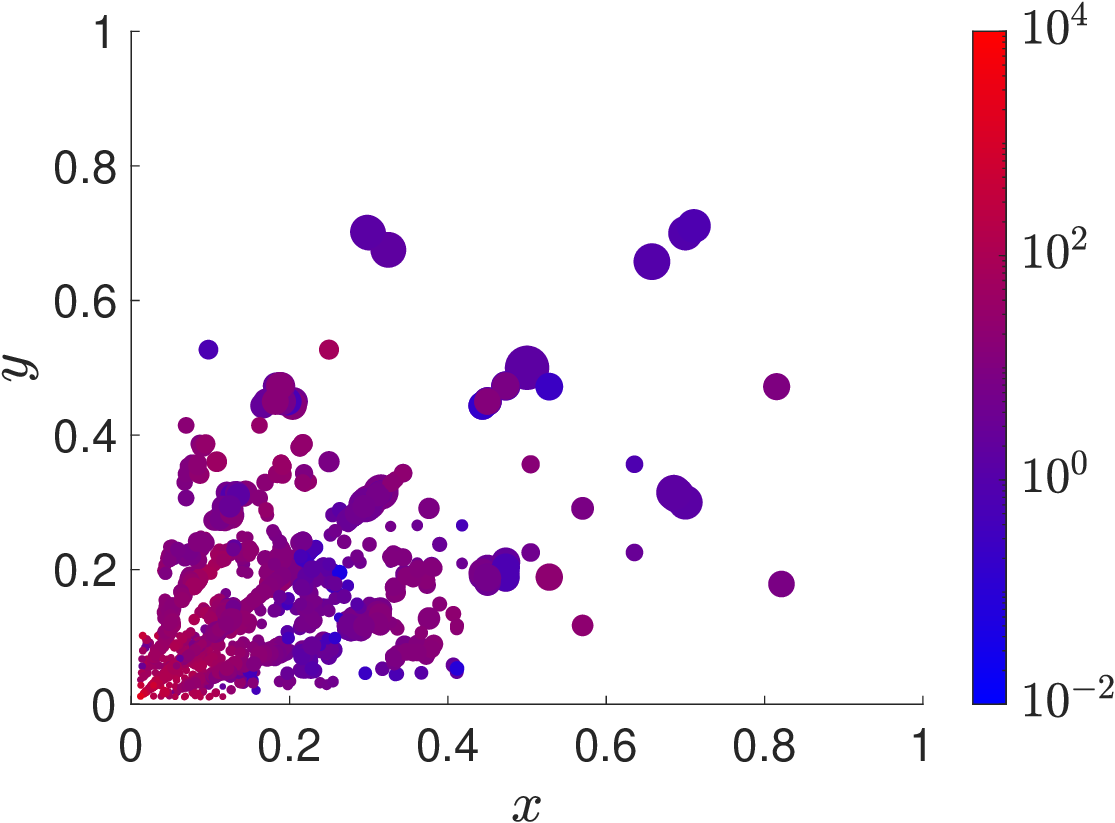}
   \subcaption{\centering}
\end{subfigure}

\medskip
  \caption{\textit{{Strategy \#2}
}: Patches used to train the MF-VPINN with $C_M=4$. Each dot represents a patch $P_i$, its position is the center $\bm{c_{P_i}}$ of the patch, its size is proportional to the patch size $h_i^2$, and its color is associated with the quantity $\eta_i^\gamma$. {(\textbf{a}) Representation of ${\cal P}_3$; (\textbf{b}) Representation of ${\cal P}_4$; (\textbf{c}) Representation of ${\cal P}_5$; (\textbf{d}) Representation of ${\cal P}_6$; (\textbf{e}) Representation of ${\cal P}_7$; (\textbf{f}) Representation of ${\cal P}_8$}. 
}
  \label{fig:struct_4_energy75_estimators}
\end{figure}

\begin{figure}[t!]
 \centering
\begin{subfigure}[H]{0.325\linewidth}
  \includegraphics[width=0.99\columnwidth,keepaspectratio,clip]{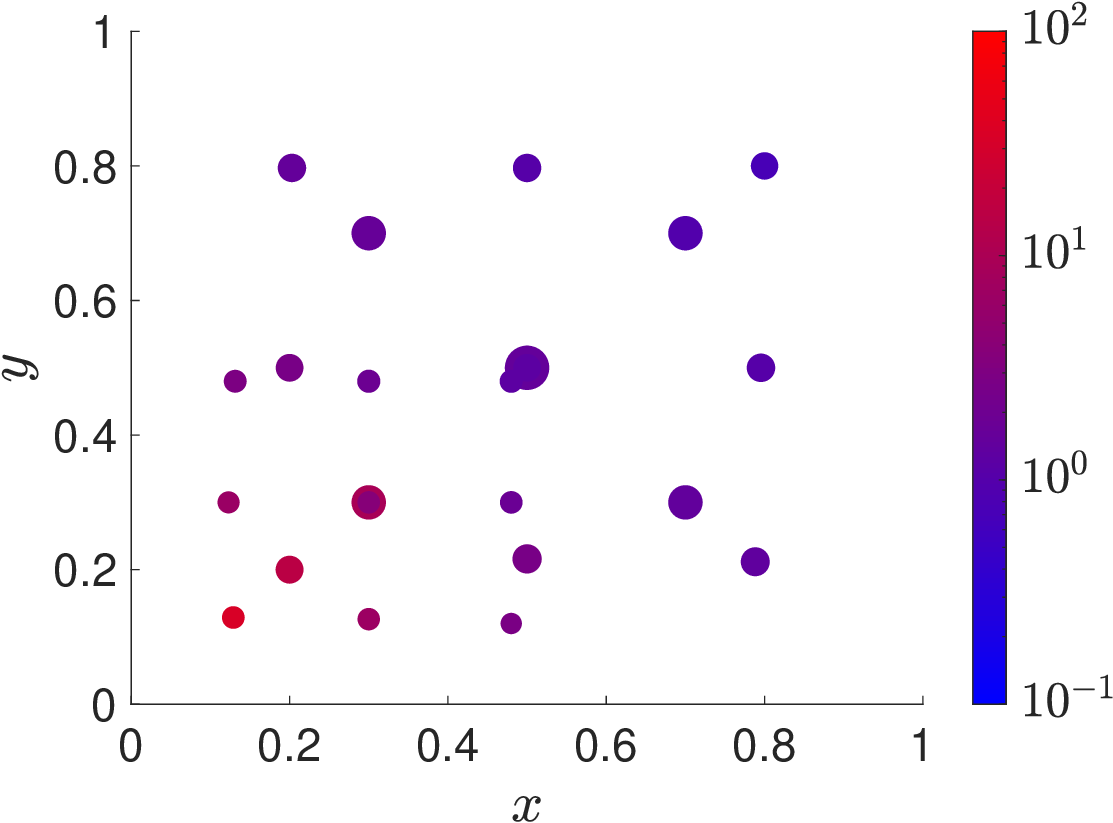} 
   \subcaption{\centering}
\end{subfigure}
\begin{subfigure}[H]{0.325\linewidth}
  \includegraphics[width=0.99\columnwidth,keepaspectratio,clip]{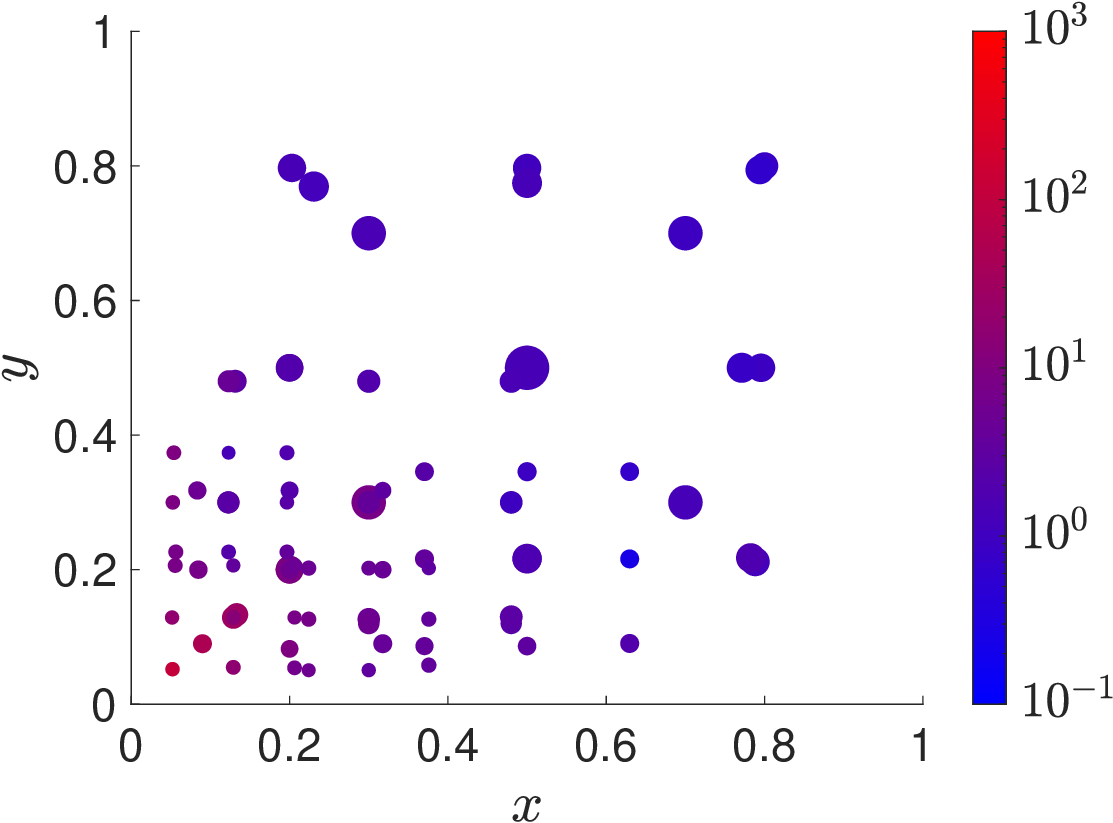} 
   \subcaption{\centering}
\end{subfigure}
\begin{subfigure}[H]{0.325\linewidth}
  \includegraphics[width=0.99\columnwidth,keepaspectratio,clip]{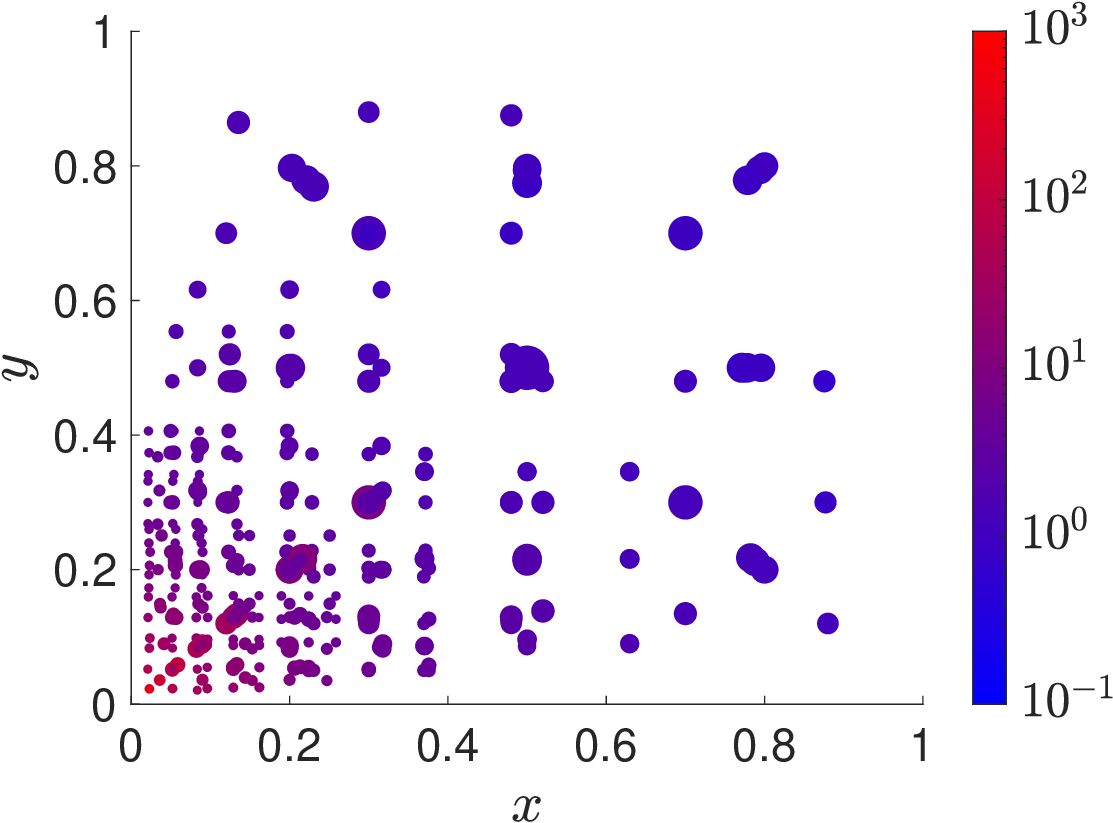}
   \subcaption{\centering}
\end{subfigure}

\medskip

\begin{subfigure}[H]{0.325\linewidth}
  \includegraphics[width=0.99\columnwidth,keepaspectratio,clip]{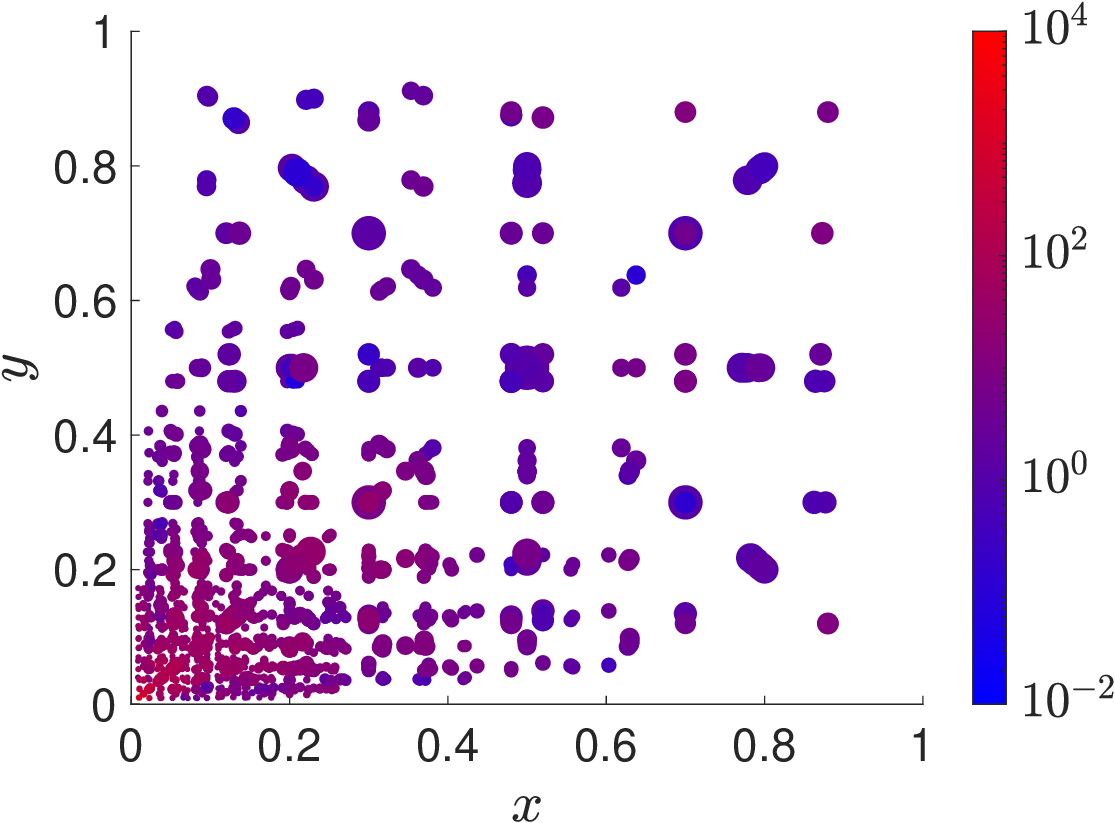} 
   \subcaption{\centering}
\end{subfigure}
\begin{subfigure}[H]{0.325\linewidth}
  \includegraphics[width=0.99\columnwidth,keepaspectratio,clip]{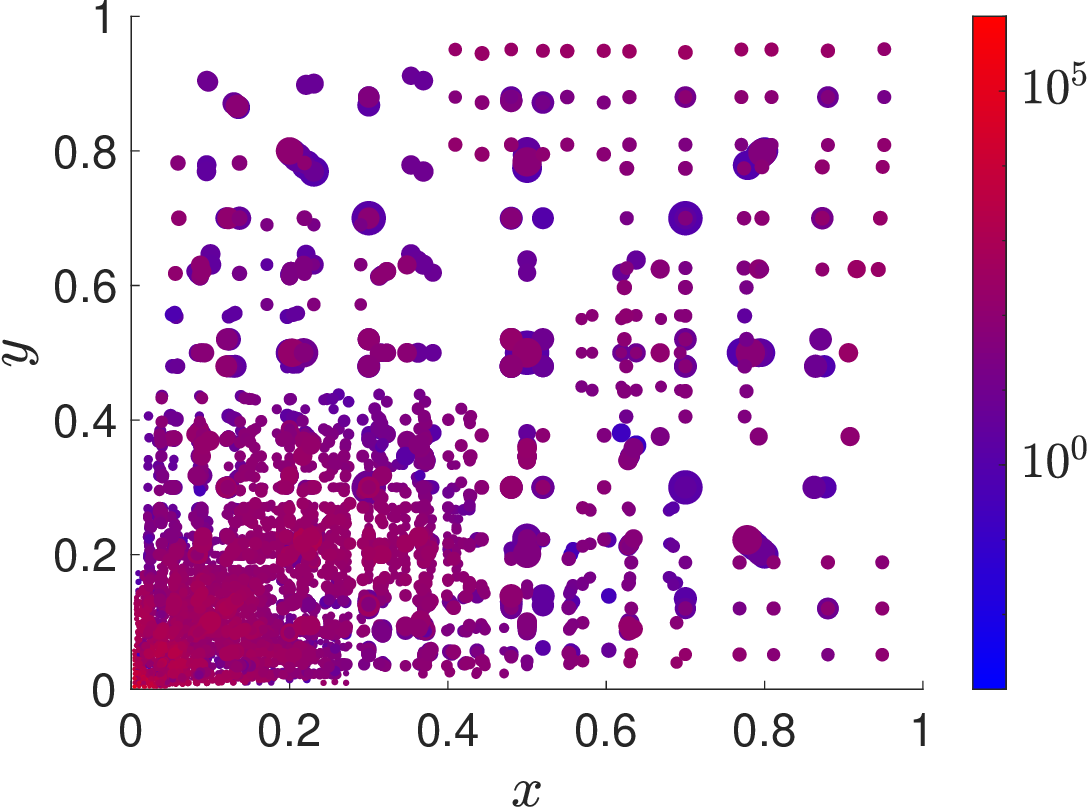} 
   \subcaption{\centering}
\end{subfigure}
\begin{subfigure}[H]{0.325\linewidth}
  \includegraphics[width=0.99\columnwidth,keepaspectratio,clip]{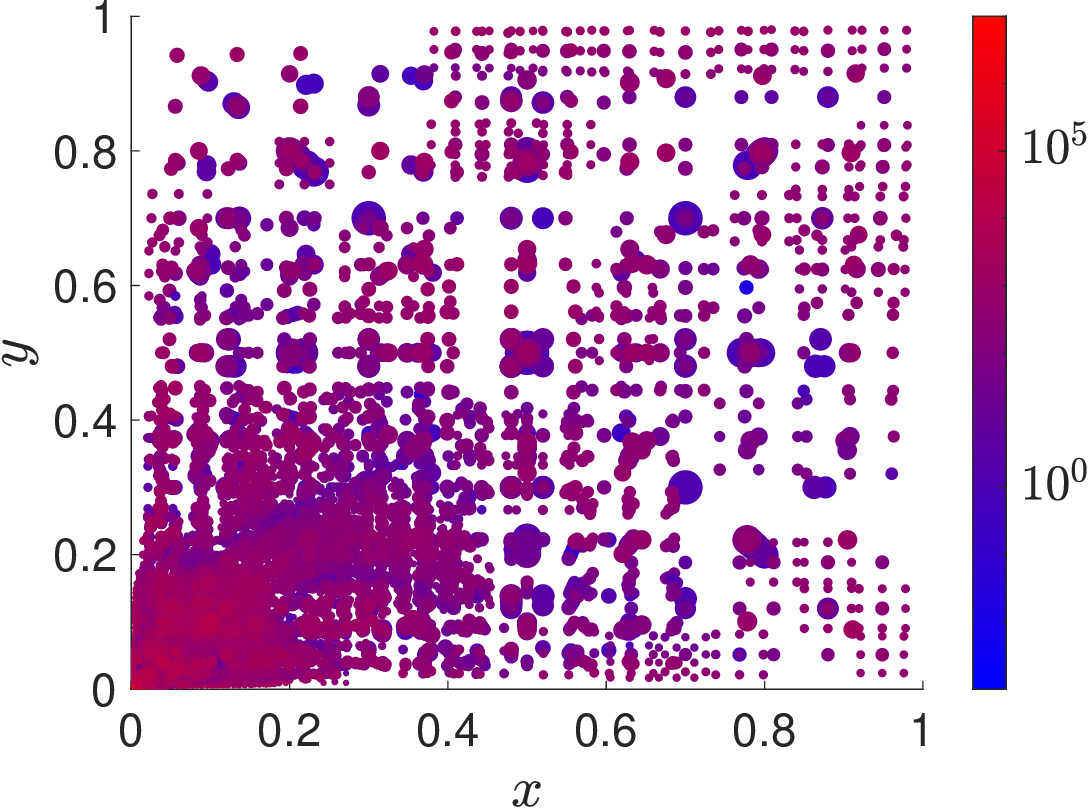}
   \subcaption{\centering}
\end{subfigure}

\medskip
  \caption{\textit{{Strategy \#2}
}: Patches used to train the MF-VPINN with $C_M=9$. Each dot represents a patch $P_i$, its position is the center $\bm{c_{P_i}}$ of the patch, its size is proportional to the patch size $h_i^2$, and its color is associated with the quantity $\eta_i^\gamma$. {(\textbf{a}) Representation of ${\cal P}_2$; (\textbf{b}) Representation of ${\cal P}_3$; (\textbf{c}) Representation of ${\cal P}_4$; (\textbf{d}) Representation of ${\cal P}_5$; (\textbf{e}) Representation of ${\cal P}_6$; (\textbf{f}) Representation of ${\cal P}_7$}. 
}
  \label{fig:struct_9_energy75_estimators}
\end{figure}

Training an MF-VPINN with such a strategy leads to more accurate results. The error decays are shown in Figure \ref{fig:error_struct_49_energy75}, whereas a comparison with the previous one will be presented in Section \ref{sec:estimator_role}. The patch distributions, for $C_M=4$ and $C_M=9$, are shown in \mbox{\cref{fig:struct_4_energy75_estimators,fig:struct_9_energy75_estimators}}, respectively. Analyzing such distributions, it can be noted that the patches still accumulate near the origin as expected. However, it is possible to observe that there are regions that are only covered by the largest patches. This phenomenon is more evident when $C_M=4$. To avoid such a phenomenon, we aim at inserting more patches far from the origin in order to train the MF-VPINN in the entire domain with a more balanced set \mbox{of patches}.

\begin{itemize}
\item[] \hspace{-0.8cm}\textit{{Strategy \#3: Fixed patch centers and small level gap strategy}}
\end{itemize}

In order to ensure better patch distributions, let us consider a new criterion to choose the position and the size of the new patches. We name this strategy the  \textit{{small-level gap strategy}} because it penalizes patch  distributions with large differences between the levels of the smallest patches and the ones of the largest patches.  

We denote by $k$-th level patch any patch $P_i$ such that $P_i\in{\cal P}_k$ and $P_i\notin{\cal P}_{k'}$ for any $k'<k$. With this notation, it is possible to group all the patches according to their level. To do so, we denote by $L_\ell$ the set of $k$-th level patches with $k\le \ell$. Let us consider the $m$-th training iteration. We define {${\bm \eta}_\text{sort}^\ell$} as the array containing the elements $\eta_i^\gamma$ of {${\bm \eta}_\text{sort}$} (maintaining the same ordering) such that $P_i\in L_\ell$. We also denote by {$\overline{\bm \eta}_{m,\ell}$} the array containing the first $\tau_m^\ell=\min\{\tau_m,\text{dim}(L_\ell)\}$ elements of {${\bm \eta}_\text{sort}^\ell$}. Note that {$\overline{\bm \eta}_{m,\ell}$} is the equivalent of {$\overline{\bm \eta}_{m}$} for patches in $L_\ell$.

In order to generate the new patches in ${\cal P}_{m+1}\backslash{\cal P}_m$, let us add $C_M$ new patches in any patch $P_i$ such that {$\eta_i^\gamma\in\overline{\bm \eta}_{m}\cup\overline{\bm \eta}_{m,\ell}$}. The centers and sizes of the new patches are chosen as in \textit{Strategy \#2}. This allows us to exploit the fact that $P_i\subsetneq\cup_{k=1}^{C_M}P_i^k$ to remove the patches $P_i$ such that {$\eta_i^\gamma\in\overline{\bm \eta}_{m}\cup\overline{\bm \eta}_{m,\ell}$} from the new set of patches ${\cal P}_{m+1}$. We remark that such patches cannot be removed when the centers are randomly chosen as in \textit{Strategy \#1} because, in that case, ${\cal P}_{m+1}$ would not be a cover of $\Omega$ anymore. 

We also highlight that, removing the patches $P_i$ such that {$\eta_i^\gamma\in\overline{\bm \eta}_{m}\cup\overline{\bm \eta}_{m,\ell}$} and choosing $A_\text{ratio}=1$, it is possible to satisfy the inequality
\[
\sum_{P_i\in{\cal P}_{m+1}}\vert P_i\vert \le C\vert \Omega\vert,
\]
for any $m\in\N$ and with $C>0$ independent of $m$. Such a bound on the sum of the area of the patches is useful to ensure that there exists a number $N_\text{patch\_per\_point}$ such that any point inside $\Omega$ belongs to at most $N_\text{patch\_per\_point}$ patches. This property is useful to derive global error indicators. We choose to maintain $A_\text{ratio}=1.25$ to compare the numerical results with the ones obtained using the previous strategies and to consider overlapping patches.

\begin{figure}[t!]
\centering 
  \includegraphics[width=0.6\linewidth]{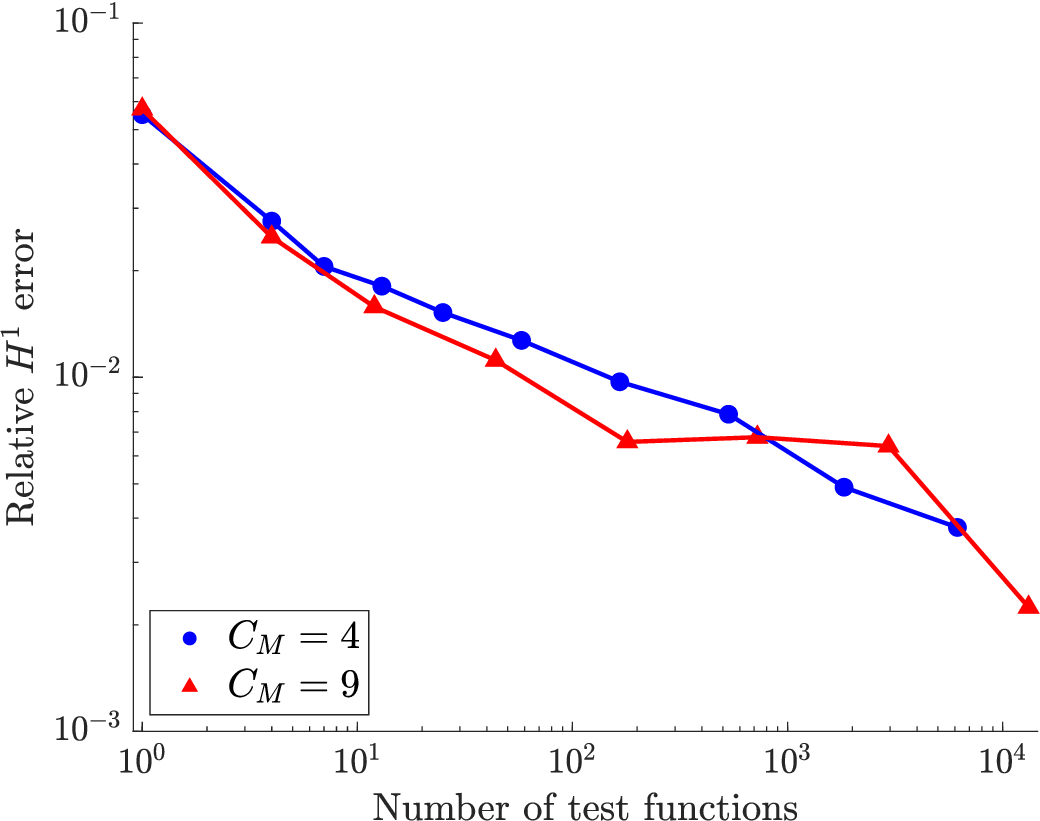} 
  \caption{\textit{{Strategy \#3}
}: Relative $H^1$ errors obtained at the end of each training iteration for $C_M=4$ (blue circles) and $C_M=9$ (red triangles).}
  \label{fig:error_struct_49_energy75_removingOld_forceSomeLevelsAndRefine}
\end{figure}

\begin{figure}[t!]
 \centering
\begin{subfigure}[H]{0.325\linewidth}
  \includegraphics[width=0.99\columnwidth,keepaspectratio,clip]{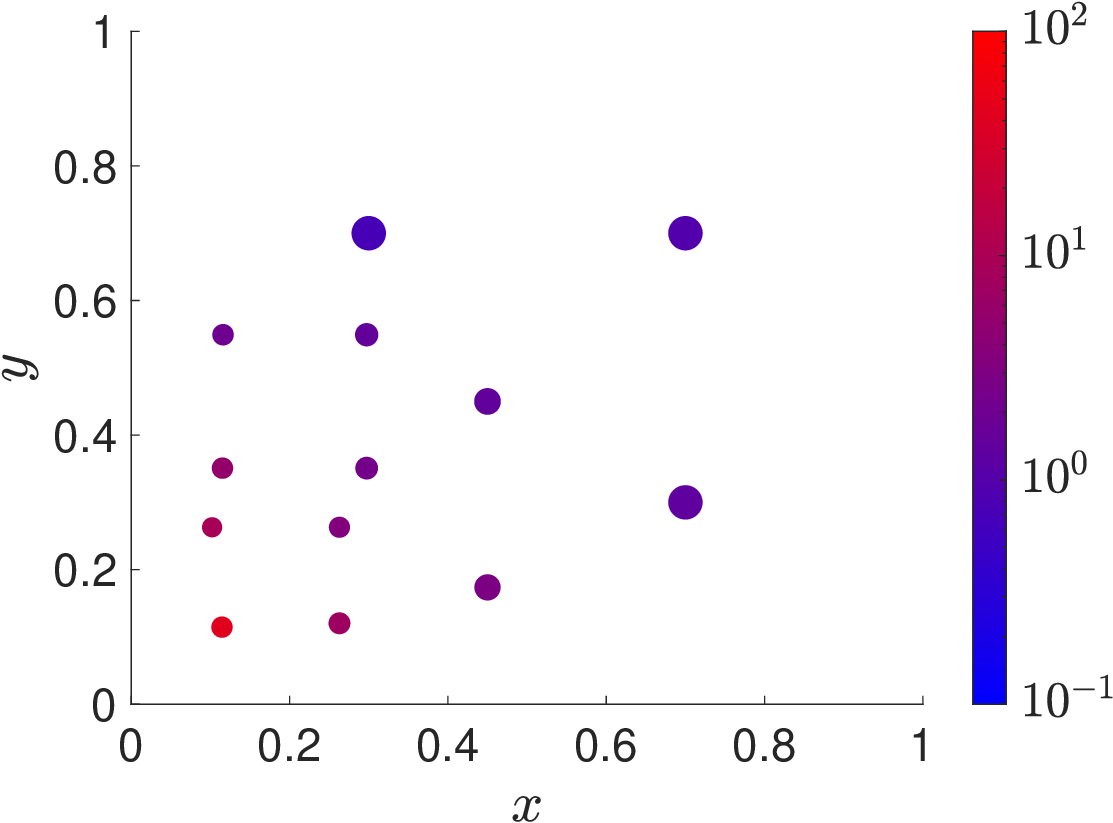} 
   \subcaption{\centering}
\end{subfigure}
\begin{subfigure}[H]{0.325\linewidth}
  \includegraphics[width=0.99\columnwidth,keepaspectratio,clip]{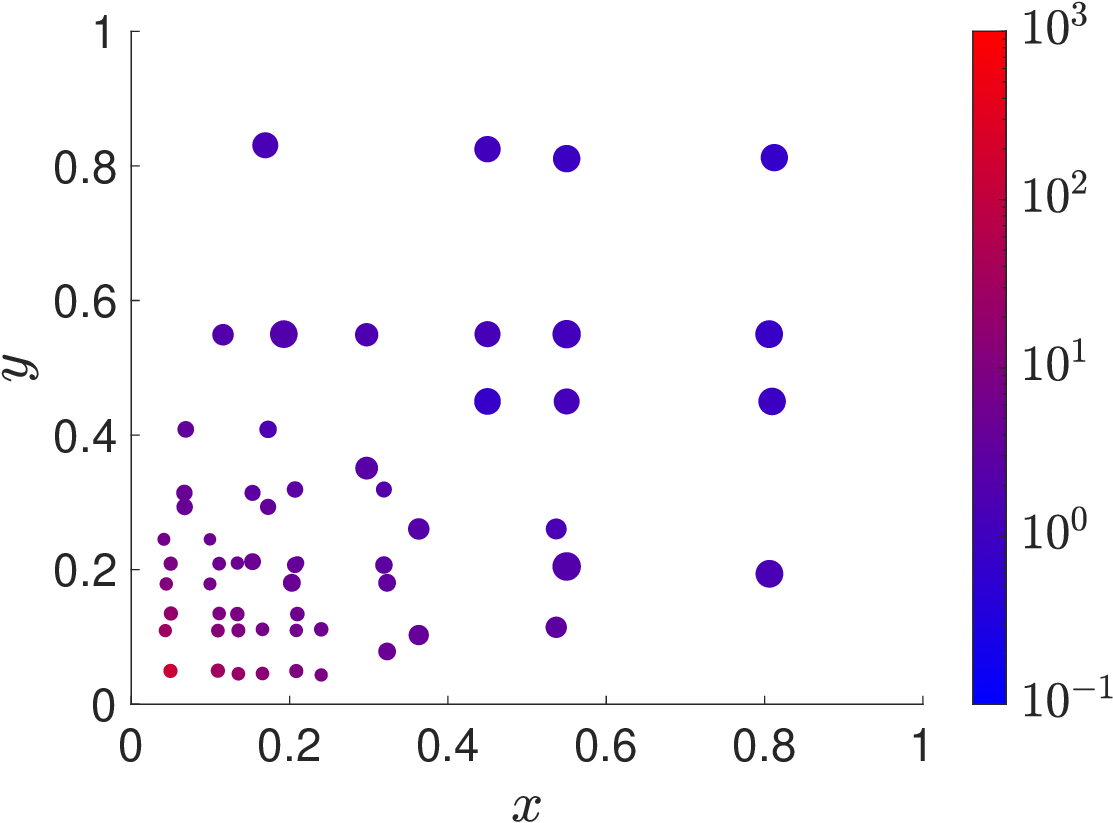} 
   \subcaption{\centering}
\end{subfigure}
\begin{subfigure}[H]{0.325\linewidth}
  \includegraphics[width=0.99\columnwidth,keepaspectratio,clip]{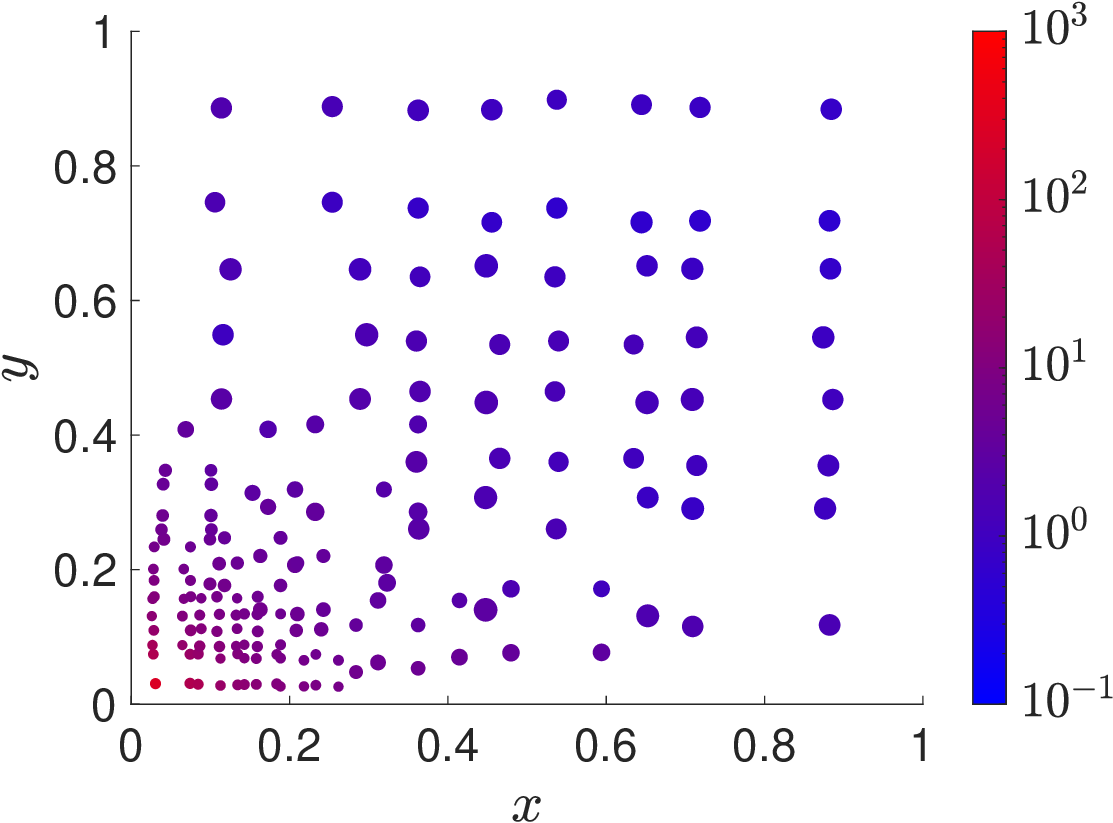} 
   \subcaption{\centering}
\end{subfigure}

\medskip

\begin{subfigure}[H]{0.325\linewidth}
  \includegraphics[width=0.99\columnwidth,keepaspectratio,clip]{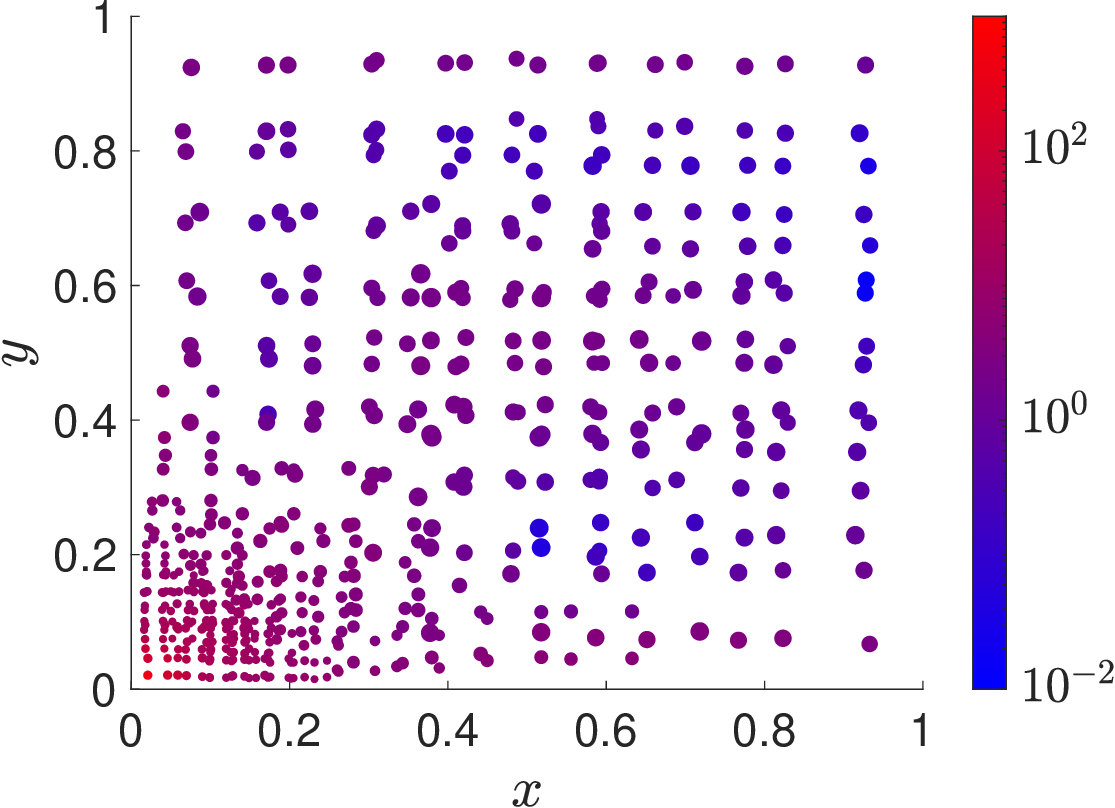} 
   \subcaption{\centering}
\end{subfigure}
\begin{subfigure}[H]{0.325\linewidth}
  \includegraphics[width=0.99\columnwidth,keepaspectratio,clip]{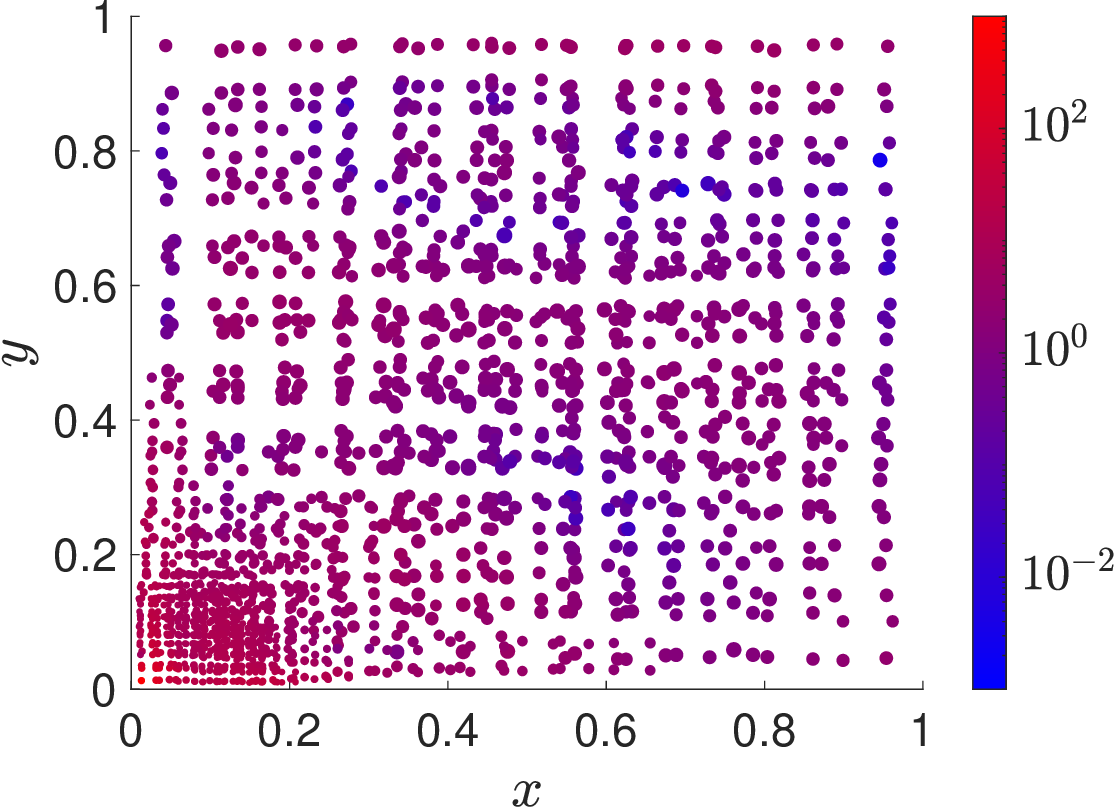} 
   \subcaption{\centering}
\end{subfigure}
\begin{subfigure}[H]{0.325\linewidth}
  \includegraphics[width=0.99\columnwidth,keepaspectratio,clip]{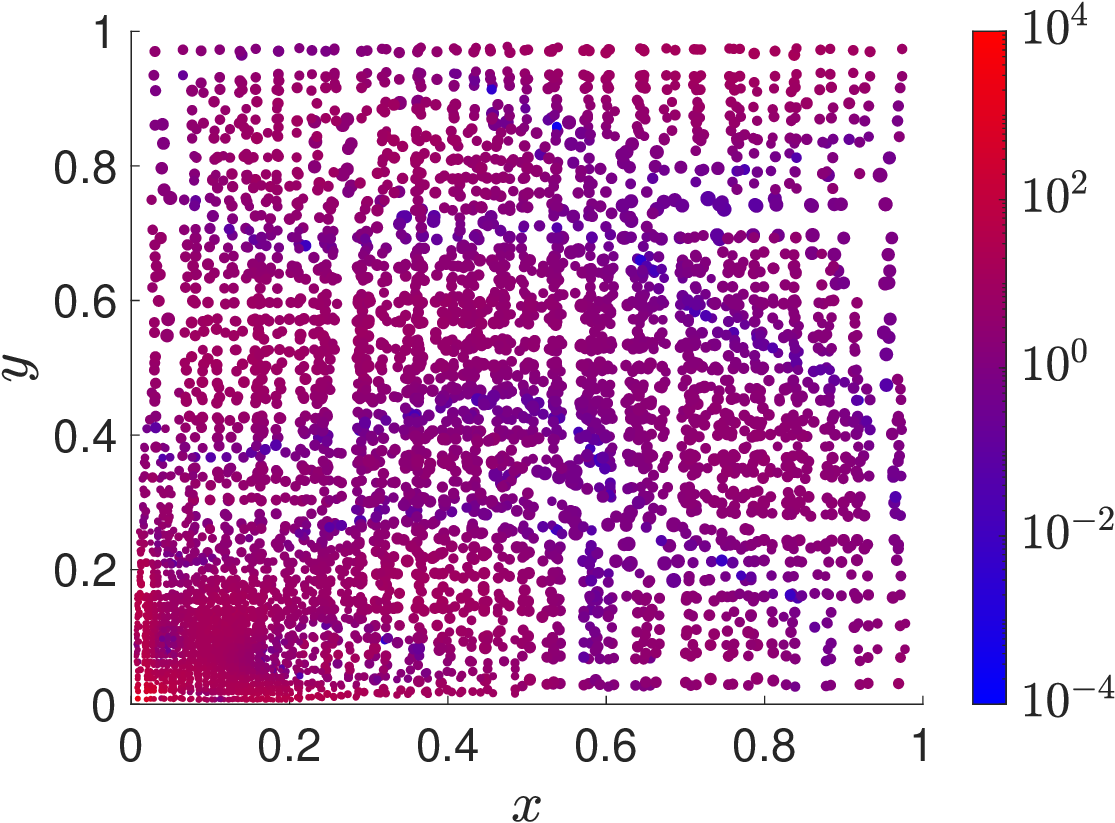} 
   \subcaption{\centering}
\end{subfigure}

\medskip
  \caption{\textit{{Strategy \#3}
}: Patches used to train the MF-VPINN with $C_M=4$. Each dot represents a patch $P_i$, its position is the center $\bm{c_{P_i}}$ of the patch, its size is proportional to the patch size $h_i^2$, and its color is associated with the quantity $\eta_i^\gamma$. {(\textbf{a}) Representation of ${\cal P}_3$; (\textbf{b}) Representation of ${\cal P}_5$; (\textbf{c}) Representation of ${\cal P}_6$; (\textbf{d}) Representation of ${\cal P}_7$; (\textbf{e}) Representation of ${\cal P}_8$; (\textbf{f}) Representation of ${\cal P}_9$}. 
}
  \label{fig:struct_4_energy75_removingOld_forceSomeLevelsAndRefine_estimators}
\end{figure}

\begin{figure}[t!]
 \centering
\begin{subfigure}[H]{0.325\linewidth}
  \includegraphics[width=0.99\columnwidth,keepaspectratio,clip]{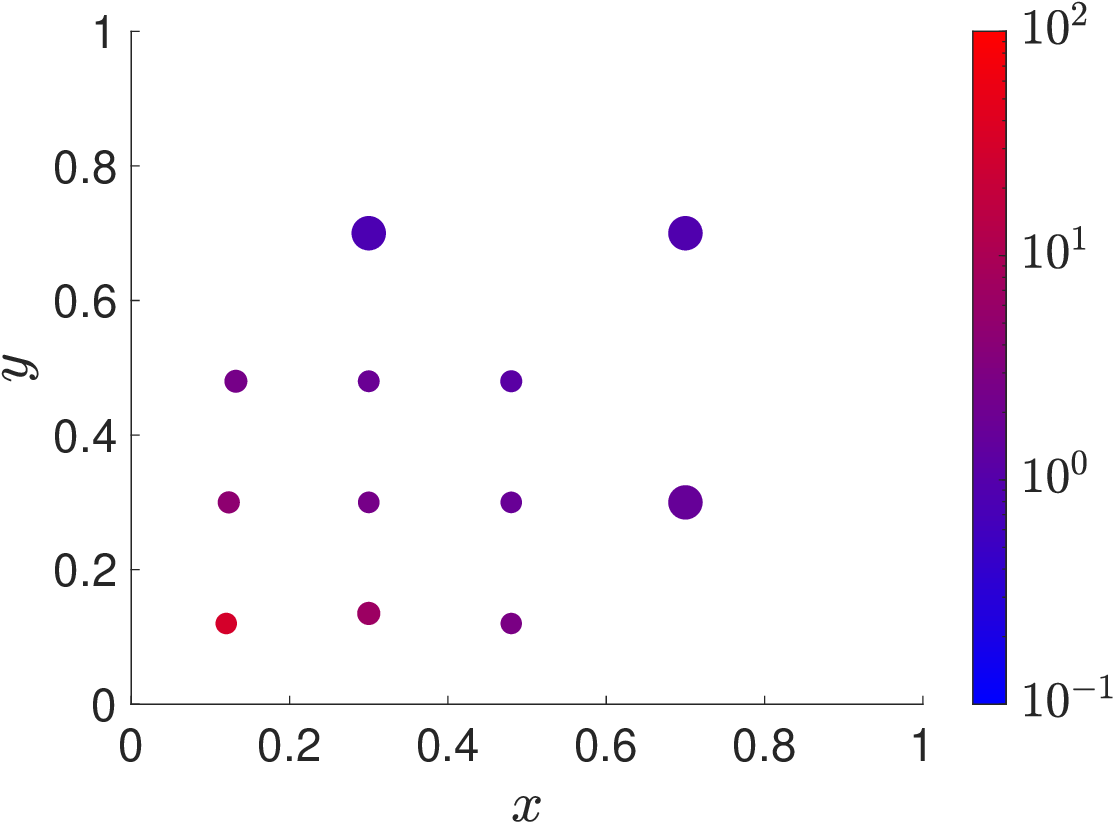}
   \subcaption{\centering}
   \label{fig:struct_9_energy75_removingOld_forceSomeLevelsAndRefine_estimators_step_1}
\end{subfigure}
\begin{subfigure}[H]{0.325\linewidth}
  \includegraphics[width=0.99\columnwidth,keepaspectratio,clip]{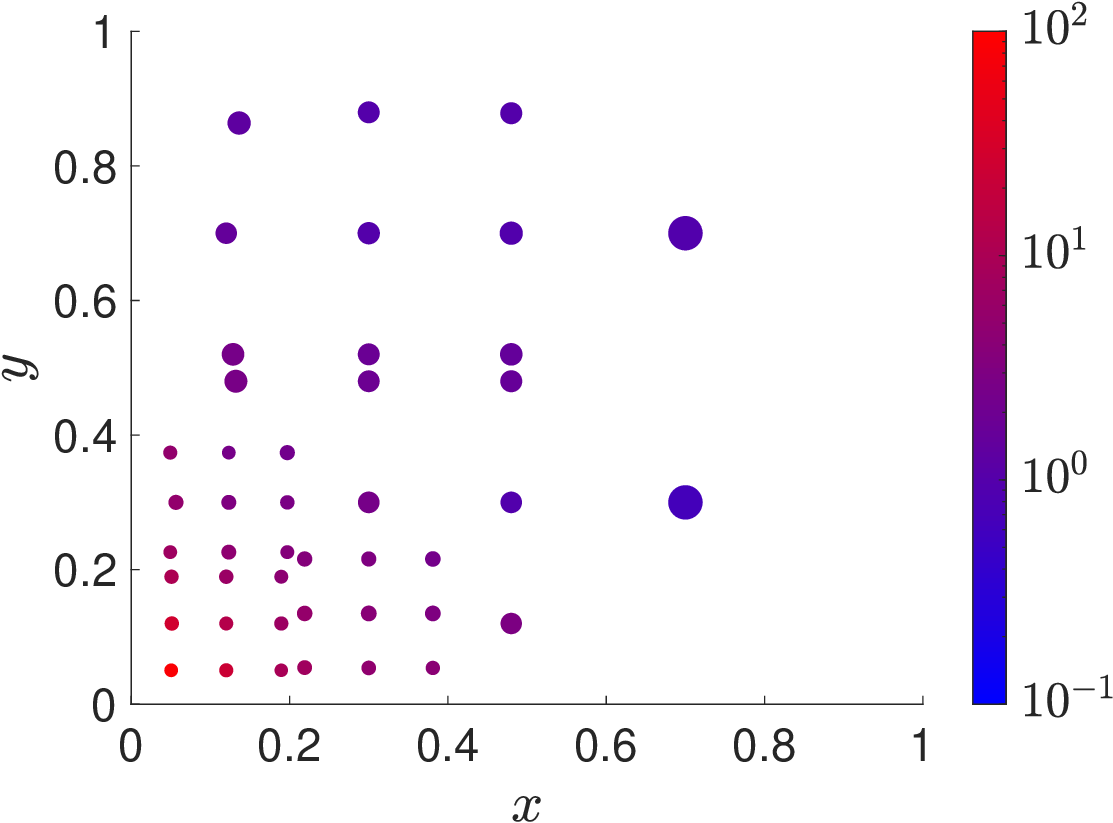}
   \subcaption{\centering}
\end{subfigure}
\begin{subfigure}[H]{0.325\linewidth}
  \includegraphics[width=0.99\columnwidth,keepaspectratio,clip]{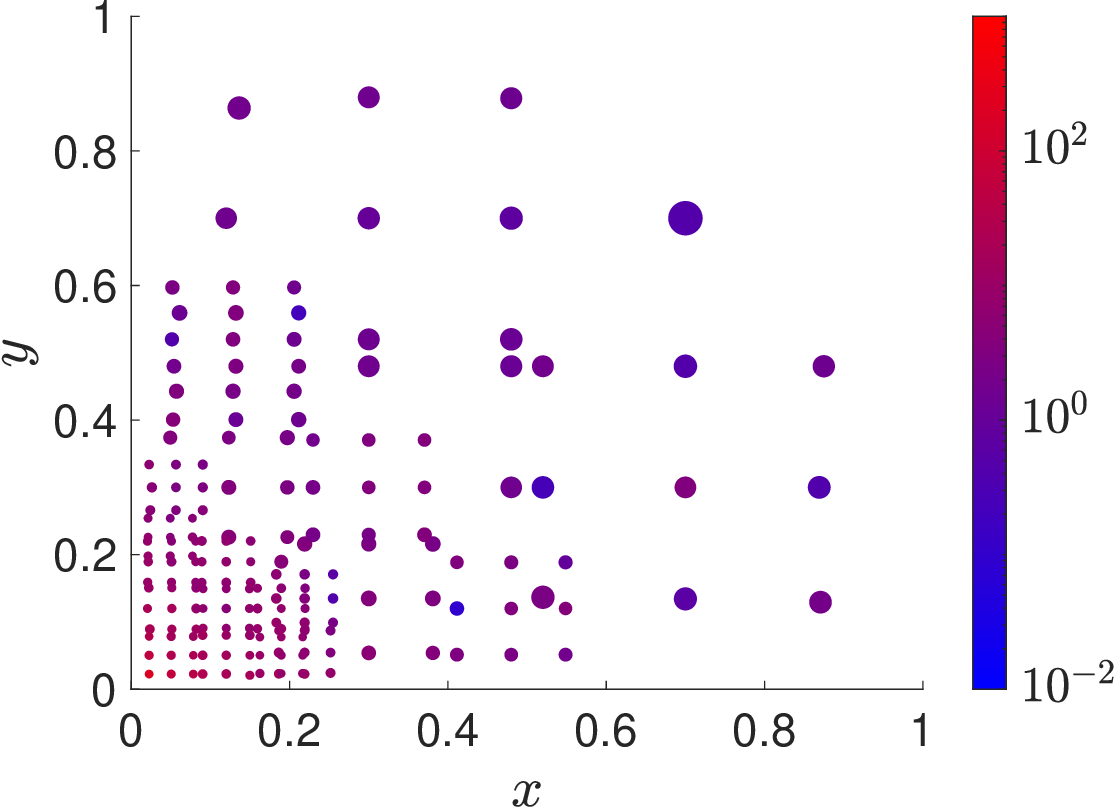} 
   \subcaption{\centering}
\end{subfigure}

\medskip

\begin{subfigure}[H]{0.325\linewidth}
  \includegraphics[width=0.99\columnwidth,keepaspectratio,clip]{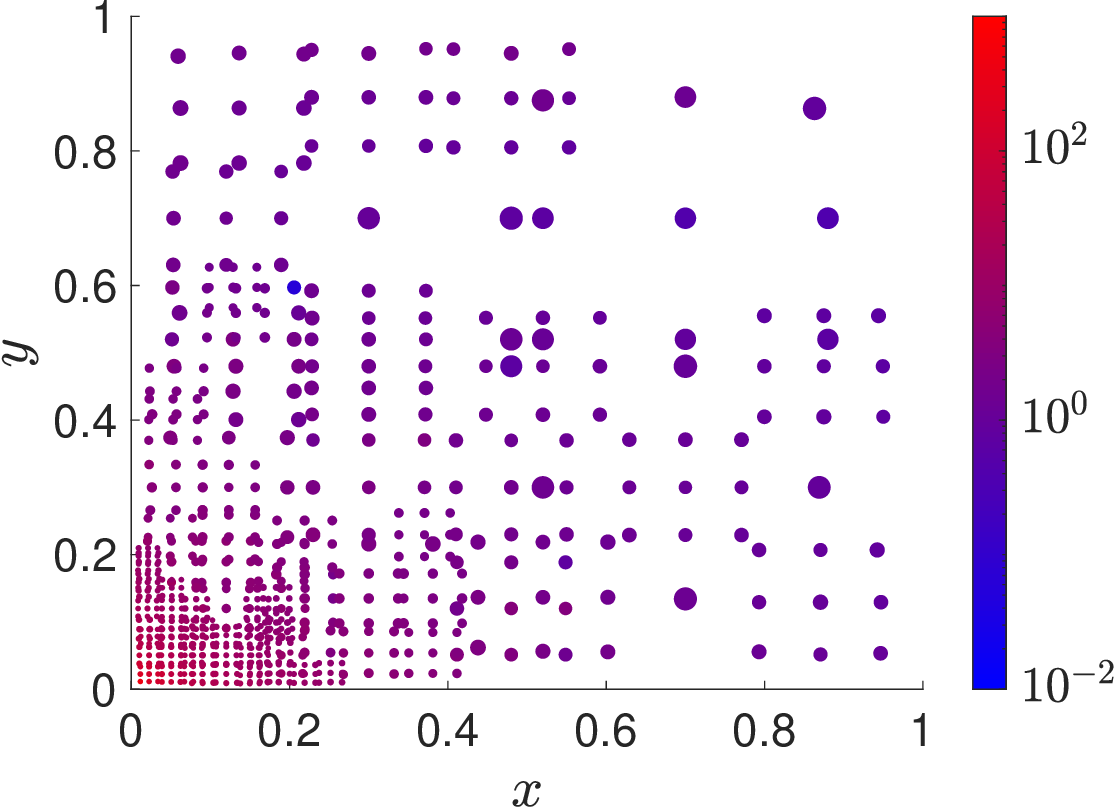}
   \subcaption{\centering}
\end{subfigure}
\begin{subfigure}[H]{0.325\linewidth}
  \includegraphics[width=0.99\columnwidth,keepaspectratio,clip]{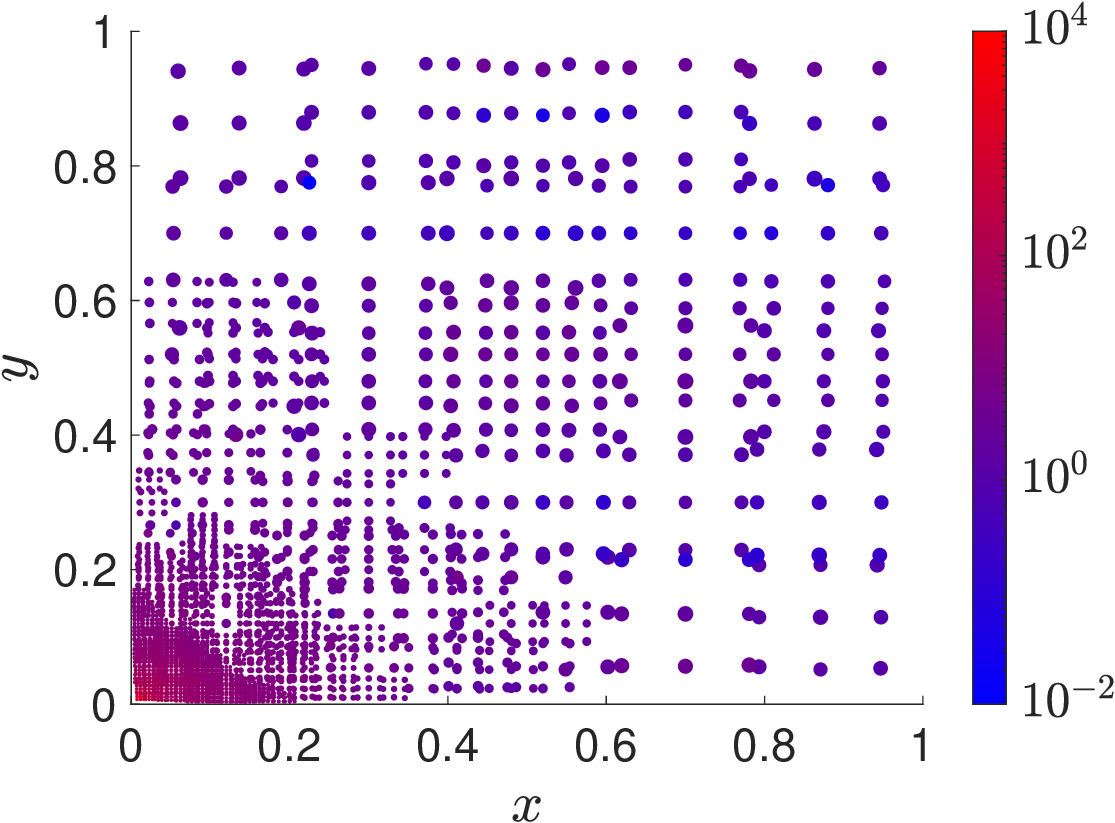}
   \subcaption{\centering}
\end{subfigure}
\begin{subfigure}[H]{0.325\linewidth}
  \includegraphics[width=0.99\columnwidth,keepaspectratio,clip]{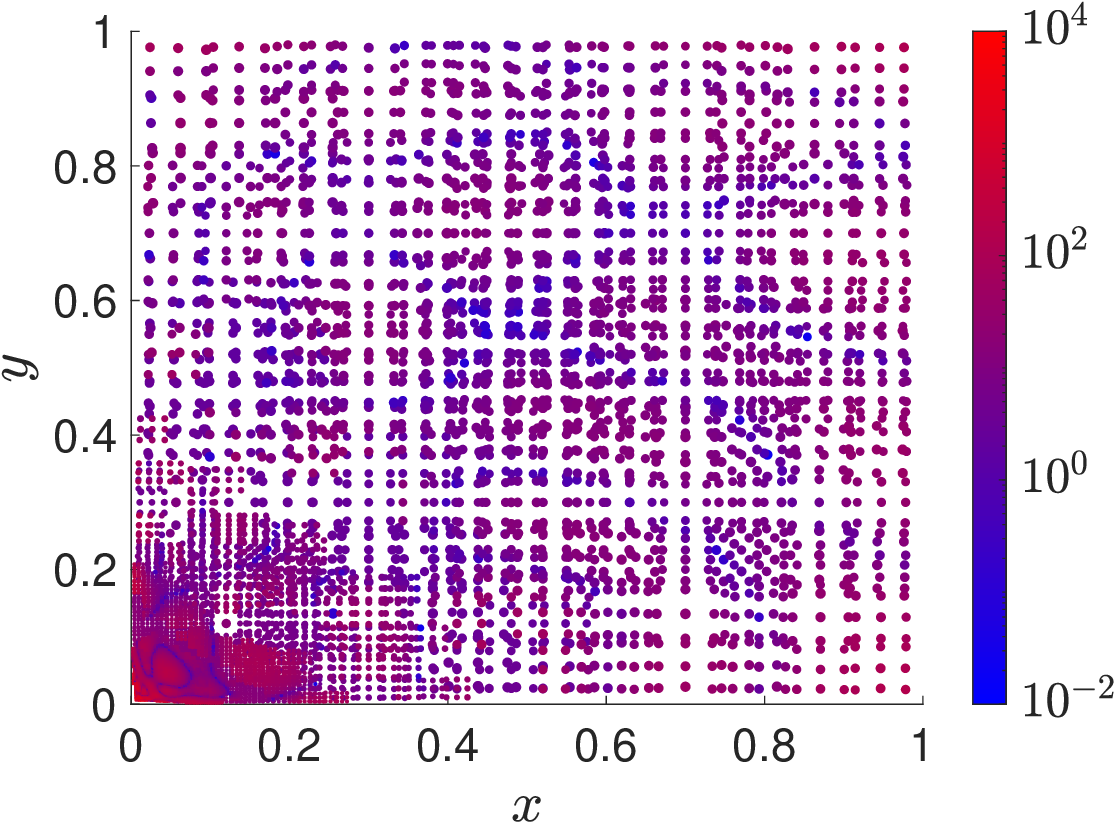}
   \subcaption{\centering}
\end{subfigure}

\medskip
  \caption{\textit{{Strategy \#3}
}: Patches used to train the MF-VPINN with $C_M=9$. Each dot represents a patch $P_i$, its position is the center $\bm{c_{P_i}}$ of the patch, its size is proportional to the patch size $h_i^2$, and its color is associated with the quantity $\eta_i^\gamma$. {(\textbf{a}) Representation of ${\cal P}_2$; (\textbf{b}) Representation of ${\cal P}_3$; (\textbf{c}) Representation of ${\cal P}_4$; (\textbf{d}) Representation of ${\cal P}_5$; (\textbf{e}) Representation of ${\cal P}_6$; (\textbf{f}) Representation of ${\cal P}_7$}. 
}
  \label{fig:struct_9_energy75_removingOld_forceSomeLevelsAndRefine_estimators}
\end{figure}

We train an MF-VPINN with $C_M=4$ and $C_M=9$ as in the previous tests. The corresponding error decays are shown in Figure \ref{fig:error_struct_49_energy75_removingOld_forceSomeLevelsAndRefine}. It can be observed that the error decreases in a smoother way and that, as in the previous tests, choosing $C_M=4$ or $C_M=9$ does not lead to significant differences in the error behavior. The patches used during the training are represented in \cref{fig:struct_4_energy75_removingOld_forceSomeLevelsAndRefine_estimators,fig:struct_9_energy75_removingOld_forceSomeLevelsAndRefine_estimators}. We highlight that, when compared with the patch  distributions in \textit{Strategy \#2}, there exist much more patches far from the origin, and, most importantly, the closer the center of a patch to the origin, the smaller its size. Even though the error decays with $C_M=4$ and $C_M=9$ are qualitatively similar, it should be noted that the patch  distribution with $C_M=9$ is more skewed. In fact. its patches can be clustered into two subgroups: the first one containing larger patches and covering most of the domain  the second one containing only small patches with centers very close to the origin. A similar distribution is obtained with $C_M=4$, even though it is characterized by a smoother transition between large and small patches.

In both cases, it can be observed that there are no large patches very close to small ones. This is in contrast with the distributions obtained in \textit{Strategy \#2} and leads to more stable solvers. Indeed, even though the test functions are not related to a global triangulation on the entire domain $\Omega$, the current loss function is very similar to the one used in a standard VPINN with a good-quality mesh, i.e., a mesh in which neighboring elements are similar in size and shape. On the other hand, in \textit{Strategy \#2}, there exist large patches that are very close to small ones; this is equivalent to training a VPINN on a very poor-quality mesh. Such meshes, in the context of FEM, are strictly related to convergence and accuracy issues.

\subsection{The Importance of the Error Indicator}\label{sec:estimator_role}
{As discussed in the previous sections, we use the error indicator described in \mbox{Section \ref{sec:aposteriori}} to interrupt the training and to decide where the new patches have to be inserted to maximize the accuracy. In this section, the advantages of such a choice are described.}

Since each set ${\cal P}_m$ is a cover of $\Omega$, the quantity $ES_m=\sum_{i=1}^{\text{dim}({\cal P}_m)}\eta_i$ is an indicator of the global $H^1$ error $\Vert u-u^\NN\Vert_1$ on the entire domain $\Omega$. Therefore, tracking its behavior during the training is equivalent to tracking that of the unknown $H^1$ error. Such information is used to implement an early stopping strategy to reduce the computational cost of the iterative training. At the beginning of the $m$-th training iteration, all the vectors and sparse matrices required to compute $ES_m$ are computed in a preprocessing phase. When such data structures are available, the error indicator can be assembled suitably combining basic algebraic operations. 

We assemble $ES_m$ every $N_\text{check}$ epochs and store the best value obtained during the training, together with the corresponding neural network trainable parameters. Then, if no improvements are obtained in $p\cdot N_\text{check}$ epochs, the training is interrupted and the neural network parameters associated with the best value of $ES_m$ are restored. Here, $p$ is a tunable parameter named \textit{patience}. The first $N_\text{negl}^m$ epochs are neglected because they are often characterized by strong oscillations due to the optimizer initialization and the different loss functions. In the numerical experiment, we use $N_\text{check}=10$, $p=10$, $N_\text{negl}^m=100(m+1)$.

\begin{figure}[t!]
 \centering
\begin{subfigure}[H]{0.4\linewidth}
  \includegraphics[width=0.89\columnwidth,keepaspectratio,clip]{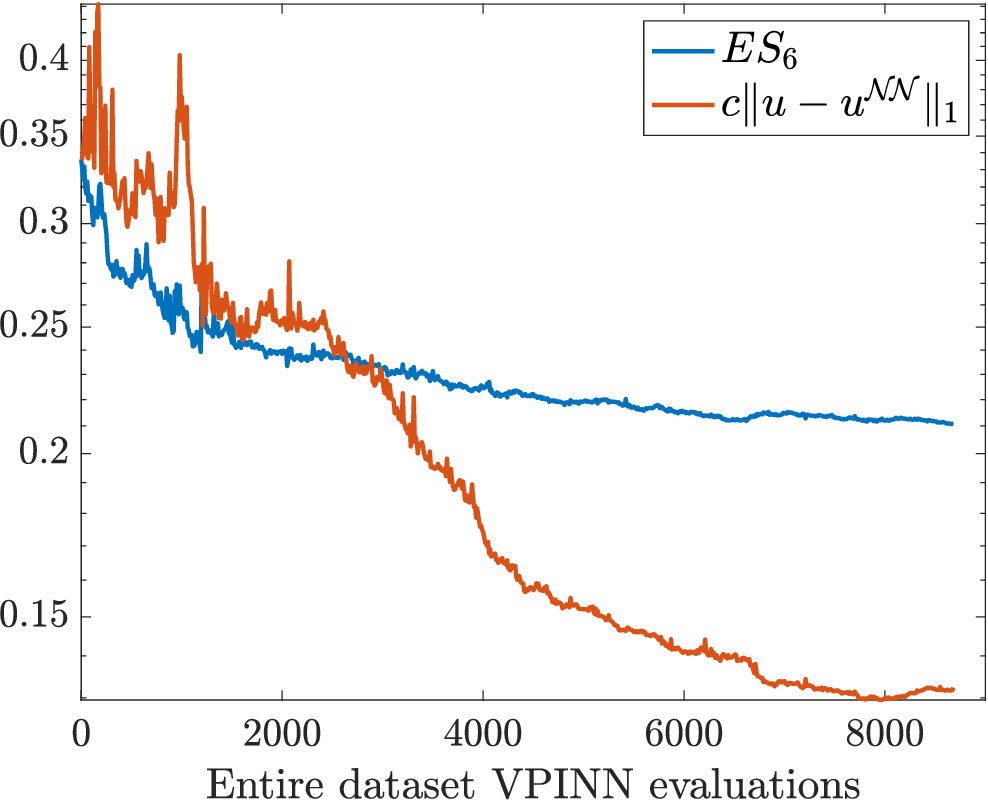}
   \subcaption{\centering}
\end{subfigure}\hspace{0.1\linewidth}
\begin{subfigure}[H]{0.4\linewidth}
  \includegraphics[width=0.89\columnwidth,keepaspectratio,clip]{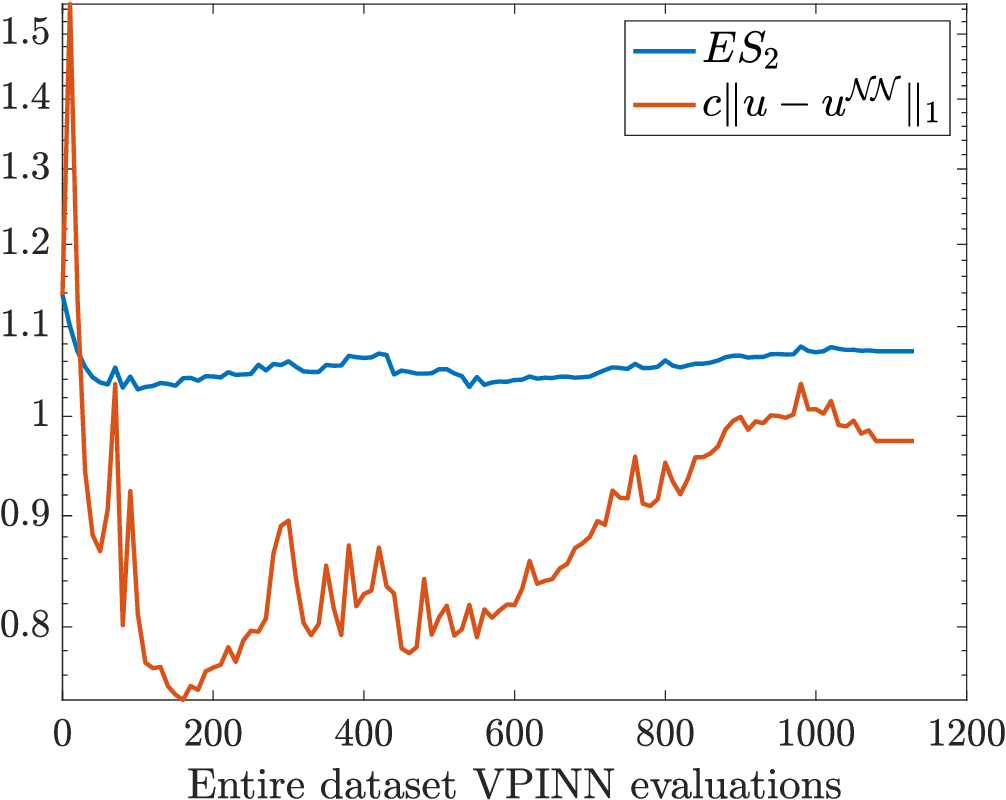} 
   \subcaption{\centering}
\end{subfigure}

\medskip

\begin{subfigure}[H]{0.4\linewidth}
  \includegraphics[width=0.95\columnwidth,keepaspectratio,clip]{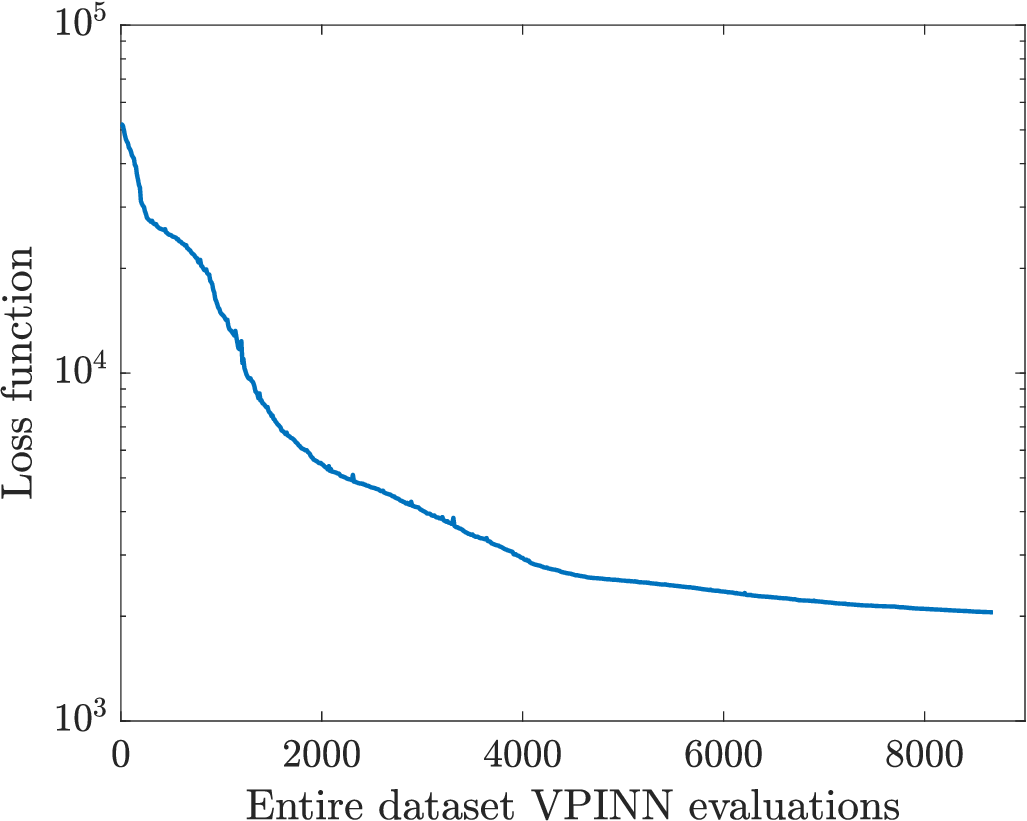}
   \subcaption{\centering}
\end{subfigure}\hspace{0.1\linewidth}
\begin{subfigure}[H]{0.4\linewidth}
  \includegraphics[width=0.89\columnwidth,keepaspectratio,clip]{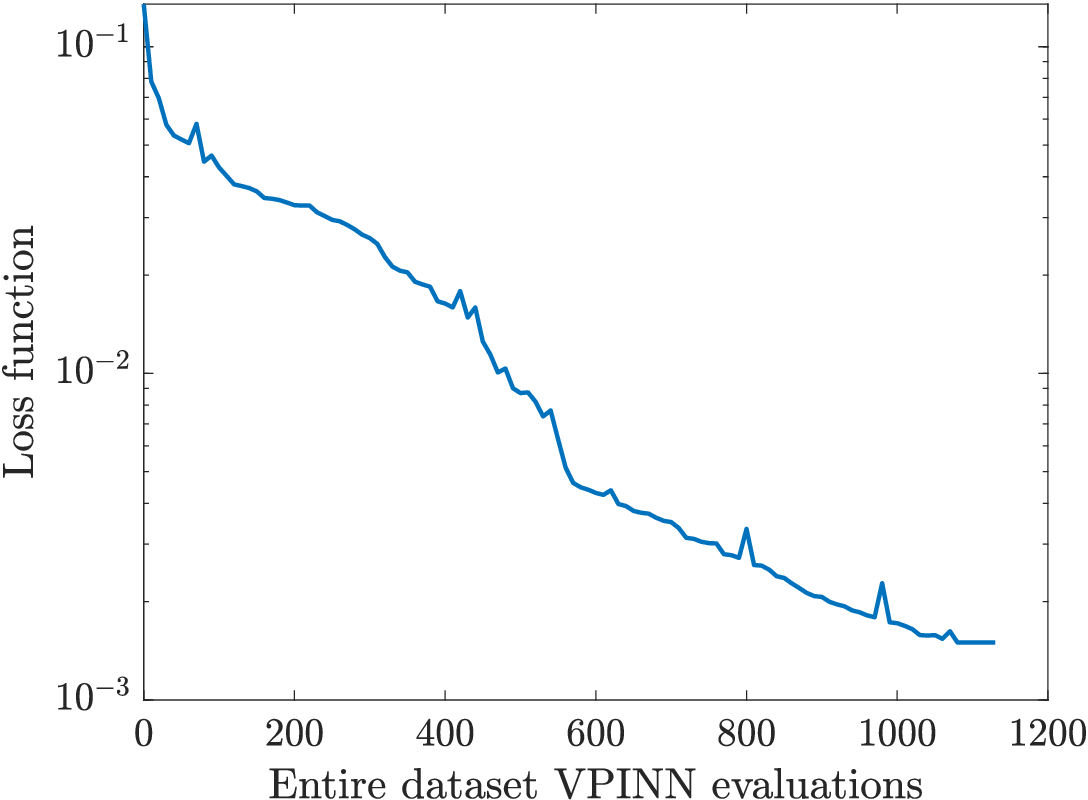}
   \subcaption{\centering}
\end{subfigure}

\medskip
  \caption{{Top row}
: error indicator $ES_m$ and rescaled $H^1$ error $c\Vert u - u^\NN\Vert$. Bottom row: loss function. Left column: curves for the training with patches in ${\cal P}_6$ shown in {Figure} \ref{fig:random_9_energy75_estimators}f. Right column: curves for the training with patches  in ${\cal P}_2$ in {Figure} \ref{fig:struct_9_energy75_removingOld_forceSomeLevelsAndRefine_estimators}a. {(\textbf{a})} 
 $ES_6$ and $c\Vert u - u^\NN\Vert$ for patches in Figure  \ref{fig:random_9_energy75_estimators}f; (\textbf{b}) $ES_2$ and $c\Vert u - u^\NN\Vert$ for patches in {Figure} \ref{fig:struct_9_energy75_removingOld_forceSomeLevelsAndRefine_estimators}a; (\textbf{c}) Loss function for patches in {Figure} \ref{fig:random_9_energy75_estimators}f; (\textbf{d}) Loss function for patches in {Figure} \ref{fig:struct_9_energy75_removingOld_forceSomeLevelsAndRefine_estimators}a.}
  \label{fig:error_loss_estim}
\end{figure}

Two typical scenarios are shown in Figure \ref{fig:error_loss_estim}. In the top row, the behaviors of $ES_m$ and of $c\Vert u - u^\NN\Vert_1$ are shown. Here, $c$ is a scaling parameter used for visualization purposes, chosen such that $ES_m$ and $c\Vert u - u^\NN\Vert_1$ coincide at the beginning of the training. Indeed, $\Vert u - u^\NN\Vert_1$ is about two orders of magnitude smaller than $ES_m$. Nevertheless, it can be noted that these two quantities display very similar behaviors during the training. In the bottom row, instead, we represent the corresponding loss function decay. The left column is associated with the training performed using the patches in ${\cal P}_6$ shown in \mbox{{Figure} 
 \ref{fig:random_9_energy75_estimators}f} and the right column with the one performed using the patches in ${\cal P}_2$ in {Figure} \ref{fig:struct_9_energy75_removingOld_forceSomeLevelsAndRefine_estimators}a. We remark that the loss function, $ES_m$ and $c\Vert u - u^\NN\Vert$ are evaluated in the same epochs and that, in real applications, it is not possible to explicitly compute $c\Vert u - u^\NN\Vert$ since $u$ is not known. Moreover, since we use the L-BFGS optimizer, the neural network is evaluated multiple times on the entire training set in each epoch. Therefore, on the $x$-axis of Figure \ref{fig:error_loss_estim} we show the number of neural network evaluations instead of the number of epochs.

It can be noted that the behavior of the quantities shown in the left column is qualitatively different from the ones in the right column. In fact, when the MF-VPINN is trained with ${\cal P}_6$ of {Figure} \ref{fig:random_9_energy75_estimators}f, the error, the error indicator, and the loss function decrease in similar ways. Therefore, there is no need to interrupt the training early since the accuracy is improving, minimizing the loss function. On the other hand, when the MF-VPINN is trained with the ${\cal P}_2$ of {Figure} \ref{fig:struct_9_energy75_removingOld_forceSomeLevelsAndRefine_estimators}a, the loss decreases even when the error and the error indicator increase or remain constant. In this case, it is convenient to interrupt the training, since minimizing the loss function further would lead to more severe overfitting phenomena and a loss in accuracy and efficiency. At the end of the training, the neural network's trainable parameters corresponding to the best value of $ES_2$ are restored. We highlight that such a phenomenon, observed in \cite{berrone2022solving} too, highlights the fact that the minimization of the loss function generates spurious oscillations that cannot be controlled and ruin the model accuracy. The issue can be partially alleviated with the adopted regularization or completely removed using inf-sup stable models as in \cite{berrone2022variational}.

\begin{itemize}
\item[] \hspace{-0.8cm}\textit{{Strategy \#4: Adaptive strategy without the error indicator}}
\end{itemize}

\begin{figure}[t!]
 \centering
\begin{subfigure}[H]{0.48\linewidth}
  \includegraphics[width=0.99\columnwidth,keepaspectratio,clip]{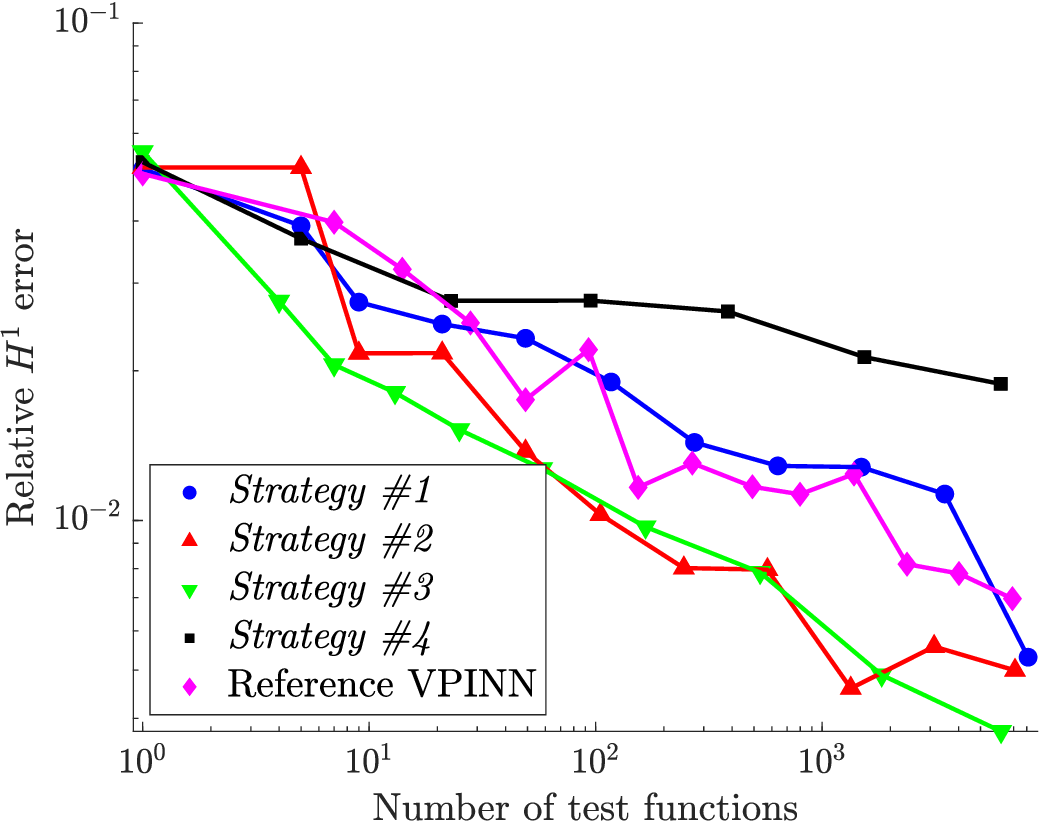}
   \subcaption{\centering}
     \label{fig:error_decay_cm4_all}
\end{subfigure}
\begin{subfigure}[H]{0.48\linewidth}
  \includegraphics[width=0.99\columnwidth,keepaspectratio,clip]{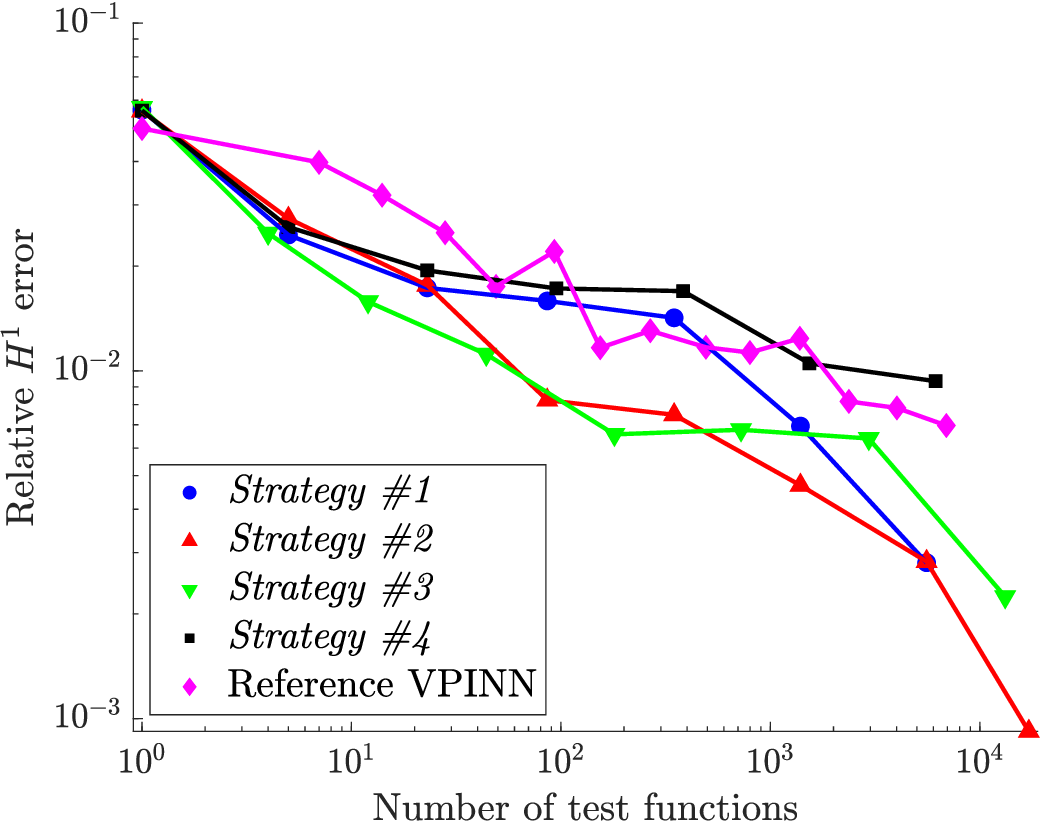}
   \subcaption{\centering}
     \label{fig:error_decay_cm9_all}
\end{subfigure}
  \caption{{Comparison between} 
 the relative $H^1$ errors obtained at the end of each training iteration with different strategies to choose the position of the new patches. {(\textbf{a}) $C_M=4$; (\textbf{b}) $C_M=9$}.}
    \label{fig:error_decay_all}
\end{figure}

Let us now analyze the consequences of choosing the position of the new patches without using the error indicator. To do so, we consider \textit{Strategy \#1} but, instead of considering the new centers inside the patches $P_i$ with the highest values of $\eta_i^\gamma$, we add them inside the patches with the highest values of $r_{h,i}^2(u^\NN)$. Using the equation residuals is a common choice in PINN adaptivity because the residuals describe how accurately the neural network satisfies the PDE at that point. The obtained error decay is shown in {Figure} \ref{fig:error_decay_all}. It can be seen that the accuracy is worse than the ones obtained with the other strategies and that the convergence rate with respect to the number of patches is lower. In such a figure, we also compare the MF-VPINN with a standard VPINN trained with test functions defined on Delaunay meshes. Note that, when \textit{Strategy \#2} or \textit{Strategy \#3} is adopted, the MF-VPINN is more accurate than a simple VPINN, even though its main advantage resides in being a meshfree method.

We highlight that, due to the low regularity of the solution, the expected convergence rate with respect to the number of test functions of an FEM solution computed on uniform refinements is $-$1/3. Note that the convergence rate of the proposed MF-VPINN method is still close to $-$1/3, even though it is a meshfree method (see Table \ref{tab:conv_rates}). For completeness, we also remark that, if an adaptive FEM is used, the rate of convergence depends on the \mbox{FEM order}.

\begin{table}
\begin{center}
\label{tab:conv_rates}
\begin{tabular}{|c|c|c|c|c|c|}
 \hline
$C_M$  & \textit{Strategy \#1} & \textit{Strategy \#2} & \textit{Strategy \#3} & \textit{Strategy \#4} & Reference VPINN\\
 \hline\hline
 4  & -0.213    & -0.295 &  -0.283 & -0.105 & -0.232\\
 9 & -0.294    &-0.376&   -0.287 & -0.182 & -0.232\\
 \hline
\end{tabular}
 \caption{Rates of convergence with respect to the number of test functions. }
\label{tab:conv_rates}
\end{center}
\end{table}

\begin{figure}[t!]
 \centering

\begin{subfigure}[H]{0.325\linewidth}
 \hspace{+0.2cm} \includegraphics[width=0.89\columnwidth,keepaspectratio,clip]{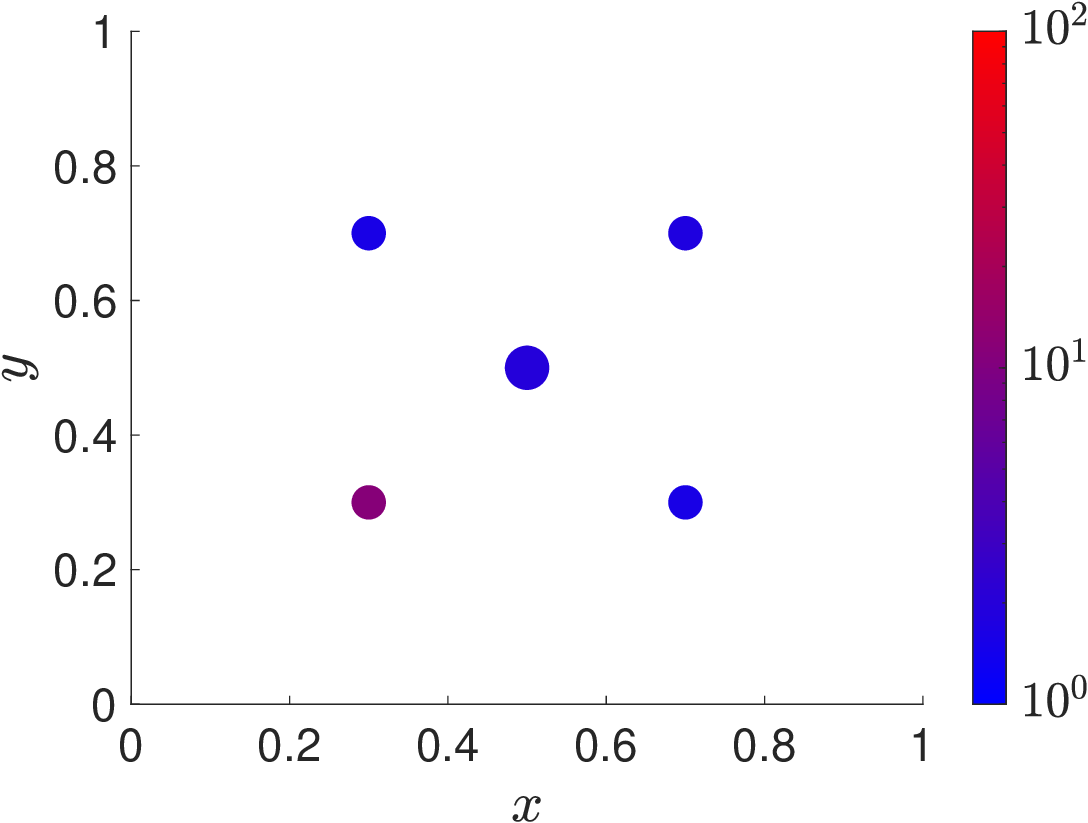}
   \subcaption{\centering}
\end{subfigure}
\begin{subfigure}[H]{0.325\linewidth}
 \hspace{+0.2cm} \includegraphics[width=0.89\columnwidth,keepaspectratio,clip]{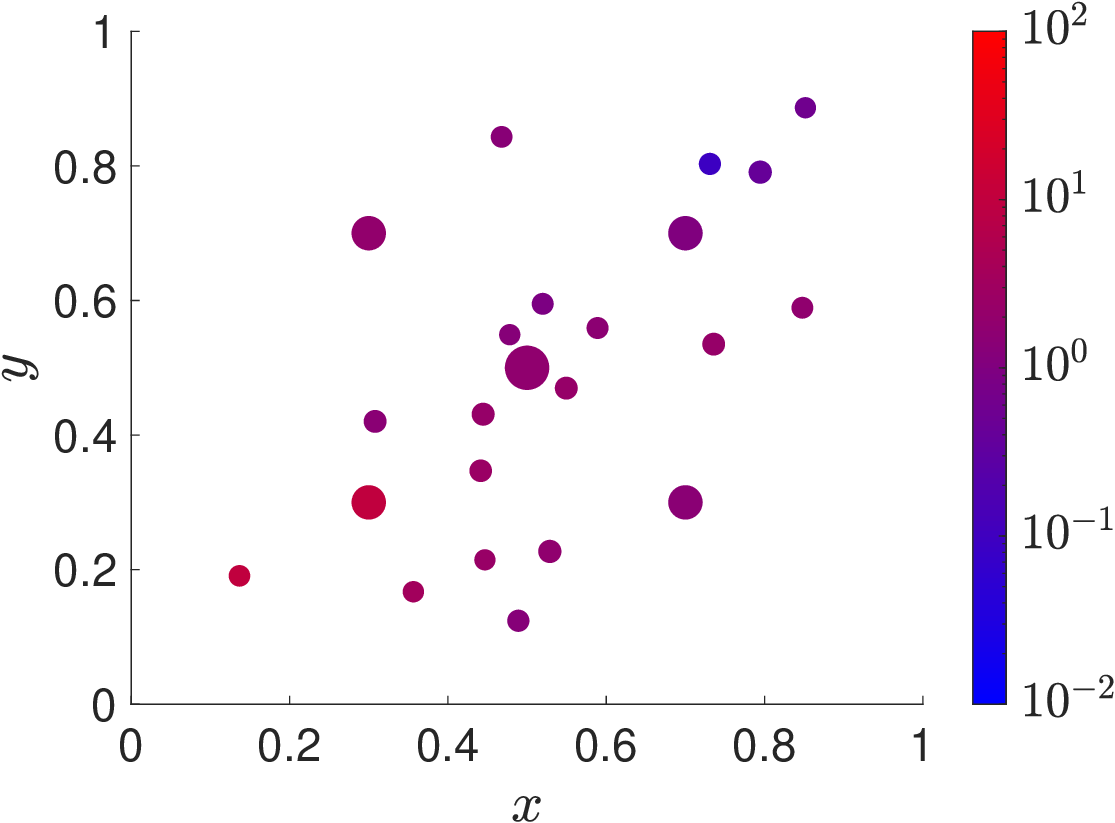}
   \subcaption{\centering}
\end{subfigure}
\begin{subfigure}[H]{0.325\linewidth}
 \hspace{+0.2cm} \includegraphics[width=0.89\columnwidth,keepaspectratio,clip]{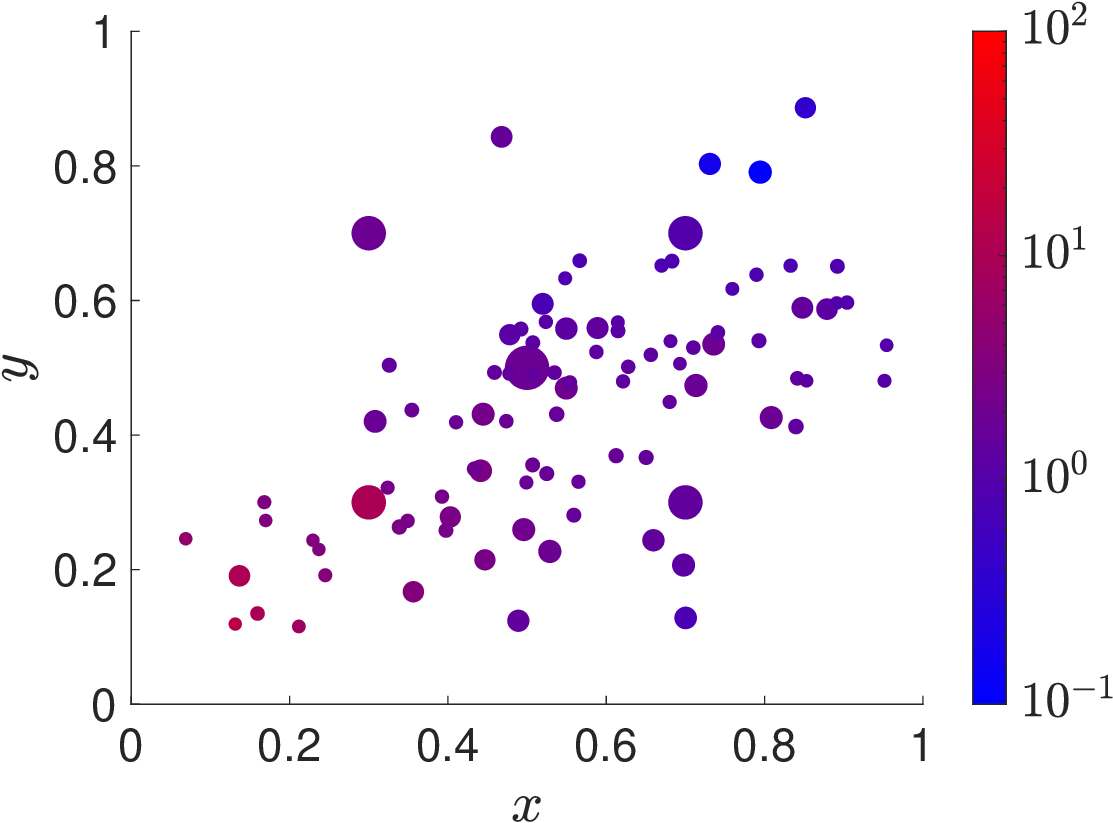} 
   \subcaption{\centering}
\end{subfigure}

\medskip

\begin{subfigure}[H]{0.325\linewidth}
 \hspace{+0.2cm} \includegraphics[width=0.89\columnwidth,keepaspectratio,clip]{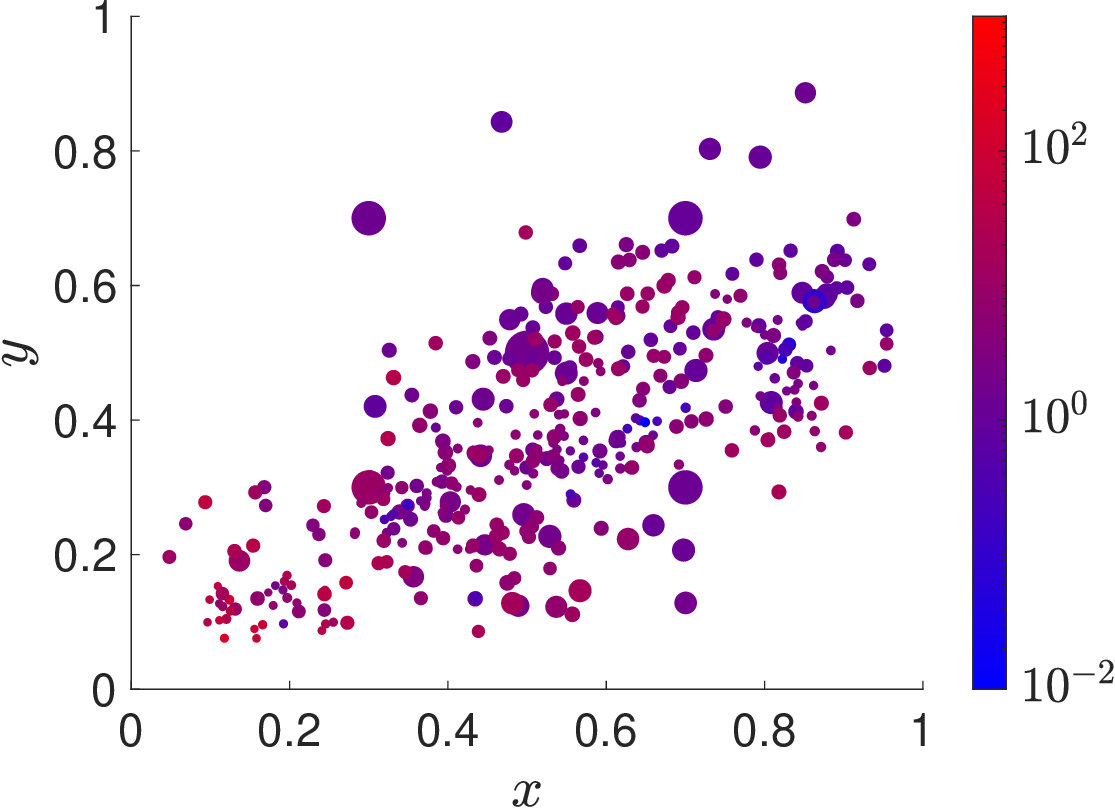}
   \subcaption{\centering}
\end{subfigure}
\begin{subfigure}[H]{0.325\linewidth}
  \hspace{+0.2cm}\includegraphics[width=0.89\columnwidth,keepaspectratio,clip]{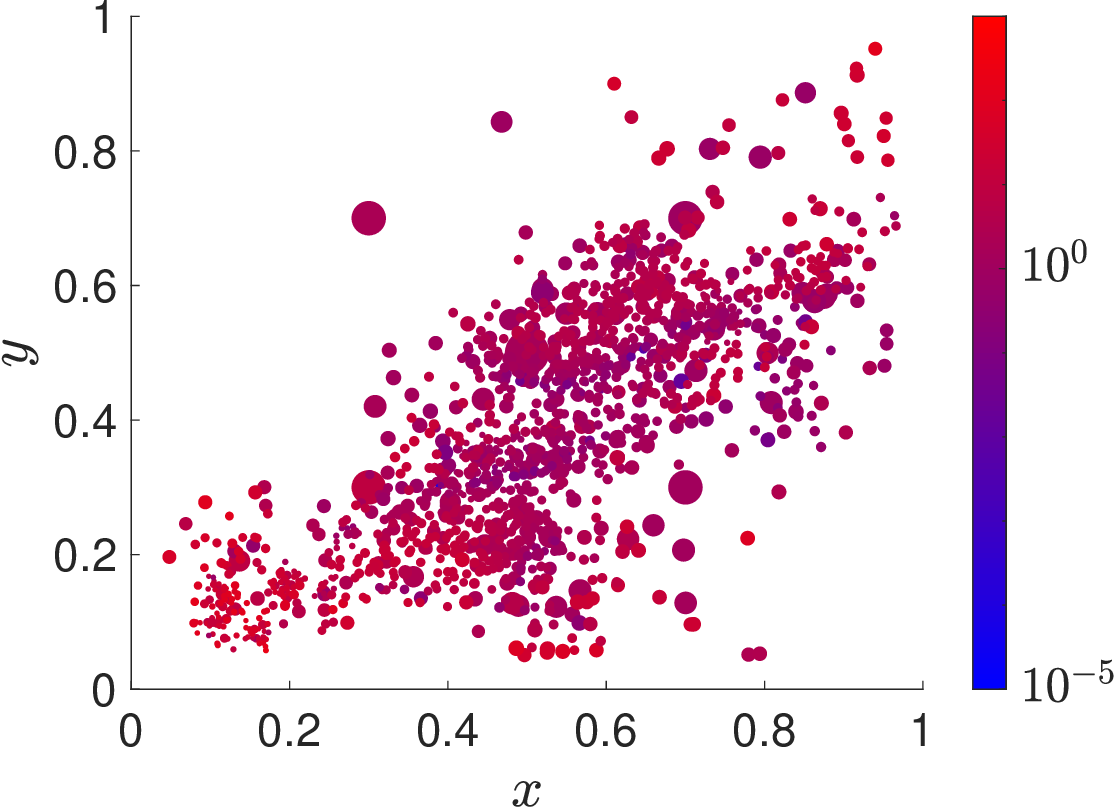}
   \subcaption{\centering}
\end{subfigure}
\begin{subfigure}[H]{0.325\linewidth}
   \hspace{+0.2cm}\includegraphics[width=0.89\columnwidth,keepaspectratio,clip]{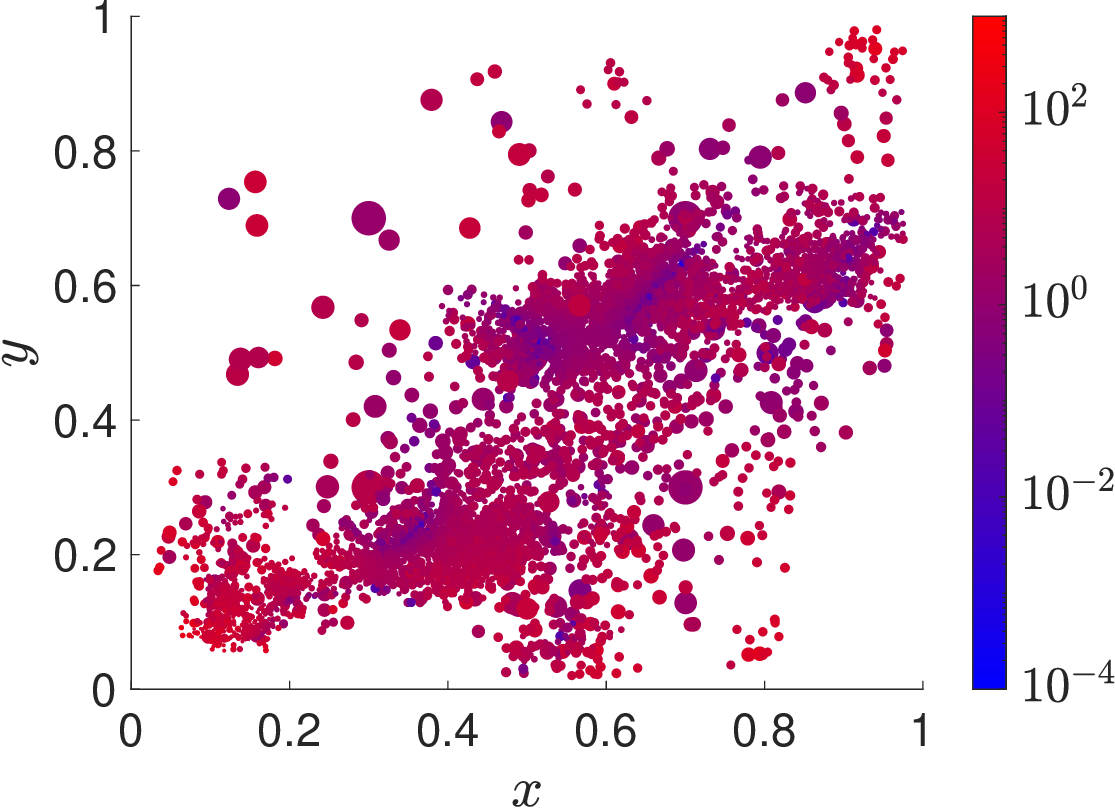}
  \subcaption{\centering}
\end{subfigure}

\medskip
  \caption{\textit{{Strategy \#4}
}: Patches used to train the MF-VPINN with $C_M=4$. Each dot represents a patch $P_i$, its position is the center $\bm{c_{P_i}}$ of the patch, its size is proportional to the patch size $h_i^2$, and its color is associated with the quantity $r_{h,i}^2(u^\NN)$. {(\textbf{a}) Representation of ${\cal P}_1$; (\textbf{b}) Representation of ${\cal P}_2$; (\textbf{c})~Representation of ${\cal P}_3$; (\textbf{d}) Representation of ${\cal P}_4$; (\textbf{e}) Representation of ${\cal P}_5$; (\textbf{f}) Representation of ${\cal P}_6$}. 
}
  \label{fig:strategy4_estimators}
\end{figure}

\begin{figure}[t!]
 \centering
\begin{subfigure}[H]{0.325\linewidth}
  \includegraphics[width=0.99\columnwidth,keepaspectratio,clip]{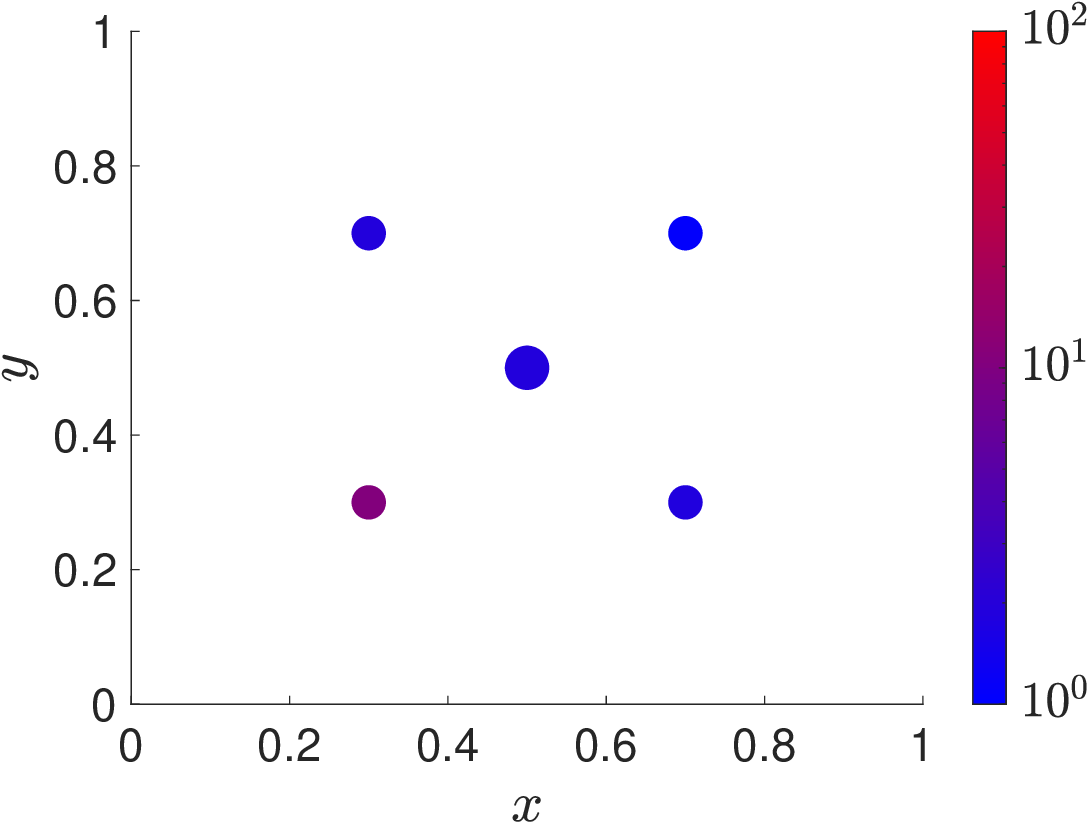}
   \subcaption{\centering}
\end{subfigure}
\begin{subfigure}[H]{0.325\linewidth}
  \includegraphics[width=0.99\columnwidth,keepaspectratio,clip]{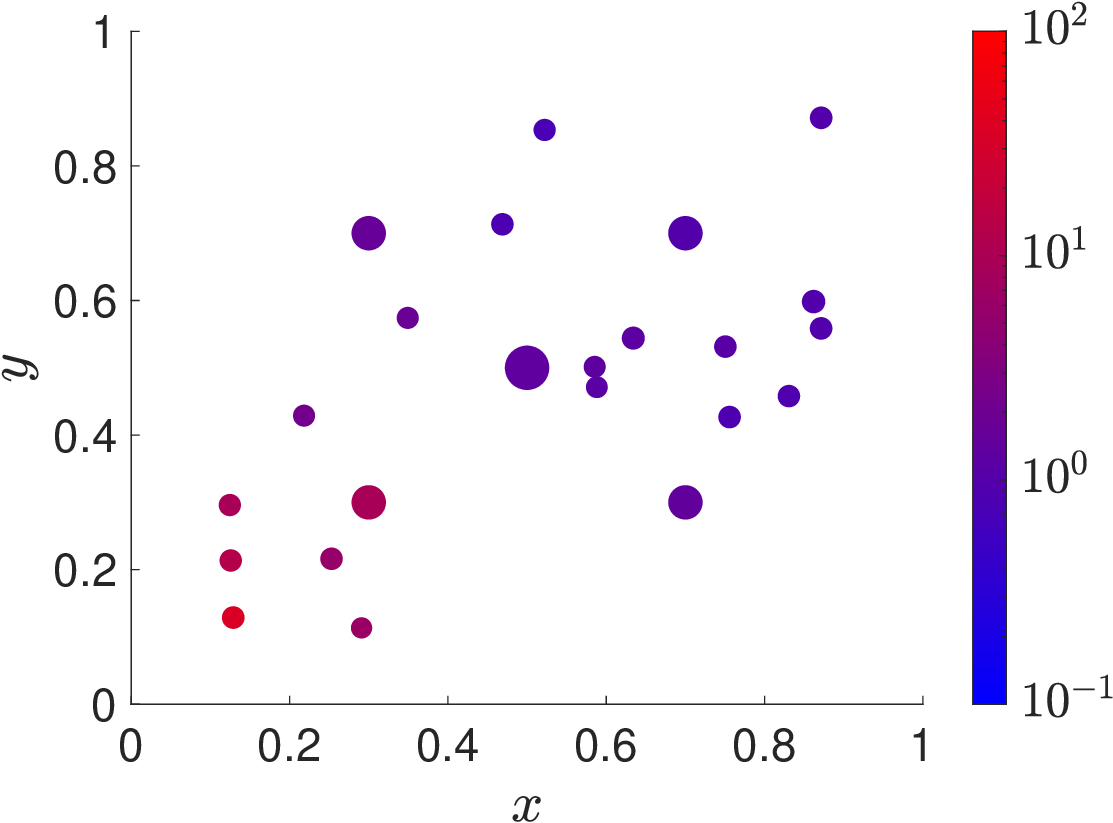}
   \subcaption{\centering}
\end{subfigure}
\begin{subfigure}[H]{0.325\linewidth}
  \includegraphics[width=0.99\columnwidth,keepaspectratio,clip]{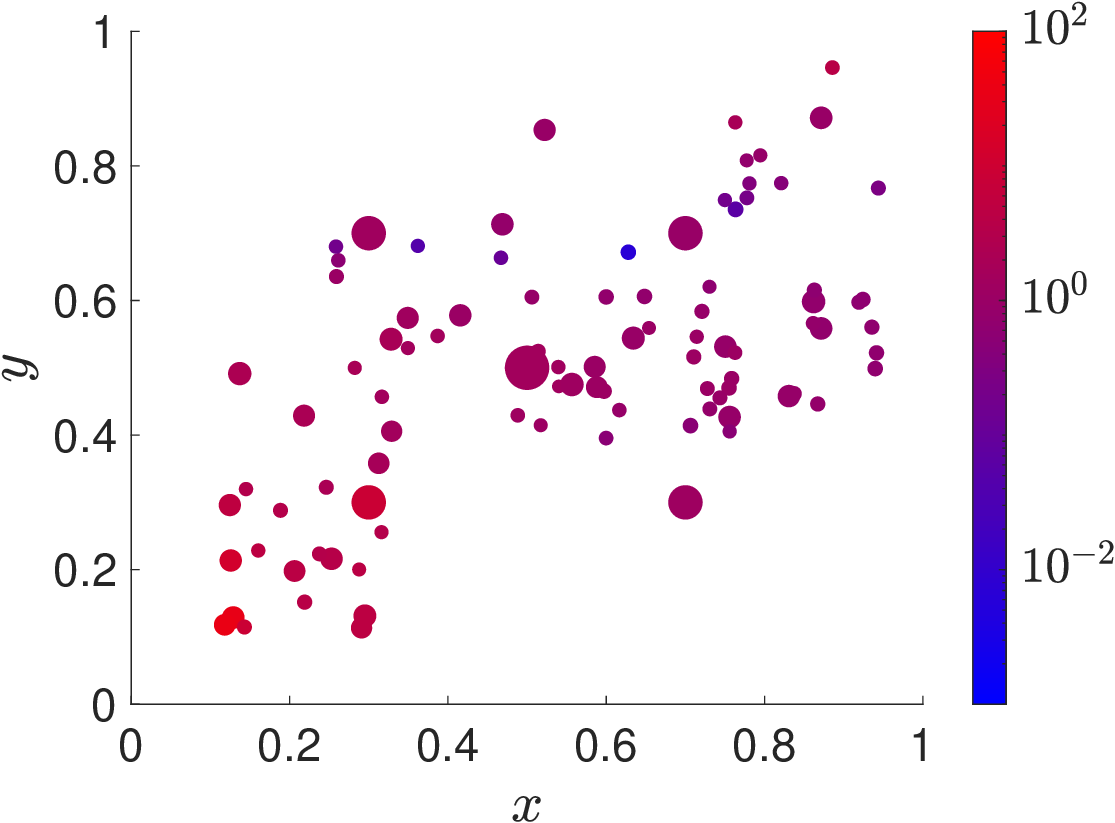} 
   \subcaption{\centering}
\end{subfigure}

\medskip

\begin{subfigure}[H]{0.325\linewidth}
  \includegraphics[width=0.99\columnwidth,keepaspectratio,clip]{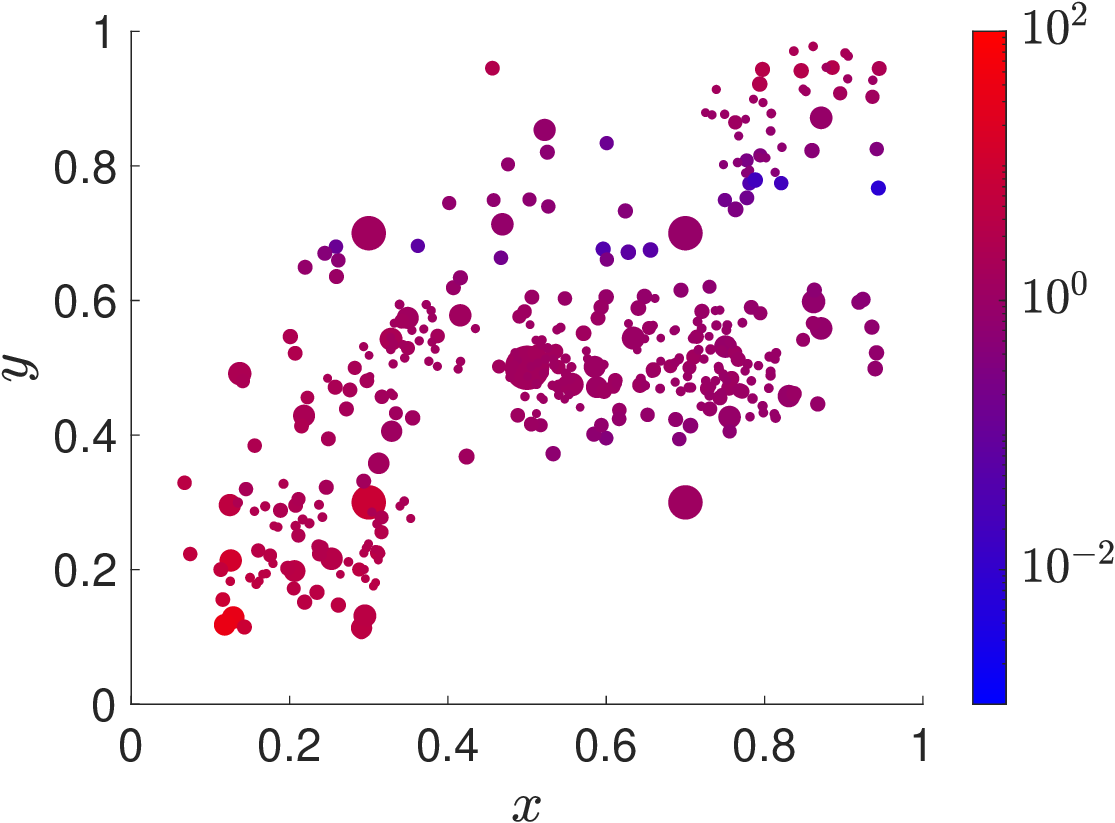}
   \subcaption{\centering}
\end{subfigure}
\begin{subfigure}[H]{0.325\linewidth}
  \includegraphics[width=0.99\columnwidth,keepaspectratio,clip]{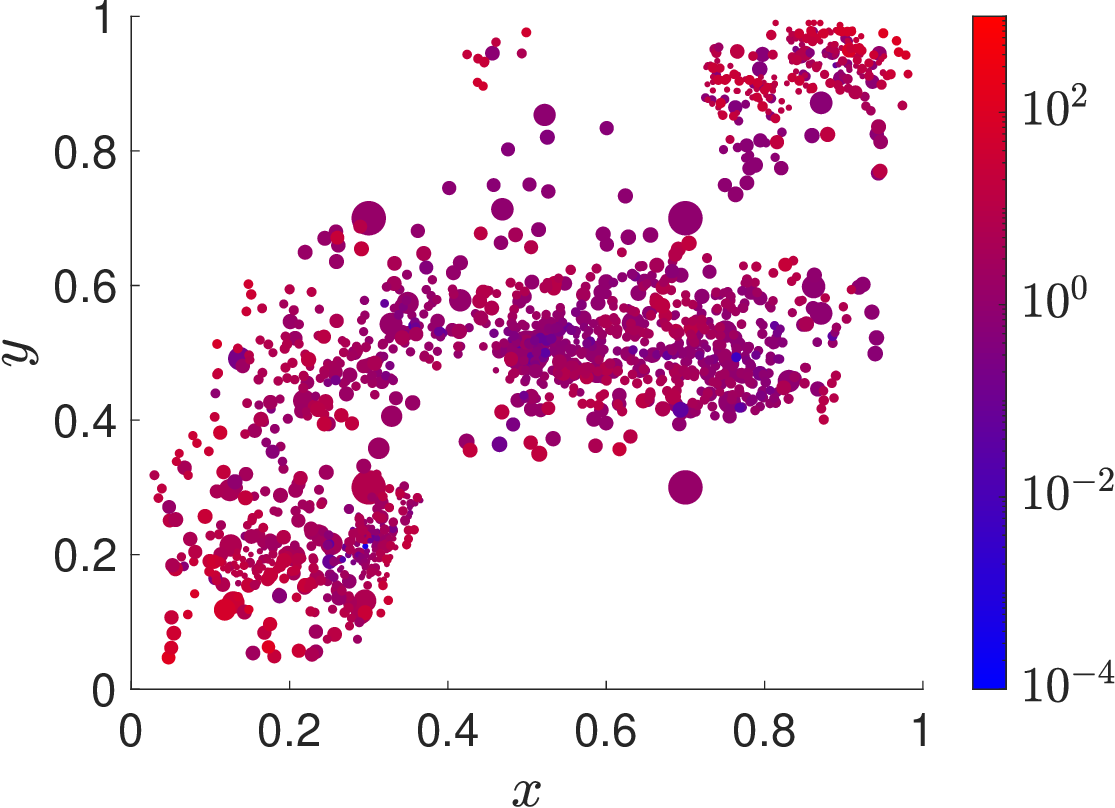}
   \subcaption{\centering}
\end{subfigure}
\begin{subfigure}[H]{0.325\linewidth}
  \includegraphics[width=0.99\columnwidth,keepaspectratio,clip]{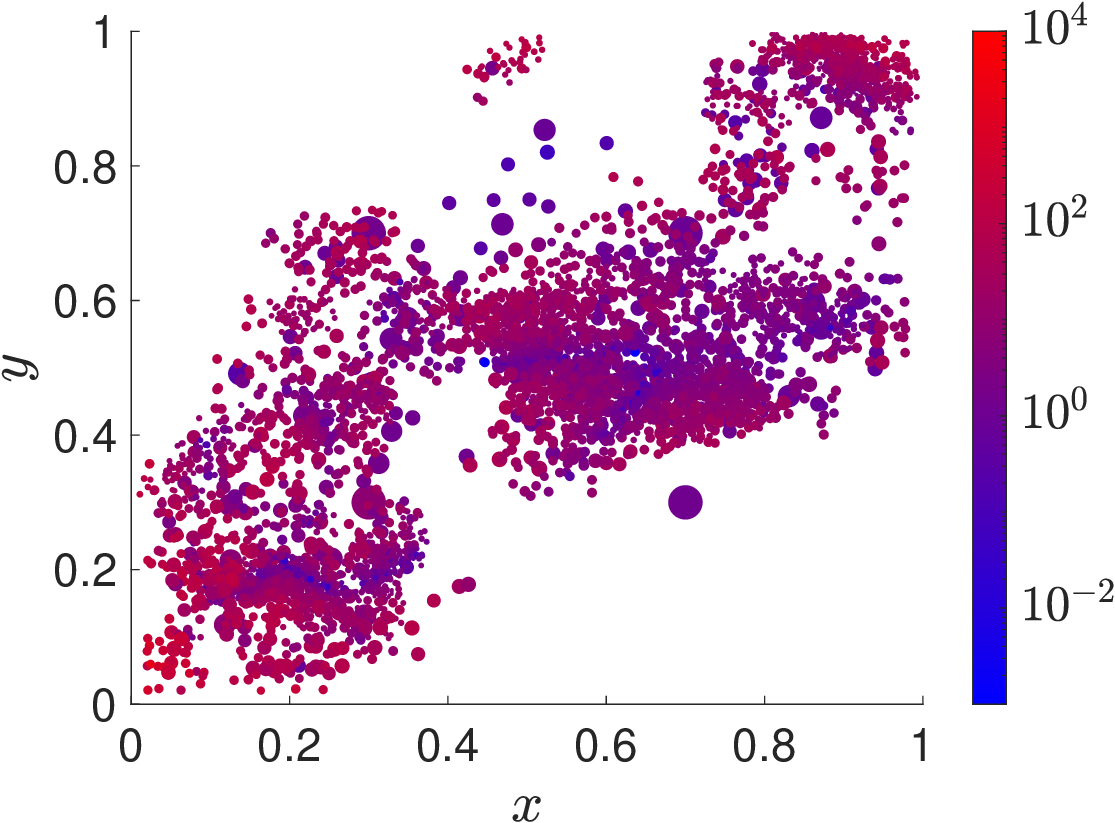}
   \subcaption{\centering}
\end{subfigure}

\medskip
  \caption{\textit{{Strategy \#4}
}: Patches used to train the MF-VPINN with $C_M=9$. Each dot represents a patch $P_i$, its position is the center $\bm{c_{P_i}}$ of the patch, its size is proportional to the patch size $h_i^2$, and its color is associated with the quantity $r_{h,i}^2(u^\NN)$. {(\textbf{a}) Representation of ${\cal P}_1$; (\textbf{b}) Representation of ${\cal P}_2$; (\textbf{c})~Representation of ${\cal P}_3$; (\textbf{d}) Representation of ${\cal P}_4$; (\textbf{e}) Representation of ${\cal P}_5$; (\textbf{f}) Representation of ${\cal P}_6$}. 
}
  \label{fig:strategy9_estimators}
\end{figure}

Coherently with Figure \ref{fig:error_decay_all}, the best strategies are \textit{Strategy \#2} and \textit{Strategy \#3}, whereas the worst one is \textit{Strategy \#4}, which does not exploit the error indicator. The poor performance of \textit{Strategy \#4} can also be explained by analyzing the corresponding patch  distribution. Such distribution is shown in Figure \ref{fig:strategy4_estimators} for $C_M=4$ and in Figure \ref{fig:strategy9_estimators} for $C_M=9$. These plots highlight that the patches do not accumulate near the origin because the residuals of the patches closer to it are not significantly higher than the other ones. For example, note the different colors in \cref{fig:random_4_energy75_estimators,fig:strategy4_estimators}, since in both cases, we randomly choose the position of $C_M=4$ centers inside the selected patches. Such a property is explained by the fact that, in order to minimize the loss function, the optimizer does not focus on specific regions of the domain. Therefore, the orders of magnitude of all the residuals with similar sizes are very close to each other regardless of the position of the corresponding patches. As discussed regarding Figure \ref{fig:error_loss_estim}, we can conclude that the value of the residuals is not a good indicator of the actual error.

{
\subsection{Extension to More a Complex Domain}\label{sec:holes}
In this section, we present some ideas that can be used to apply the method to more complex domains. }

{
Let us consider a domain $\Omega_2$ with some internal holes and boundary $\partial\Omega_2=\Gamma_2$. In particular, $\Omega_2 = (0,1)^2 \backslash \left( \cup_{i=1}^4 H_i\right)$, where $H_i$, $i=1,2,3,4$ are rectangular holes with centers {$\mathbf{c}_{H_i}$} defined as
\[
\mathbf{c}_{H_1} = \left(\frac{9}{26}, \frac{9}{34}\right),  \hspace{1cm} \mathbf{c}_{H_2} = \left(\frac{17}{26}, \frac{9}{34}\right),
\]\[
\mathbf{c}_{H_3} = \left(\frac{9}{26}, \frac{25}{34}\right),  \hspace{1cm} \mathbf{c}_{H_4} = \left(\frac{17}{26}, \frac{25}{34}\right),
\]
and basis and height equal   $\frac{1}{26}$ and $\frac{1}{34}$, respectively.
}

\begin{figure}[t!]
\centering 
\includegraphics[width=0.75\linewidth]{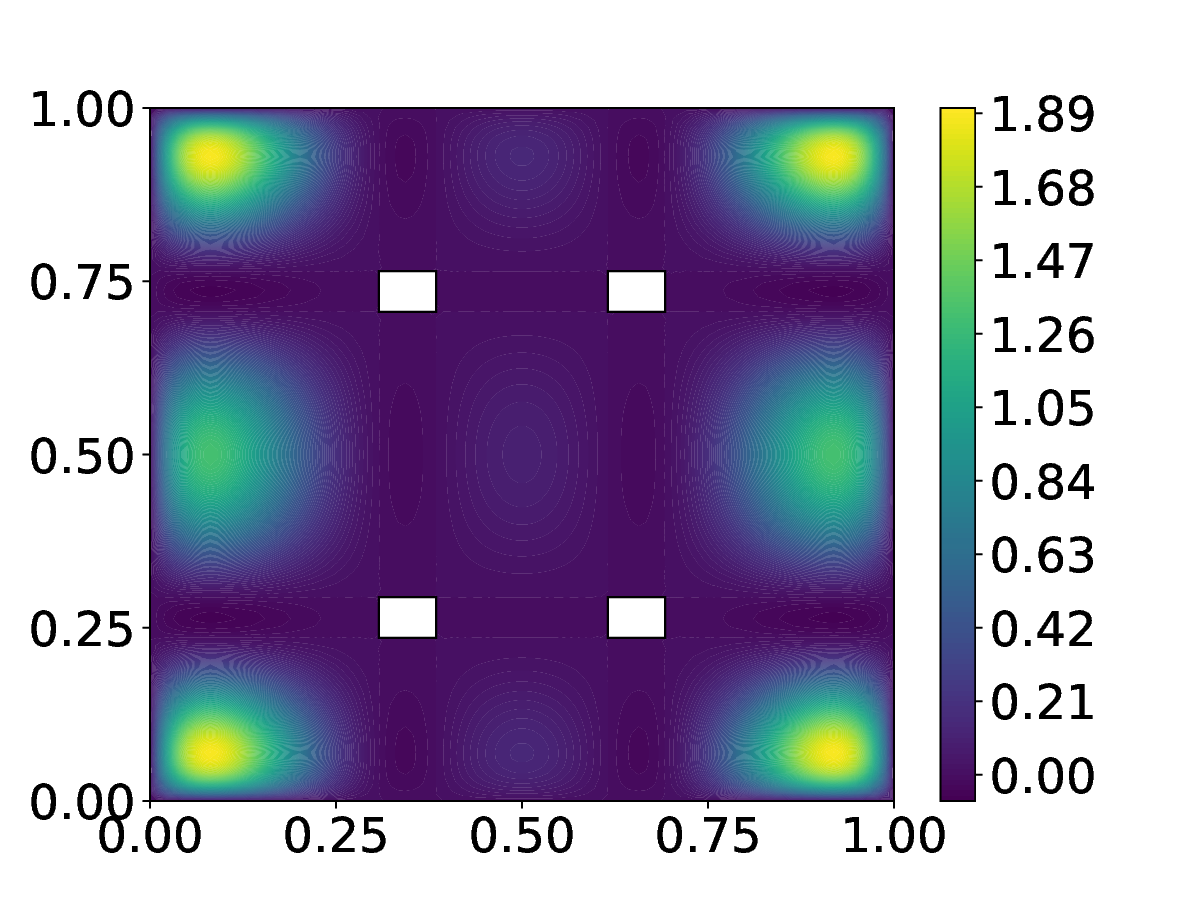} 
  \caption{{Graphical representation of the solution $u$ in \eqref{eq:sol20}.}}
  \label{fig:solution_holes}
\end{figure}

{
In this domain, we consider the Poisson problem:
\begin{equation}\label{eq:model-pb-poisson-2}
\begin{cases}
-\Delta u = f & \text{in \ } \Omega_2\,, \\
\ \ \, u=g & \text{on \ } \Gamma_2 \,, \end{cases}
\end{equation}
with $f$ and $g$ such that the exact solution is
\begin{equation}\label{eq:sol20}
\begin{aligned}
u(x,y) = \frac{1}{C_u}&\left[x(x - 1)\left(x - \dfrac{4}{13}\right)\left(x - \frac{5}{13}\right)\left(x - \frac{8}{13}\right)\left(x - \frac{9}{13}\right)\cdot\right.\\
&\left.y(y - 1)\left(y - \frac{4}{17}\right)\left(y - \frac{5}{17}\right)\left(y - \frac{12}{17}\right)\left(y - \frac{13}{17}\right)\right] ,
\end{aligned}
\end{equation}
normalized through the constant $\frac{1}{C_u}$ to assume value 1 in $\left(\frac{2}{13},\frac{2}{17}\right)$. This function is represented in Figure \ref{fig:solution_holes}.
}

{
We extend the approaches proposed in Section \ref{sec:performance} by adding a cutting procedure after the generation of the new patches. Note that, in particular, all the patches are already completely inside the square $[0,1]^2$ when we apply the cutting procedure, and we can thus focus only on the holes. When a patch intersects more than one hole, we recursively remove it from ${\cal P}_m$, we subdivide the corresponding region in 4 overlapping patches, and we add them to ${\cal P}_m$ until all the generated patches intersect at most one hole. Moreover, we observe that the region $P_i \backslash H_j$ inside the patch $P_i\in{\cal P}_m$ and outside the hole $H_j$, $j=1,2,3,4$, can always be covered by the union of at most four rectangles. When a generated patch intersects a hole, we thus remove the patch and generate the minimum number of patches (at most four) that are as large possible and whose union is the region $P_i \backslash H_j$.}

{To avoid numerical instabilities, when this cutting procedure generates a patch with an aspect ratio larger than 100 or with an area more than 100 times smaller than the original uncutted patch, the new patches are removed from ${\cal P}_m$. This implies that it is not possible to remove the patches associated with the highest error indicators as in \textit{Strategy \#3} because otherwise, the union of all the patches would not cover the entire domain. We thus present numerical results only for \textit{Strategy \#1} and \textit{Strategy \#2}.
}

\begin{figure}[t!]
\centering 
  \includegraphics[width=0.6\linewidth]{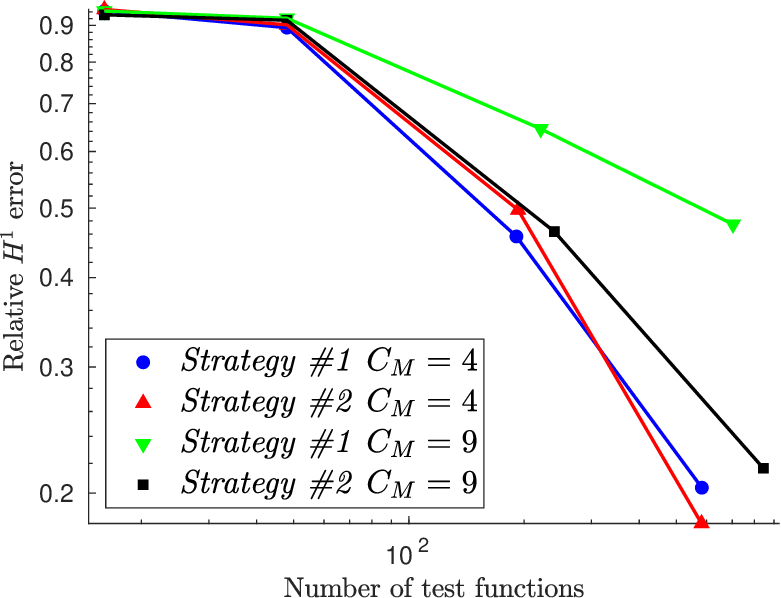} 
  \caption{{Relative $H^1$ errors obtained by solving problem \eqref{eq:model-pb-poisson-2}.}}
  \label{fig:error_holes}
\end{figure}

\begin{figure}[t!]
 \centering
\begin{subfigure}[H]{0.49\linewidth}
  \includegraphics[width=0.99\columnwidth,keepaspectratio,clip]{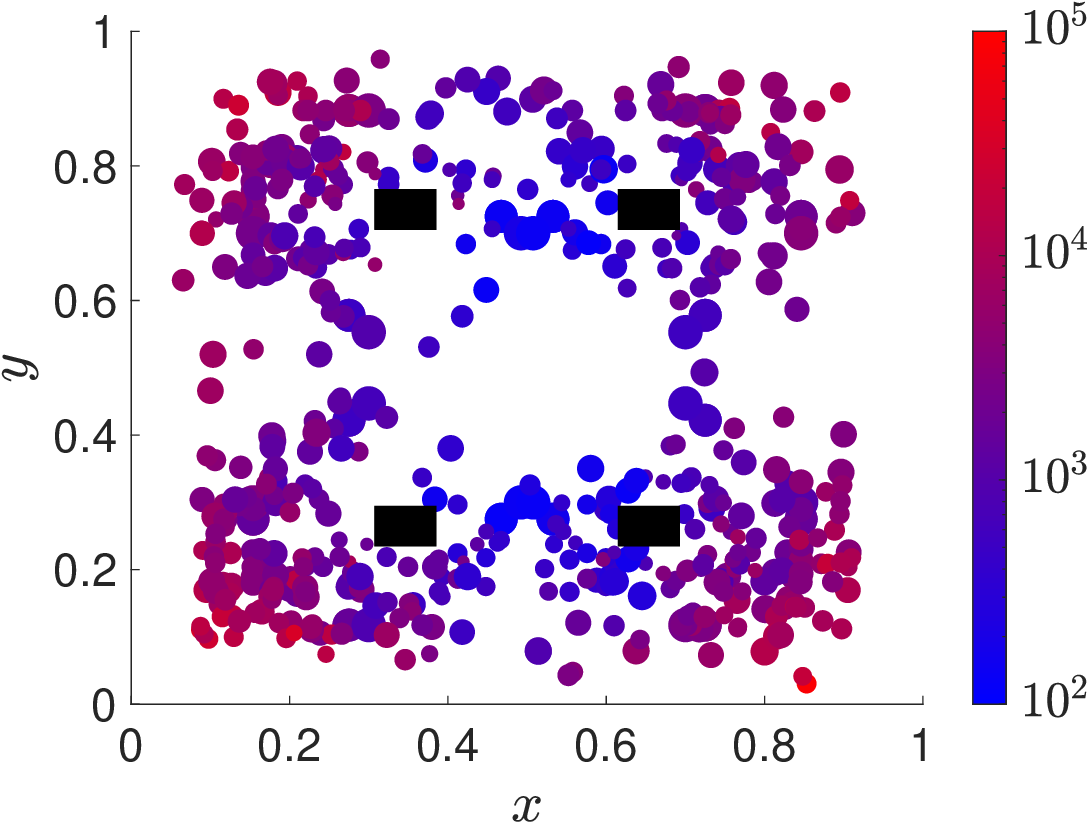}
   \subcaption{\centering}
\end{subfigure}
\begin{subfigure}[H]{0.49\linewidth}
  \includegraphics[width=0.99\columnwidth,keepaspectratio,clip]{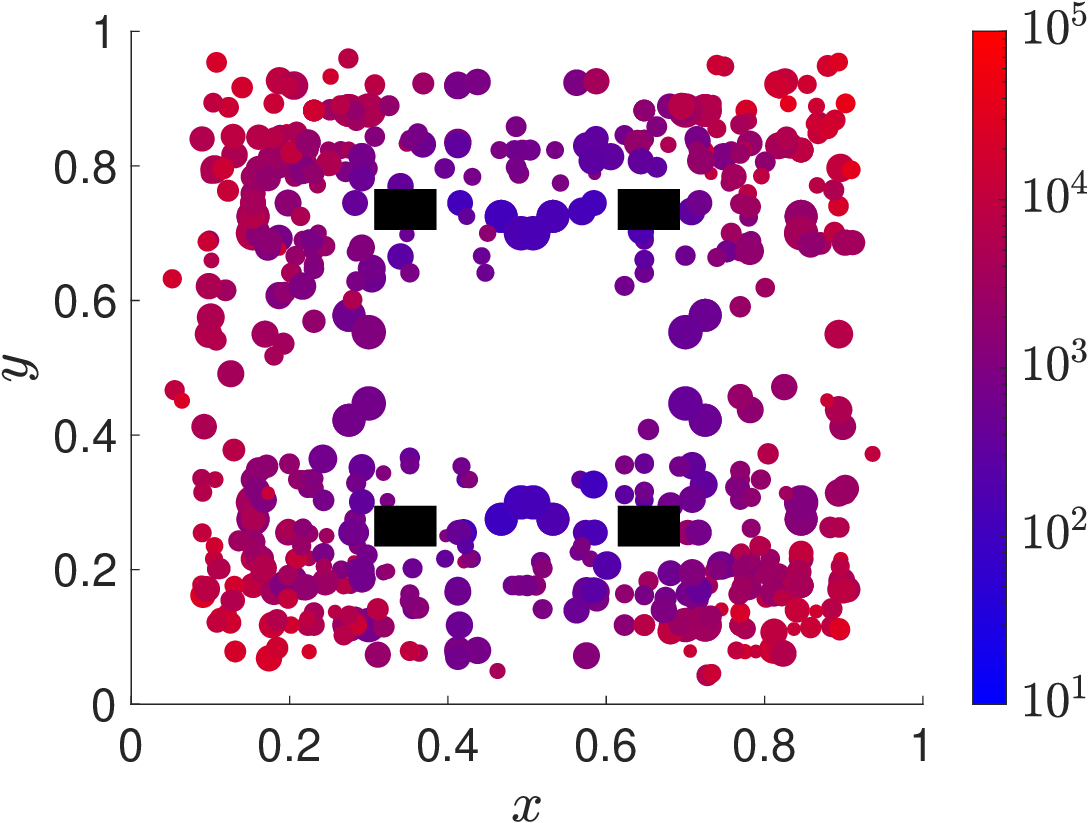}
   \subcaption{\centering}
\end{subfigure}

\medskip
\begin{subfigure}[H]{0.49\linewidth}
  \includegraphics[width=0.99\columnwidth,keepaspectratio,clip]{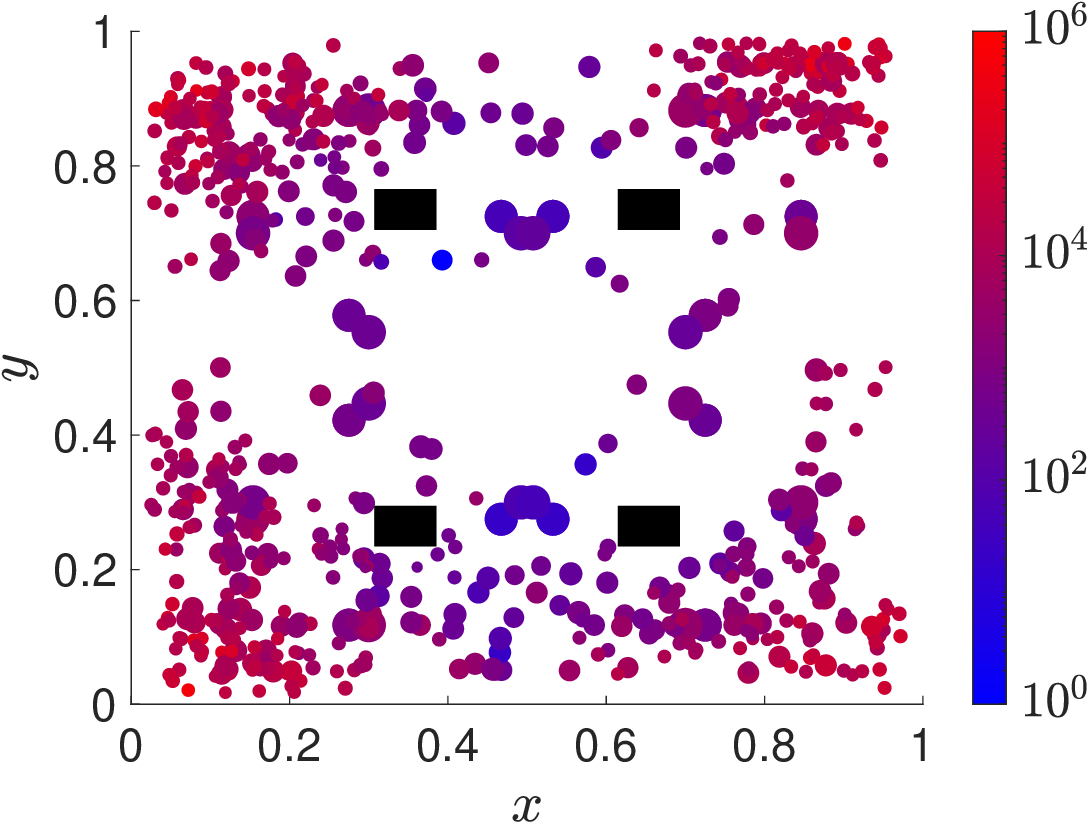}
   \subcaption{\centering}
\end{subfigure}
\begin{subfigure}[H]{0.49\linewidth}
  \includegraphics[width=0.99\columnwidth,keepaspectratio,clip]{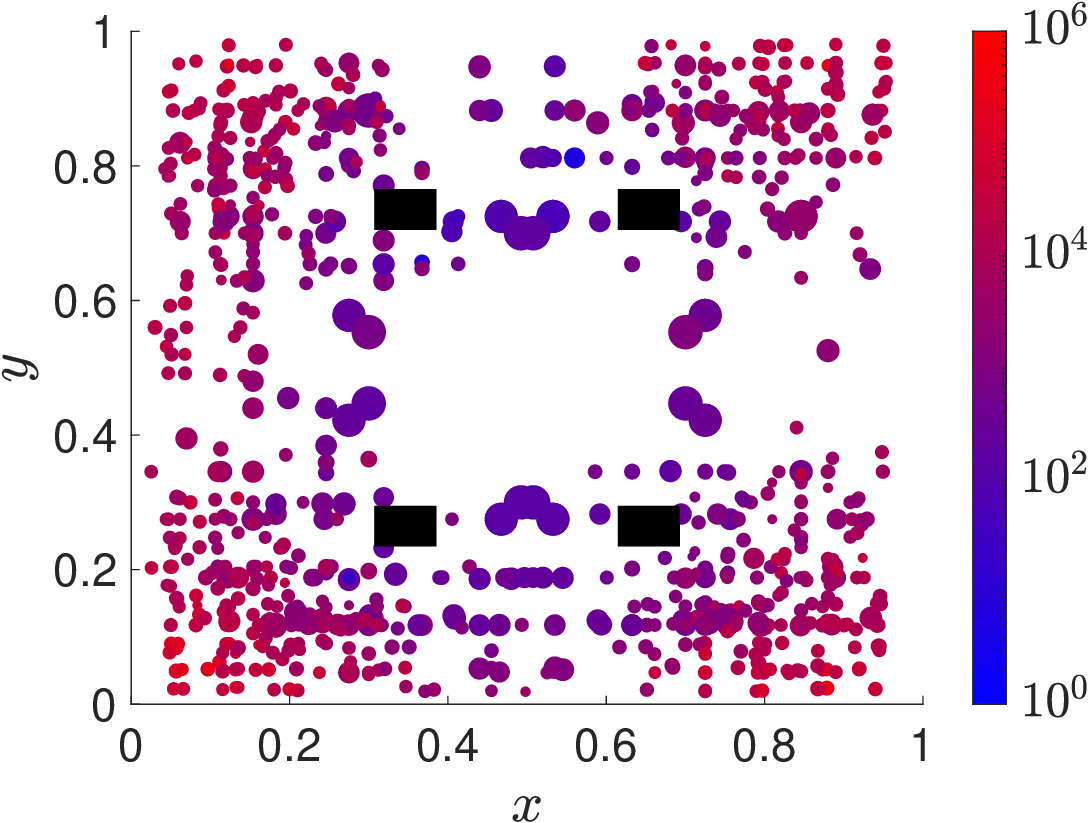}
   \subcaption{\centering}
\end{subfigure}

\medskip
  \caption{{{Problem} 
 \eqref{eq:model-pb-poisson-2}: Representation of the last set of patches obtained with the different strategies. Each dot represents a patch $P_i$, its position is the center {$\bm{c_{P_i}}$} of the patch, its size is proportional to the patch size $h_i^2$, and its color is associated with the quantity $r_{h,i}^2(u^\NN)$. The black rectangles represent the holes $H_i$, $i=1,2,3,4$.} {(\textbf{a}) \textit{Strategy \#1}, $C_M=4$; (\textbf{b}) \textit{Strategy \#2}, $C_M=4$; (\textbf{c}) \textit{Strategy \#1}, $C_M=9$; (\textbf{d})~\textit{Strategy \#2}, $C_M=9$.} 
}
  \label{fig:holes_estimators}
\end{figure}

{The obtained error decays are shown in Figure \ref{fig:error_holes} for \textit{Strategy \#1} and \textit{Strategy \#2} with $C_M=4$ and $C_M=9$.  The first and second errors are computed with the patches generated by cutting the patches in ${\cal P}_0$ and ${\cal P}_1$, respectively, whereas the third and fourth errors are obtained by refining the previous patches with the error indicator as previously described. Note that the first and second errors are very close for all the curves since the strategy and the value of $C_M$ does not influence the training and that both strategies converge better with $C_M=4$. The final patch  distributions are displayed in Figure \ref{fig:holes_estimators}. Here, we can see that the inner part of the domain is covered by a few large patches, whereas the distribution is denser closer to the external boundary of $\Omega_2$, where the solution is more oscillating.}

\section{Conclusions and Discussion}\label{sec:conclusion}
In this work, we presented a Meshfree Variational-Physics-Informed Neural Network (MF-VPINN). It is a PINN trained using the PDE variational formulation that does not require the generation of a global triangulation of the entire domain. In order to generate the test functions involved in the loss computation, we use an a posteriori error estimator based on the one discussed in \cite{berrone2022solving}. Using such an error estimator, it is possible to add test functions only in regions in which the error is higher, thus increasing the efficiency of \mbox{the method}. 

{We highlight that the main advantages of the method are that it is meshfree, as it requires only a covering of the domain with patches that can be of different shapes and that it automatically improves the solution with the application of local patches without requiring a global mesh manipulation. It can be therefore used in domains where it is expensive or impossible to generate a mesh. On the other hand, if a mesh suitable to describe the solution can be generated, a standard VPINN is preferable since the implementation is simpler and the convergence rate with respect to the number of test functions is higher.}

We discuss several strategies to generate the set of test functions. We observe that adding a few test functions inside the patches associated with higher errors while ensuring a smooth transition between regions with large patches and regions with small patches is the best way to obtain accurate solutions. We also show that, if the a posteriori error indicator is not used, the model's accuracy decreases and the training is slower.

In this paper, we only focus on second-order elliptic problem even though VPINNs can be used to solve more complex problems. In a forthcoming paper, we will adapt the a posteriori error estimator and analyze the MF-VPINN performance on other PDEs. Moreover, we are interested in the analysis of the approach in more complex domains (in which the patches have to be suitably deformed) and in high-dimensional problems, where using a standard VPINN is not practical.

\section*{Acknowledgements}
The author S.B. kindly acknowledges partial financial support provided by PRIN project "Advanced polyhedral discretisations of heterogeneous PDEs for multiphysics problems" (No. 20204LN5N5\_003) and by PNRR M4C2 project of CN00000013 National center for HPC, Big Data and Quantum Computing (HPC) (CUP: E13C22000990001). The author M.P. kindly acknowledges the financial support provided by the Politecnico di Torino where the research was carried out.



\begin{thebibliography}{999}

\providecommand{\natexlab}[1]{#1}

\bibitem{lagaris1997artificial}
Lagaris, I.; Likas, A.; Fotiadis, D.
\newblock Artificial neural network methods in quantum mechanics.
\newblock {\em Comput. Phys. Commun.} {\bf 1997}, {\em 104},~1--14.

\bibitem{lagaris1998artificial}
Lagaris, I.; Likas, A.; Fotiadis, D.
\newblock Artificial neural networks for solving ordinary and partial
  differential equations.
\newblock {\em IEEE Trans. Neural Netw.} {\bf 1998}, {\em
  9},~987--1000.

\bibitem{lagaris2000neural}
Lagaris, I.; Likas, A.; Papageorgiou, D.
\newblock Neural-network methods for boundary value problems with irregular
  boundaries.
\newblock {\em IEEE Trans. Neural Netw.} {\bf 2000}, {\em
  11},~1041--1049.

\bibitem{tensorflow2015-whitepaper}
Abadi, M.; {Agarwal, A.; Barham, P.; Brevdo, E; Chen, Z.; Citro, C.; Corrado, G.S.; Davis, A.; Dean, J.; Devin, M}
.; et al.
\newblock TensorFlow: Large-Scale Machine Learning on Heterogeneous Systems.
  2015.
\newblock {Available online}
: \url{https://www.tensorflow.org} ({accessed on 15 September 2024}
).

\bibitem{NEURIPS2019_9015}
Paszke, A.; {Gross, S.; Massa, F.; Lerer, A.; Bradbury, J.; Chanan, G.; Killeen, T.; Lin, Z.; Gimelshein, N.; Antiga, L}.; et~al.
\newblock PyTorch: An Imperative Style, High-Performance Deep Learning Library.
  In {\em Advances in Neural Information Processing Systems 32}; Curran
  Associates, Inc.: {Red Hook, NY, USA}, 2019; pp. 8024--8035.

\bibitem{jax2018github}
Bradbury, J.; Frostig, R.; Hawkins, P.; Johnson, M.J.; Leary, C.; Maclaurin,
  D.; Necula, G.; Paszke, A.; Vander{P}las, J.; \mbox{Wanderman-{M}ilne, S.;  et~al}. {JAX}: Composable Transformations of {P}ython+{N}um{P}y Programs. 2018. Available online: \url{http://github.com/google/jax} ({accessed on 15 September 2024}).

\bibitem{raissi2017physics}
Raissi, M.; Perdikaris, P.; Karniadakis, G.
\newblock Physics informed deep learning (part i): Data-driven solutions of
  nonlinear partial differential equations.
\newblock {\em arXiv} {\bf 2017}, arXiv:1711.10561 .

\bibitem{raissi2017bphysics}
Raissi, M.; Perdikaris, P.; Karniadakis, G.
\newblock Physics informed deep learning (part ii): Data-driven solutions of
  nonlinear partial differential equations.
\newblock {\em arXiv} {\bf 2017}, arXiv:1711.10566 .

\bibitem{raissi2019physics}
Raissi, M.; Perdikaris, P.; Karniadakis, G.
\newblock Physics-informed neural networks: A deep learning framework for
  solving forward and inverse problems involving nonlinear partial differential
  equations.
\newblock {\em J. Comput. Phys.} {\bf 2019}, {\em
  378},~686--707.
\newblock {{https://doi.org/10.1016/j.jcp.2018.10.045}}.

\bibitem{pu2021solving}
Pu, J.; Li, J.; Chen, Y.
\newblock Solving localized wave solutions of the derivative nonlinear
  Schr{\"o}dinger equation using an improved PINN method.
\newblock {\em Nonlinear Dyn.} {\bf 2021}, {\em 105},~1723--1739.

\bibitem{yuan2022pinn}
Yuan, L.; Ni, Y.; Deng, X.; Hao, S.
\newblock A-PINN: Auxiliary physics informed neural networks for forward and
  inverse problems of nonlinear integro-differential equations.
\newblock {\em J. Comput. Phys.} {\bf 2022}, {\em 462},~111260.

\bibitem{guo2023high}
Guo, Q.; Zhao, Y.; Lu, C.; Luo, J.
\newblock High-dimensional inverse modeling of hydraulic tomography by physics
  informed neural network (HT-PINN).
\newblock {\em J. Hydrol.} {\bf 2023}, {\em 616},~128828.

\bibitem{demo2021extended}
Demo, N.; Strazzullo, M.; Rozza, G.
\newblock An extended physics informed neural network for preliminary analysis
  of parametric optimal control problems.
\newblock {\em Comput. Math. Appl.} {\bf 2023}, {\em
  143},~383--396.
\newblock {{https://doi.org/10.1016/j.camwa.2023.05.004}}.

\bibitem{gao2021phygeonet}
Gao, H.; Sun, L.; Wang, J.
\newblock PhyGeoNet: Physics-informed geometry-adaptive convolutional neural
  networks for solving parameterized steady-state PDEs on irregular domain.
\newblock {\em J. Comput. Phys.} {\bf 2021}, {\em 428},~110079.

\bibitem{chen2020electromagnetic}
Yuyao, C.; Lu, L.; Karniadakis, G.; Dal~Negro, L.
\newblock Physics-informed neural networks for inverse problems in nano-optics
  and metamaterials.
\newblock {\em Opt. Express} {\bf 2020}, {\em 28},~11618--11633.
\newblock {{https://doi.org/10.1364/OE.384875}}.

\bibitem{tartakovsky2018learning}
Tartakovsky, A.; Marrero, C.; Perdikaris, P.; Tartakovsky, G.; Barajas-Solano,
  D.
\newblock Learning parameters and constitutive relationships with physics
  informed deep neural networks.
\newblock {\em arXiv} {\bf 2018}, arXiv:1808.03398.

\bibitem{zhao2021physics}
Chen, Z.; Liu, Y.; Sun, H.
\newblock Physics-informed learning of governing equations from scarce data.
\newblock {\em Nat. Commun.} {\bf 2021}, {\em 12}, {6136}.

\bibitem{yu2018deep}
Weinan, E.; Yu, B.
\newblock {The Deep Ritz method: A deep learning-based numerical algorithm for
  solving variational problems}.
\newblock {\em Commun. Math. Stat.} {\bf 2018}, {\em
  6},~1--12.

\bibitem{muller2022error}
M{\"u}ller, J.; Zeinhofer, M.
\newblock Error estimates for the deep Ritz method with boundary penalty.
\newblock In Proceedings of the Mathematical and Scientific Machine Learning.
  PMLR, {Beijing, China, 15--17 August} 2022; pp. 215--230.

\bibitem{lu2021priori}
Lu, Y.; Lu, J.; Wang, M.
\newblock A priori generalization analysis of the deep Ritz method for solving
  high dimensional elliptic partial differential equations.
\newblock In Proceedings of the Conference on Learning Theory. PMLR, {Boulder, CO, USA, 15--19 August 2021}; \mbox{pp. 3196--3241}.

\bibitem{sirignano2018dgm}
Sirignano, J.; Spiliopoulos, K.
\newblock {DGM: A deep learning algorithm for solving partial differential
  equations}.
\newblock {\em J. Comput. Phys.} {\bf 2018}, {\em
  375},~1339--1364.
\newblock {{https://doi.org/10.1016/j.jcp.2018.08.029}}.

\bibitem{al2022extensions}
Al-Aradi, A.; Correia, A.; Jardim, G.; de~Freitas~Naiff, D.; Saporito, Y.
\newblock Extensions of the deep Galerkin method.
\newblock {\em Appl. Math. Comput.} {\bf 2022}, {\em
  430},~127287.

\bibitem{li2021deep}
Li, J.; Zhang, W.; Yue, J.
\newblock A deep learning Galerkin method for the second-order linear elliptic
  equations.
\newblock {\em Int. J. Numer. Anal. Model.} {\bf
  2021}, {\em 18}, {427--441}.

\bibitem{smith1997domain}
Smith, B.F.
\newblock Domain decomposition methods for partial differential equations. In
  {\em Parallel Numerical Algorithms}; Springer:  {Berlin/Heidelberg, Germany,} 
 1997; pp. 225--243.

\bibitem{toselli2006domain}
Toselli, A.; Widlund, O.
\newblock {\em Domain Decomposition Methods-Algorithms and Theory}; 
  Springer Science \& Business Media: {Berlin/Heidelberg, Germany,} 2006; Volume 34.

\bibitem{ameya2020conservative}
Jagtap, A.; Kharazmi, E.; Karniadakis, G.
\newblock Conservative physics-informed neural networks on discrete domains for
  conservation laws: Applications to forward and inverse problems.
\newblock {\em Comput. Methods Appl. Mech. Eng.} {\bf
  2020}, {\em 365},~113028.
\newblock {{https://doi.org/10.1016/j.cma.2020.113028}}.

\bibitem{shukla2021parallel}
Shukla, K.; Jagtap, A.D.; Karniadakis, G.E.
\newblock Parallel physics-informed neural networks via domain decomposition.
\newblock {\em J. Comput. Phys.} {\bf 2021}, {\em 447},~110683.

\bibitem{jagtap2020extended}
Jagtap, A.; Karniadakis, G.
\newblock {Extended physics-informed neural networks (XPINNs): A generalized
  space-time domain decomposition based deep learning framework for nonlinear
  partial differential equations}.
\newblock {\em Commun. Comput. Phys.} {\bf 2020}, {\em
  28},~2002--2041.

\bibitem{moseley2023finite}
Moseley, B.; Markham, A.; Nissen-Meyer, T.
\newblock Finite Basis Physics-Informed Neural Networks (FBPINNs): A scalable
  domain decomposition approach for solving differential equations.
\newblock {\em Adv. Comput. Math.} {\bf 2023}, {\em 49},~62.

\bibitem{viana2021estimating}
Viana, F.; Nascimento, R.; Dourado, A.; Yucesan, Y.
\newblock Estimating model inadequacy in ordinary differential equations with
  physics-informed neural networks.
\newblock {\em Comput. Struct.} {\bf 2021}, {\em 245},~106458.
\newblock
  {{https://doi.org/10.1016/j.compstruc.2020.106458}}.

\bibitem{yang2021bpinns}
{Yang, L.; Meng, X.; Karniadakis, G.}
\newblock {{B-PINNs: Bayesian physics-informed neural networks for forward and
  inverse PDE problems with noisy data}}.
\newblock {\em {J. Comput. Phys.}} {\bf {2021}}, {\em {425}},~{109913}
.
\newblock {{https://doi.org/10.1016/j.jcp.2020.109913}}.

\bibitem{yang2020physics}
Yang, L.; Zhang, D.; Karniadakis, G.
\newblock Physics-Informed Generative Adversarial Networks for Stochastic
  Differential Equations.
\newblock {\em SIAM J. Sci. Comput.} {\bf 2020}, {\em
  42},~A292--A317.
\newblock {{https://doi.org/10.1137/18M1225409}}.

\bibitem{yucesan2021hybrid}
Yucesan, Y.; Viana, F.
\newblock Hybrid physics-informed neural networks for main bearing fatigue
  prognosis with visual grease inspection.
\newblock {\em Comput. Ind.} {\bf 2021}, {\em 125},~103386.
\newblock
  {{https://doi.org/10.1016/j.compind.2020.103386}}.

\bibitem{zhu2019physics}
Zhu, Y.; Zabaras, N.; Koutsourelakis, P.; Perdikaris, P.
\newblock Physics-constrained deep learning for high-dimensional surrogate
  modeling and uncertainty quantification without labeled data.
\newblock {\em J. Comput. Phys.} {\bf 2019}, {\em 394},~56--81.
\newblock {{https://doi.org/10.1016/j.jcp.2019.05.024}}.

\bibitem{pang2019fpinns}
Pang, G.; Lu, L.; Karniadakis, G.E.
\newblock fPINNs: Fractional physics-informed neural networks.
\newblock {\em SIAM J. Sci. Comput.} {\bf 2019}, {\em
  41},~A2603--A2626.

\bibitem{liu2024kan}
Liu, Z.; Wang, Y.; Vaidya, S.; Ruehle, F.; Halverson, J.; Solja{\v{c}}i{\'c},
  M.; Hou, T.Y.; Tegmark, M.
\newblock Kan: Kolmogorov-arnold networks.
\newblock {\em arXiv} {\bf 2024}, arXiv:2404.19756.

\bibitem{koenig2024kan}
Koenig, B.C.; Kim, S.; Deng, S.
\newblock KAN-ODEs: Kolmogorov-Arnold Network Ordinary Differential Equations
  for Learning Dynamical Systems and Hidden Physics.
\newblock {\em arXiv} {\bf 2024}, arXiv:2407.04192.

\bibitem{qian2024investigating}
Qian, K.; Kheir, M.
\newblock {Investigating KAN-Based Physics-Informed Neural Networks for EMI/EMC
  Simulations}.
\newblock {\em arXiv} {\bf 2024}, arXiv:2405.11383.

\bibitem{kumar2023mycrunchgpt}
Kumar, V.; Gleyzer, L.; Kahana, A.; Shukla, K.; Karniadakis, G.E.
\newblock Mycrunchgpt: A llm assisted framework for scientific machine
  learning.
\newblock {\em J. Mach. Learn. Model. Comput.} {\bf
  2023}, {\em 4}, {41--72}.

\bibitem{beck2022overview}
Beck, C.; Hutzenthaler, M.; Jentzen, A.; Kuckuck, B.
\newblock An overview on deep learning-based approximation methods for partial
  differential equations.
\newblock {\em {Discret. Contin. Dyn. Syst. B}} {\bf 2022}, \emph{{28}}, {3697--3746}.
\newblock {{https://doi.org/10.3934/dcdsb.2022238}}.

\bibitem{cuomo2022scientific}
Cuomo, S.; Di~Cola, V.S.; Giampaolo, F.; Rozza, G.; Raissi, M.; Piccialli, F.
\newblock Scientific Machine Learning Through Physics-Informed Neural Networks:
  Where we are and What's Next.
\newblock {\em J. Sci. Comput.} {\bf 2022}, {\em 92}, {88}.
\newblock {{https://doi.org/10.1007/s10915-022-01939-z}}.

\bibitem{lawal2022physics}
Lawal, Z.; Yassin, H.; Lai, D.; Che~Idris, A.
\newblock {Physics-Informed Neural Network (PINN) Evolution and Beyond: A
  Systematic Literature Review and Bibliometric Analysis}.
\newblock {\em Big Data Cogn. Comput.} {\bf 2022}, {\em 6}, {140}.
\newblock {{https://doi.org/10.3390/bdcc6040140}}.

\bibitem{viana2021survey}
Viana, F.A.; Subramaniyan, A.K.
\newblock A survey of Bayesian calibration and physics-informed neural networks
  in scientific modeling.
\newblock {\em Arch. Comput. Methods Eng.} {\bf 2021},
  {\em 28},~3801--3830.

\bibitem{kharazmi2019variational}
Kharazmi, E.; Zhang, Z.; Karniadakis, G.
\newblock {VPINNs: Variational physics-informed neural networks for solving
  partial differential equations}.
\newblock {\em arXiv} {\bf 2019}, arXiv:1912.00873.

\bibitem{kharazmi2021hp}
Kharazmi, E.; Zhang, Z.; Karniadakis, G.
\newblock {$hp$-VPINNs: Variational physics-informed neural networks with
  domain decomposition}.
\newblock {\em Comput. Methods Appl. Mech. Eng.} {\bf
  2021}, {\em 374},~113547.

\bibitem{berrone2022solving}
Berrone, S.; Canuto, C.; Pintore, M.
\newblock {Solving PDEs by variational physics-informed neural networks: An a
  posteriori error analysis}.
\newblock {\em  Ann. Dell'Universita' Ferrara} {\bf 2022}, {\em
  68},~575--595.
\newblock {{https://doi.org/10.1007/s11565-022-00441-6}}.

\bibitem{berrone2022variational}
Berrone, S.; Canuto, C.; Pintore, M.
\newblock Variational-Physics-Informed Neural Networks: The role of quadratures
  and test functions.
\newblock {\em J. Sci. Comput.} {\bf 2022}, {\em 92},~{100}.

\bibitem{berrone2016towards}
Berrone, S.; Pieraccini, S.; Scial\`o, S.
\newblock Towards effective flow simulations in realistic discrete fracture
  networks.
\newblock {\em J. Comput. Phys.} {\bf 2016}, {\em
  310},~181--201.

\bibitem{sukumar2022exact}
Sukumar, N.; Srivastava, A.
\newblock Exact imposition of boundary conditions with distance functions in
  physics-informed deep neural networks.
\newblock {\em Comput. Methods Appl. Mech. Eng.} {\bf
  2022}, {\em 389},~114333.

\bibitem{berrone2022enforcing}
Berrone, S.; Canuto, C.; Pintore, M.; Sukumar, N.
\newblock Enforcing Dirichlet boundary conditions in physics-informed neural
  networks and variational physics-informed neural networks.
\newblock {\em Heliyon} {\bf 2023}, {\em 9},~e18820.

\bibitem{kingma2014adam}
Kingma, D.; Ba, J.
\newblock {Adam: A method for stochastic optimization}.
\newblock {\em arXiv} {\bf 2014}, {arXiv:1412.6980}
.

\bibitem{wright1999numerical}
Wright, S.; Nocedal, J.
\newblock {\em {Numerical Optimization}}; Springer: {Berlin/Heidelberg, Germany,} 1999; Volume 35, p. 7.

\bibitem{baydin2018automatic}
Baydin, A.; Pearlmutter, B.; Radul, A.; Siskind, J.
\newblock Automatic differentiation in machine learning: A survey.
\newblock {\em J. Mach. Learn. Res.} {\bf 2018}, {\em 18}, {5595--5637}.

\bibitem{prechelt1998early}
{Prechelt, L.
\newblock Early stopping-but when? In {\em Neural Networks: Tricks of the
  Trade}; Springer: {Berlin/Heidelberg, Germany,} 1998; pp. 55--69.}

\end{thebibliography}
\end{document}